\documentclass[11pt]{amsart}
\usepackage{amssymb,amsmath,txfonts,mathrsfs,mathtools, titletoc, tikz, stmaryrd} 
\usepackage[alphabetic]{amsrefs}

\usepackage{graphicx}

\newtheorem{theorem}{Theorem}[section]
\newtheorem*{theorem*}{Theorem}
\newtheorem{prop}[theorem]{Proposition}
\newtheorem{lemma}[theorem]{Lemma}
\newtheorem{remark}[theorem]{Remark}

\newtheorem{question}[theorem]{Question}
\newtheorem{definition}[theorem]{Definition}
\newtheorem{cor}[theorem]{Corollary}

\newtheorem{conj}[theorem]{Conjecture}

\numberwithin{equation}{section}

\def\pf{{\it Proof:}~}

\begin{document}

\title[Local estimate of fundamental groups]{Local estimate of fundamental groups}

\author{Guoyi Xu}
\address{Department of Mathematical Sciences,
\\Tsinghua University, Beijing\\P. R. China, 100084}
\email{guoyixu@tsinghua.edu.cn}
\date{\today}

\begin{abstract}
For any complete $n$-dim Riemannian manifold $M^n$ with nonnegative Ricci curvature, Kapovitch and Wilking proved that any finitely generated subgroup of the fundamental group $\pi_1(M^n)$ can be generated by $C(n)$ generators. Inspired by their work, we give a quantitative proof of the above theorem and show that $C(n)\leq n^{n^{20n}} $. Our main tools are quantitative Cheeger-Colding's almost splitting theory, and the squeeze lemma for covering groups between two Riemannian manifolds with nonnegative Ricci curvature.

%We studied the covering group between two complete Riemannian manifolds with nonnegative Ricci curvature, and estimate the number of generators for any finitely generated subgroup of the covering group. Our method is based on the existence of `partial' harmonic Gromov-Hausdorff approximation locally defined on the manifolds, which was established in Cheeger-Colding's study on almost rigidity result. As an application, we gave an alternative proof of Kapovitch and Wilking's theorem: any finitely generated subgroup of the fundamental group of complete Riemannian manifolds with nonnegative Ricci curvature, can be generated by uniformly finite generators. The original proof of Kapovitch and Wilking reduced to Cheeger-Colding's theory of Ricci limit space by contradiction. Our argument can be thought as the direct version of their indirect proof.
\end{abstract}

%\keywords{$p$-sweepout, Min-Max theory, $p$-width} \subjclass[2010]{53A10, 49Q05} 
\thanks{The author was partially supported by NSFC-11771230}

\maketitle

\titlecontents{section}[0em]{}{\hspace{.5em}}{}{\titlerule*[1pc]{.}\contentspage}
\titlecontents{subsection}[1.5em]{}{\hspace{.5em}}{}{\titlerule*[1pc]{.}\contentspage}
%\contentsline{section}{\numberline{1}something}{}
%\contentsline{subsection}{\numberline{}1.1 other}{}
\tableofcontents

\section*{Introduction}

It is well-known that any compact Riemannian manifold $(M^n, g)$ has finitely generated fundamental group. For non-compact complete Riemannian manifolds, the conclusion is not always right. For example, the surface with infinite genus has infinitely generated fundamental group. A natural question is:
\begin{question}
{For which complete Riemannian manifold $(M^n, g)$, the fundamental group $\pi_1(M^n)$ is finitely generated?
}
\end{question}

Note the above example has no non-negative sectional curvature metric, to obtain the finitely generatedness of fundamental group, we may consider adding some curvature assumption on complete Riemannian manifolds. 

For the group $\mathbf{G}$, we define that $\mathfrak{ng}(\mathbf{G})$ is the minimal number of generators needed. In $1911$ Bieberbach proved: For complete flat Riemannian manifold $M^n$, the fundamental group $\pi_1(M^n)$ is finitely generated and $\mathfrak{ng}\big(\pi_1(M^n)\big)\leq C(n)$.

Bieberbach reduced the study of fundamental groups of flat manifolds to the study of the discrete subgroup of the isometry group of $\mathbb{R}^n$, which was later developed into a more general theory about the discrete subgroups of Lie groups (see \cite{Rag}). 

Later for hyperbolic manifolds, the following theorem was obtained (see \cite{BGS}): 
Any finite volume hyperbolic manifold has finitely generated fundamental group.

\begin{remark}
{The above theorem is a corollary of the following general result in Lie group theory: Any lattice $\Gamma$ of a Lie group $\mathbf{G}$ is finitely generated, where lattice means that $\Gamma$ is a discrete subgroup of $\mathbf{G}$ and $\mathbf{G}/\Gamma$ has finite volume.
}
\end{remark}

The study of the above non-positive curvature case, has more algebraic flavor, which has strong contrast to the following non-negative curvature case.

In $1972$, Cheeger and Gromoll \cite{Soul} studied non-negative sectional curvature Riemannian manifolds. Among other things, they obtained the following result: If $M^n$ has sectional curvature $sec\geq 0$, then $\pi_1(M^n)$ is finitely generated.

In fact, in \cite{Soul} it was showed that $M^n$ is homotopic to a compact totally geodesic submanifold of $M^n$, through the deformation by the gradient flow of Busemann function. The above theorem follows as a corollary of this more general structure result.

In $1978$, Gromov \cite{Gromov-almost-flat} used Toponogov Comparison Theorem to study the fundamental group directly, and proved: For $M^n$ with sectional curvature $sec\geq 0$, $\mathfrak{ng}\big(\pi_1(M^n)\big)\leq C(n)= \frac{\mathrm{Vol}\big(\mathbb{S}^{n- 1}\big)}{\mathrm{Vol}\big(\mathbb{D}^{n- 1}(\frac{\pi}{6})\big)}$, where $\mathbb{S}^{n- 1}$ is the unit sphere in $\mathbb{R}^n$, $\mathbb{D}^{n- 1}(r)$ is the geodesic ball with radius $r$ in $\mathbb{S}^{n- 1}$.

For compact Riemannian manifolds with $Rc\geq 0$, under the additional assumption of conjugate radius, Guofang Wei \cite{Wei} gave a uniform estimate on the generators of the fundamental group similar as in Gromov's result.

Back in $1968$, Milnor \cite{Milnor} proved that for complete Riemannian manifold $M^n$ with $Rc\geq 0$, any finitely generated subgroup of $\pi_1(M^n)$, is polynomial growth of order $\leq n$. Furthermore, he posed the following conjecture:
\begin{conj}[Milnor]\label{conj Milnor}
{For complete Riemannian manifold $(M^n, g)$ with $Rc\geq 0$, $\pi_1(M^n)$ is finitely generated.
}
\end{conj}

In the $1980s$, Peter Li \cite{Li} used the heat kernel and Anderson \cite{Anderson} used the property of covering maps, to prove Milnor conjecture for Euclidean (maximal) volume growth case respectively (moreover, they proved the fundamental group is finite in fact).

In $2000$, Sormani \cite{Sormani} used the excess estimate of Abresch and Gromoll \cite{AG} on the universal cover of manifolds, successfully proved Milnor conjecture for linear (minimal) volume growth case.

Also in $2000$, B. Wilking \cite{Wilking} used Milnor's Theorem and the theory of discrete subgroup in Lie group to prove: If for any complete Riemannian manifold $(M^n, g)$ with $Rc(g)\geq 0$ and $\pi_1(M^n)$ is an abelian group, we have $\pi_1(M^n)$ is finitely generated; then for any complete Riemannian manifold $(N^n, \tilde{g})$ with $Rc(\tilde{g})\geq 0$, $\pi_1(N^n)$ is finitely generated. 

Note the proof of the above Wilking's Theorem do not need to use Bochner formula in Riemannian geometry, only relies on the Bishop-Gromov Volume Comparison Theorem. 

Recently, in $2011$, V. Kapovitch and B. Wilking \cite{KW} proved the following local estimate of fundamental groups among other things :
\begin{theorem}[Kapovitch and Wilking]\label{thm KW}
{For complete Riemannian manifold $(M^n, g)$ with $Rc\geq 0$, any finitely generated subgroup $\Gamma$ of $\pi_1(M^n)$ satisfies $\mathfrak{ng}(\Gamma)\leq C(n)$.
}
\end{theorem}

Their proof was inspired by Fukaya and Yamaguchi's work \cite{FY}, started by contradiction, used the equivariant pointed Gromov-Hausdorff convergence and reduced the problem to the study of Ricci limit space with group actions, then do the induction on the dimension of Ricci limit space. The main technical tools are Cheeger-Colding's theory of Ricci limit space and almost rigidity results. However this proof can not give the explicit estimate of the above $C(n)$.

One main purpose of this paper is to give an explicit uniform estimate of $\mathfrak{ng}(\Gamma)$ for finitely generated subgroup $\Gamma\subseteq  \pi_1(M^n)$. More precisely, we prove the following theorem, which can be thought as the quantitative version of Theorem \ref{thm KW}.
\begin{theorem}\label{thm main of xu}
{For complete Riemannian manifold $(M^n, g)$ with $Rc\geq 0$, any finitely generated subgroup $\Gamma$ of $\pi_1(M^n)$ satisfies $\mathfrak{ng}(\Gamma)\leq n^{n^{20n}} $.
}
\end{theorem}

\begin{remark}
{The above upper bound $n^{n^{20n}} $ is not sharp, and our method can not provide the sharp bound either. Also, comparing the concise proof of Theorem \ref{thm KW} in \cite{KW}, our proof is sort of lengthy. 

One advantage of our proof is self-contained. Basically, our argument only use Bishop-Gromov Volume Comparison Theorem and the Bochner formula for complete Riemannian manifolds with $Rc\geq 0$. We try to reveal the relation between discrete isometry group actions and $Rc\geq 0$ in an intrinsic way. Even the Compactness Theorem of Gromov-Hausdorff convergence is not used, let alone the theory of Ricci limit spaces, we only need the concept of $\epsilon$-Gromov-Hausdorff approximation. 
}
\end{remark}

There are three key ingredients of our proof. The first one tells us how to transfer from Gromov-Hausdorff approximation (geometry assumption) to almost orthonormal linear harmonic functions (analysis result). The second one is doing reverse argument, i.e. transferring analysis to geometry. And the third one studies the generating set of the covering groups through analysis and geometry. We will describe those three key ingredients in the rest of this section, and conclude the section with a sketchy description of our proof.

Firstly, we hope to reveal more explicit relation between group actions and the assumption $Rc\geq 0$, through the concrete analysis on distance function, in the similar spirit of Gromov's proof about estimate of fundamental group's generators. This hope starts from the following classical global result of Cheeger and Gromoll \cite{Cheeger-Gromoll} obtained in $1971$:
\begin{theorem}[Cheeger and Gromoll]\label{thm Cheeger and Gromoll Ricci}
{If $(M^n, g)$ is a complete $n$-dim Riemannian manifold with $Rc\geq 0$, and $M^n$ contains a line, then $M^n$ is isometric to $N^{n- 1}\times \mathbb{R}$, where $N^{n- 1}$ is a complete $(n- 1)$-dim Riemannian manifold with $Rc\geq 0$.
}
\end{theorem}

The proof of Theorem \ref{thm Cheeger and Gromoll Ricci} used the harmonic function $\mathbf{b}$ globally defined on $M^n$, which is sort of the limit of the distance function of $M^n$ (more precisely, the Busemann function). In fact, the above manifold $N^{n- 1}$ is a level set of $\mathbf{b}$ and the splitting lines $\mathbb{R}$ are the gradient flow lines of $\mathbf{b}$. One crucial technical point is to get the gradient estimate of the harmonic function $\mathbf{b}$, which is $|\nabla \mathbf{b}|\equiv 1$.

In $1975$, Cheng and Yau \cite{CY} established the well-known local gradient estimate of harmonic functions on complete Riemannian manifolds with $Rc\geq 0$, which enables us to obtain the gradient estimate of harmonic functions from the $C^0$-bound of harmonic functions locally. On such manifolds, in $1990$, Abresch and Gromoll \cite{AG} obtained the excess estimate, which gives the local estimate of the sum of two local Busemann functions with respect to a segment (such sum is $0$ for two Busemann  functions with reverse directions in the proof of Theorem \ref{thm Cheeger and Gromoll Ricci}).

Then in the $1990$s, based on the gradient estimate and the excess estimate, Colding initiated the study of local properties of distance function, through harmonic function locally defined on manifolds with Ricci lower bound, in a series paper \cite{Colding-shape}, \cite{Colding-large}, \cite{Colding-volume}; while solving several important problems in the theory of Gromov-Hausdorff convergence, which was established by Gromov in the $1980$s (for more details see \cite{Gromov-book}).  More precisely, among other things, Colding constructed the locally defined harmonic function based on the local Busemann function with respect to a segment, and proved such harmonic functions are `almost linear' harmonic functions with bounded gradient, where the `almost linear' is in average integral sense. 

The \textbf{first} ingredient of our proof is the above existence of almost orthonormal linear harmonic functions, which was originally established by Colding in \cite{Colding-shape}, \cite{Colding-large}, \cite{Colding-volume}. For our purpose, we need the explicit quantitative version, so we give all the details of the proof here. Although some calculation is sort of tedious, these explicit estimates possibly give an intrinsic expression of the transfer from geometry (Gromov-Hausdorff approximation) to analysis (existence of almost orthonormal linear harmonic functions).

In $1996$, Cheeger and Colding \cite{CC-Ann} established the almost splitting theorem among other things, which transfers from analysis (existence of almost orthonormal linear harmonic functions) to geometry (Gromov-Hausdorff approximation). Very roughly, they proved that if there exist $k$ almost orthonormal linear harmonic functions on a geodesic ball $B_r(p)\subset M^n$ with $Rc\geq 0$, then a smaller concentric geodesic ball is close to a ball of $\mathbb{R}^k\times \mathbf{X}_k$ in the sense of Gromov-Hausdorff distance. 

This almost splitting theorem is the \textbf{second} ingredient of our proof. For the same reason as the above, we need the quantitative version. During the proof of the almost splitting theorem, one crucial thing is the existence of a suitable Gromov-Hausdorff approximation under the suitable Hessian integral bound assumption. Because the original proof of this existence result in \cite{CC-Ann} is concise, one of our contribution is a different detailed proof by modifying some argument in \cite{CN}, for more details see Section \ref{SEC splitting on Ricci limit spaces} of this paper.

It is a natural question whether the topology of two metric spaces are the same when they are very close in Gromov-Hausdorff distance sense. Generally, the answer is no, although there is an intrinsic Reifenberg type theorem when one of the spaces is $n$-dim Euclidean space and the other space is $n$-dim Riemannian manifolds (for details, see \cite[Appendix]{CC1}).

However, for a family of Riemannian manifolds converging to a metric space in Gromov-Hausdorff distance sense, using the equivariant Gromov-Hausdorff convergence theory developed by Fukaya and Yamaguchi in the $1980$s (see \cite{Fukaya-Japan}, \cite{FY}), Kapovitch and Wilking \cite{KW} obtained the results, which relate the fundamental groups of those converging Riemannian manifolds, to the limit group, which acts on the limit space of the universal covers of those Riemannian manifolds. 

The study of the fundamental groups can be put into a more general context, i.e. the covering group of a covering map between two Riemannian manifolds with $Rc\geq 0$. The \textbf{third} ingredient of our proof is to characterize the change of the covering groups by the Gromov-Hausdorff distance between two metric spaces in quantitative form, where one metric space is a geodesic ball in Riemannian manifolds with $Rc\geq 0$ and the other one is a ball in a product metric space $\mathbb{R}^k\times \mathbf{X}_k$. Our squeeze lemma provides a bridge linking analysis with group actions and geometry. Very roughly, if there are $k$ almost orthonormal linear harmonic functions on a geodesic ball $B_r(p)\subset M^n$ with $Rc\geq 0$ (analysis), then from the almost splitting theorem in the second ingredient, we know that $B_{cr}(p)$ is close to $B_{cr}(0, \hat{p})\subset \mathbb{R}^k\times \mathbf{X}_k$ in Gromov-Hausdorff sense (geometry). Then the group actions on $B_{cr}(p)$ are almost generated by the group actions on $B_{cr}(\hat{p})\subset \mathbf{X}_k$, for more details see Lemma \ref{lem first squeeze big Eucliean directions}.

Now we describe our proof in a rough way. We start with a geodesic ball $B_r(p)\subset M^n$ with $Rc\geq 0$, firstly we use the first ingredient tool to find one almost linear harmonic function, then apply Squeeze Lemma to shrink the group action to a group action on a ball $B_{cr}(0, \hat{p})\subset \mathbf{X}_{n- 1}$. Now we apply the first ingredient tool again to find two almost orthonormal linear harmonic functions on a smaller geodesic ball $B_{r_1}(p_1)\subset M^n$, and the group action on $B_r(p)$ can be `controlled' or generated by the group action on $B_{r_1}(p_1)$. We repeat the above procedure by induction on the dimension of almost orthonormal linear harmonic functions. When the dimension is $n$, the group action is shown to be trivial, and we get our conclusion.  

%If we have `almost linear' harmonic function locally defined on manifolds, we have the squeeze lemma to squeeze the group action on the original manifolds, almost to the group action on the lower dimension metric space $\mathbf{X}$ in the almost splitting theorem. This squeeze lemma is one of the key technical tools of our argument, and we squeeze the group action gradually, by induction on the dimension of `almost orthonormal linear' harmonic functions locally defined on gradually shrinking domain of the manifolds. When the dimension of such locally defined harmonic functions is $n$, we showed the group action is almost trivial directly, without applying the intrinsic Reifenberg theorem mentioned above. 

\section*{Part I. G-H approximation yields A.O.L. harmonic functions}

In Part I of this paper, we will prove the following version result about analytic characterization of Gromov-Hausdorff approximation, which was implied in the argument of Cheeger and Colding in a series of papers, \cite{Colding-shape}, \cite{Colding-large}, \cite{Colding-volume}, \cite{CC-Ann} and \cite{Cheeger-GAFA}.

Let us recall the definition of the pointed $\epsilon$-Gromov-Hausdorff approximation.
\begin{definition}\label{def pted G-H appro}
{Let $(\mathbf{X}, d_\mathbf{X}, x_0)$ and $(\mathbf{Y}, d_\mathbf{Y}, y_0)$ be two pointed metric spaces. For $\epsilon> 0$, a map $f: \big(\mathbf{X}, x_0\big)\rightarrow (\mathbf{Y}, y_0)$ is called a \textbf{pointed $\epsilon$-Gromov-Hausdorff approximation} if 
\begin{align}
&f(x_0)= y_0\ , \quad \quad \quad \quad  \mathbf{Y}\subset \mathbf{U}_{\epsilon}\big(f(\mathbf{X})\big),\nonumber \\
&\Big|d_\mathbf{Y}\big(f(x_1), f(x_2)\big)- d_\mathbf{X}(x_1, x_2)\Big|< \epsilon \ , \quad \quad \quad \quad \forall x_1, x_2\in \mathbf{X}, \nonumber 
\end{align}
where $\mathbf{U}_{\epsilon}\big(f(\mathbf{X})\big)\vcentcolon = \Big\{z\in \mathbf{Y}: d\big(z, f(\mathbf{X})\big)\leq \epsilon\Big\}$. For simplicity, we also call that $f$ is an $\epsilon$-Gromov-Hausdorff approximation when the base points are fixed and clear.
}
\end{definition}

In the rest of the paper, unless otherwise mentioned, we use $\mathbf{X}, \mathbf{Y}, \mathbf{X}_k$ to denote metric spaces. Also, we are only interested in geometry and analysis on $n$-dim manifolds with $n\geq 3$, so we will always assume the dimension of any manifolds $\mathbf{n\geq 3}$ in the rest of this paper.

For our application, we need the quantitative estimate, which relate the existence of the Gromov-Hausdorff approximation to the existence of almost orthonormal linear harmonic functions. Although many results in Part I are well-known to some experts in this field, we made the contribution to establish the suitable statement and the quantitative estimates. Also we elaborate the concise argument of Cheeger-Colding to present the proof of some known results in all the details for self-contained reason, and also hope to provide a backup reference for future study, besides the original works of Cheeger-Colding.

\begin{theorem}\label{thm AG-triple imply one more splitting-1}
{For $B_{10r}(q)\subset (M^n, g)$ with $Rc(g)\geq 0$, any $0< \epsilon< 1, \delta= n^{-3400n^3}\epsilon^{110n}$ and integer $0\leq k\leq n$, assume there is an $(\delta r)$-Gromov-Hausdorff approximation, $f: B_{10r}(q)\rightarrow B_{10r}(0, \hat{q})\subset \mathbb{R}^k\times \mathbf{X}_k$, and $\mathrm{diam}\big(B_r(\hat{q})\big)\geq \frac{1}{4}r$ where $B_r(\hat{q})\subset \mathbf{X}_k$. Then there are harmonic functions $\big\{\mathbf{b}_i\big\}_{i= 1}^{k+ 1}$ defined on $B_{s}(p)\subset B_{10r}(q)$, where $s= n^{-320n^3}\epsilon^{10n}r$, such that 
\begin{align}
\sup_{B_{s}(p)\atop i= 1, \cdots, k+ 1} |\nabla \mathbf{b}_i|\leq 1+ \epsilon \quad \quad and \quad \quad \sup_{t\leq s}\fint_{B_{t}(p)}\sum_{i, j= 1}^{k+ 1} \big|\langle \nabla \mathbf{b}_i, \nabla \mathbf{b}_j\rangle- \delta_{ij}\big|\leq \epsilon . \nonumber 
\end{align}
}
\end{theorem}

%\section{Colding's integral Toponogov theorem}

\section{Almost orthonormal local Busemann functions}

Cheng-Yau's gradient estimates was originally proved in \cite{CY}, the following form include an explicit form of constant, which is needed in the later proof. The proof is the same as in \cite{CY}, so we omit it.

\begin{theorem}[Cheng-Yau's gradient estimates]\label{thm Cheng-Yau's lemma}
{Assume $Rc(M^n)\geq 0, p\in M^n$, $B_{2R}(p)\subset M^n$, 
$f: B_{2R}(p)\rightarrow \mathbb{R}$ is a harmonic function and $f\in C\big(\overline{B_{2R}(p)}\big)$, then 
\begin{align}
\sup_{B_R(p)} |\nabla f|\leq \frac{60n}{R} \sup_{B_{2R}(p)} |f| .\nonumber 
\end{align}
If $\Delta f= c_0\geq 0$ for $f\in C\big(\overline{B_{2R}(p)}\big)$, then $\sup\limits_{B_R(p)} |\nabla f|\leq \frac{200n(R+ 1)}{R} \Big[\sup_{B_{2R}(p)} |f|+ c_0\Big]$.
}
\end{theorem}\qed

Set $b^{+}(\cdot)= d(q, \cdot)- d(q, p): M^n\rightarrow \mathbb{R}$, the function $b^+$ is called the \textbf{local Busemann function} with respect to the couple points $[p, q]$. And we define the positive part of a function $f$ as 
\begin{align}
f_+(x)= \max\{f(x), 0\} .\nonumber 
\end{align}

\begin{lemma}\label{lem estimate of Lap b}
{Given $R> 0$, $Rc(M^n)\geq 0$ and $p, q\in M$ with $d(p, q)> 2R$, then
\begin{align}
\frac{1}{V\big(B_R(p)\big)}\int_{B_R(p)} |\Delta b^{+}|\leq \frac{3n}{R} , \nonumber 
\end{align}
where $b^{+}(\cdot)= d(q, \cdot)- d(q, p): M^n\rightarrow \mathbb{R}$.
}
\end{lemma}

\pf
{Note for any $x\in B_R(p)$, 
\begin{align}
d(q, x)\geq d(q, p)- d(p, x)> 2R- R= R .\nonumber 
\end{align}
By Laplace Comparison Theorem, 
\begin{align}
\Delta b^{+}= \Delta d(q, \cdot)\leq \frac{n- 1}{d(q, \cdot)}< \frac{n- 1}{R} \ ,\quad \quad \quad \quad on \ B_R(p) \nonumber 
\end{align}
hence the positive part of $\Delta b^+$ satisfies
\begin{align}
\sup_{x\in B_R(p)}\big(\Delta b^+\big)_{+}(x)\leq \frac{n- 1}{R} . \label{Delta posi part}
\end{align}

By Bishop-Gromov Comparison Theorem, 
\begin{align}
\frac{V\big(\partial B_R(p)\big)}{V\big(B_R(p)\big)}\leq \frac{n}{R} .\label{eq 1.1.2}
\end{align}

Now from (\ref{Delta posi part}), (\ref{eq 1.1.2}), we have
\begin{align}
\int_{B_R(p)} |\Delta b^+| &= 2\int_{B_R(p)} \big(\Delta b^+\big)_+- \int_{B_R(p)} \Delta b^+  \nonumber \\
&\leq 2V\big(B_R(p)\big)\max_{B_R(p)} \big(\Delta b^+\big)_+ + \Big|\int_{B_R(p)} \Delta b^+\Big| \nonumber \\
&\leq \frac{2(n- 1)}{R}V\big(B_R(p)\big)+ \Big|\int_{\partial B_R(p)} \frac{\partial b^+}{\partial \vec{n}}\Big| \nonumber \\
&\leq \frac{2(n- 1)}{R}V\big(B_R(p)\big)+ V\big(\partial B_R(p)\big) \nonumber \\
&\leq \frac{3n}{R}V\big(B_R(p)\big) .\nonumber 
\end{align}
}
\qed

\begin{lemma}\label{lem harmonic b}
{Suppose $Rc(M^n)\geq 0$, $x\in M$, and the function $\mathbf{b}$ satisfies 
\begin{equation}\nonumber 
\left\{
\begin{array}{rl}
\Delta \mathbf{b}&= 0 \quad \quad \quad \quad \quad on\ B_{4R}(x) \\
\mathbf{b}&= b^+ \quad \quad \quad \quad \quad on\ \partial B_{4R}(x) ,
\end{array} \right.
\end{equation}
where $b^{+}(\cdot)= d(q, \cdot)- d(q, p): M^n\rightarrow \mathbb{R}$. Then 
\begin{align}
\fint_{B_{4R}(x)} |\nabla (\mathbf{b}- b^+)|^2 &\leq 8R\cdot \fint_{B_{4R}(x)} |\Delta b^+| , \label{eq 1.11} \\
\fint_{B_{2R}(x)} |\nabla^2 \mathbf{b}|^2 &\leq \frac{10^8\cdot n^3}{R^2}\Big[1+ \frac{d(p, x)}{R}\Big]^2 . \label{eq 1.13}
\end{align}
}
\end{lemma}

\pf
{From Maximum principle, one have $z_1, z_2\in \partial B_{4R}(x)$ such that 
\begin{align}
\mathbf{b}(z_1)= \sup_{z\in B_{4R}(x)} \mathbf{b}(z)\ , \quad \quad \quad \quad \mathbf{b}(z_2)= \min_{z\in B_{4R}(x)}\mathbf{b}(z) . \nonumber 
\end{align}
Then for any $y\in B_{4R}(x)$, 
\begin{align}
\mathbf{b}(y)- b^+(y)&\leq \mathbf{b}(z_1)- b^+(y)= b^+(z_1)- b^+(y)= d(q, z_1)- d(q, y) \nonumber  \\
&\leq d(z_1, y)\leq d(z_1, x)+ d(x, y)\leq 8R . \nonumber 
\end{align}
Similarly, we have $\mathbf{b}(y)- b^+(y)\geq -d(z_2, y)\geq -8R$. Hence 
\begin{align}
\sup_{B_{4R}(x)} |\mathbf{b}- b^+|\leq 8R . \label{eqa 1.2}
\end{align}

From integration by parts, we get 
\begin{align}
\int_{B_{4R}(x)} \big|\nabla (\mathbf{b}- b^+)\big|^2&= -\int_{B_{4R}(x)} (\mathbf{b}- b^+)\Delta (\mathbf{b}- b^+)= \int_{B_{4R}(x)} (\mathbf{b}- b^+)\Delta b^+ \nonumber \\
&\leq  \sup_{B_{4R}(x)} |\mathbf{b}- b^+| \int_{B_{4R}(x)} |\Delta b^+| \leq 8R \int_{B_{4R}(x)} |\Delta b^+| .\nonumber 
\end{align}

From $\Delta \mathbf{b}= 0$ and Bochner's formula, 
\begin{align}
\frac{1}{2}\Delta\big(|\nabla \mathbf{b}|^2\big)= |\nabla^2 \mathbf{b}|^2+ \langle \nabla \Delta \mathbf{b}, \nabla \mathbf{b}\rangle+ Rc(\nabla \mathbf{b}, \nabla \mathbf{b}) \geq |\nabla^2 \mathbf{b}|^2 . \label{eq 1.2.-2}
\end{align}

Let $\phi\in C^{\infty}(R_+)$ be a nonnegative cut-off function such that 
\begin{equation}\nonumber 
\phi(s)= \left\{
\begin{array}{rl}
&1 \quad \quad \quad \quad \quad 0\leq s\leq 2R \\
&0 \quad \quad \quad \quad \quad s\geq 3R  ,
\end{array} \right.
\end{equation}
and 
\begin{align}
-\frac{2}{R}\leq \phi '\leq 0\ , \quad \quad \quad \quad |\phi ''|\leq \frac{16}{R^2} . \label{eq 1.2.-1}
\end{align}

We use the notation $d_p(x)= d(p, x)$, similarly $d_x(\cdot)= d(x, \cdot)$. From Theorem \ref{thm Cheng-Yau's lemma} and Maximum principle,
\begin{align}
\sup_{B_{3R}(x)} |\nabla \mathbf{b}|\leq \frac{60n}{R} \sup_{B_{4R}(x)} |\mathbf{b}|\leq \frac{60n}{R}\sup_{\partial B_{4R}(x)} |b^+|\leq 60n\cdot \Big[4+ \frac{d_p(x)}{R}\Big] .\label{eq 1.2.-3} 
\end{align}
Then from (\ref{eq 1.2.-2}) and (\ref{eq 1.2.-3}),
\begin{align}
\int_{B_{2R}(x)} |\nabla^2 \mathbf{b}|^2&\leq \int_{B_{3R}(x)} |\nabla^2 \mathbf{b}|^2\cdot (\phi\circ d_x)\leq \int_{B_{3R}(x)} \frac{1}{2} \Delta \big(|\nabla \mathbf{b}|^2\big) \cdot (\phi\circ d_x) \nonumber \\
&=  \frac{1}{2} \int_{B_{3R}(x)}|\nabla \mathbf{b}|^2\cdot \Delta \big(\phi\circ d_x\big) \leq 10^6n^2\Big[1+ \frac{d_p(x)}{R}\Big]^2 \int_{B_{3R}(x)} \Big|\Delta \big(\phi\circ d_x\big)\Big|. \label{eq 1.2.0}
\end{align}

On the other hand, from (\ref{eq 1.2.-1}) and Laplace Comparison Theorem, 
\begin{align}
\int_{B_{3R}(x)} \big|\Delta (\phi\circ d_x)\big|&= 2\int_{B_{3R}(x)} \big(\Delta (-\phi\circ d_x)\big)_{+}- \int_{B_{3R}(x)} \Delta (-\phi\circ d_x) \nonumber \\
&\leq 2V\big(B_{3R}(x)\big)\max_{B_{3R}(x)- B_{2R}(x)} \big(\Delta (-\phi\circ d_x)\big)_{+}+ \Big|\int_{B_{3R}(x)} \Delta (\phi\circ d_x)\Big| \nonumber \\
&\leq 2V\big(B_{3R}(x)\big)\max_{B_{3R}(x)- B_{2R}(x)} \big(-\phi '\Delta d_x- \phi ''\big)_{+}+ \int_{\partial B_{3R}(x)} \Big|\frac{\partial (\phi\circ d_x)}{\partial \vec{n}}\Big| \nonumber \\
&\leq 2V\big(B_{3R}(x)\big)\max_{B_{3R}(x)- B_{2R}(x)} \Big[\frac{2}{R}\cdot \frac{n- 1}{d_x}+ |\phi ''|\Big]+ \sup_{\partial B_{3R}(x)}|\phi '|\cdot V\big(\partial B_{3R}(x)\big) \nonumber \\
&\leq 2V\big(B_{3R}(x)\big)\cdot \frac{20n}{R^2}+ \frac{2n}{R^2}V\big(B_{3R}(x)\big) \nonumber \\
&\leq \frac{50n}{R^2}V\big(B_{2R}(x)\big) .\label{eq 1.2.2}
\end{align}

The inequality (\ref{eq 1.13}) follows from (\ref{eq 1.2.0}) and (\ref{eq 1.2.2}). 
}
\qed

\begin{lemma}\label{lem good choice of balls}
{Let $0< \delta< \frac{1}{2}$, $R> 0$, $D> 1$ be given, suppose $Rc(M^n)\geq 0$, $\tilde{p}\in M$, and $f, h\in L^1\big(B_R(\tilde{p})\big)$ with $\fint_{B_R(\tilde{p})} |f| \leq K$ , $\fint_{B_R(\tilde{p})} |h| \leq \hat{K}$. Then there exists finite many disjoint balls $\{B_{\frac{1}{2}\delta R}(x_j)\}_{j\in \mathbf{J}}$ such that we have the following:
\begin{enumerate}
\item[(I)] For any $j\in \mathbf{J}$, $B_{\frac{1}{2}\delta R}(x_j)\subset B_{(1- \frac{3}{2}\delta)R}(\tilde{p})$ and $B_{2\delta R}(x_j)\subset B_R(\tilde{p})$.
\item[(II)] 
\begin{equation}\nonumber 
\left\{
\begin{array}{rl}
V\big(\cup_{i\in \mathbf{J}}B_{\delta R}(x_i)\big)&\geq \big[(1- 2\delta)^n- \frac{2\cdot 10^n}{D}\big]V\big(B_R(\tilde{p})\big)  \\
\sup\limits_{i\in \mathbf{J}}\fint_{B_{2\delta R}(x_i)} |f| &\leq DK\quad \quad and \quad \quad \sup\limits_{i\in \mathbf{J}}\fint_{B_{2\delta R}(x_i)} |h| \leq D\hat{K} .
\end{array} \right.
\end{equation}
\end{enumerate}
}
\end{lemma}

\pf
{Choose a maximal set of disjoint balls $\{B_{\frac{1}{2}\delta R}(x_i)\}_{i\in \mathbf{S}}$ contained in $B_{(1- \frac{3}{2}\delta)R} (\tilde{p})$ and $B_{2\delta R}(x_i)\subset B_R(\tilde{p})$. Assume $\tilde{q}\in \cap_{i= 1}^{\lambda} B_{2\delta R}(x_i)$, where $\{1, \cdots, \lambda\}\subset \mathbf{S}$, and  note $B_{\frac{1}{2}\delta R}(x_i)\subset B_{\frac{5}{2}\delta R}(\tilde{q})$, 
\begin{align}
V\big(B_{\frac{5}{2}\delta R}(\tilde{q})\big)&\geq \sum_{i= 1}^{\lambda} V\big(B_{\frac{1}{2}\delta R}(x_i)\big)\geq \sum_{i= 1}^{\lambda} \Big(\frac{\frac{1}{2}\delta R}{5\delta R}\Big)^n V\big(B_{5\delta R}(x_i)\big) \nonumber \\
&\geq \sum_{i= 1}^{\lambda} \Big[ 10^{-n} V\big(B_{\frac{5}{2}\delta R}(\tilde{q})\big) \Big]= \lambda \cdot 10^{- n} V\big(B_{\frac{5}{2}\delta R}(\tilde{q})\big) . \nonumber 
\end{align}
Hence $\lambda\leq 10^n$, and for all $\tilde{q}\in B_R(\tilde{p})$, there exists at most $10^n$-many $i$ with $\tilde{q}\in B_{2\delta R}(x_i)$.

Note $B_{(1- 2\delta)R}(\tilde{p})\subset \cup_{i\in \mathbf{S}} B_{\delta R}(x_i)$ (otherwise, if $y\in B_{(1- 2\delta)R}(\tilde{p})- \cup_{i\in \mathbf{S}} B_{\delta R}(x_i)$, then $B_{\frac{1}{2}\delta R}(y)\cap B_{\frac{1}{2}\delta R}(x_i)= \emptyset$ for any $i\in \mathbf{S}$ and $B_{\frac{1}{2}\delta R}(y)\subset B_{(1- \frac{3}{2}\delta) R}(\tilde{p})$, contradicting the choice of $\big\{B_{\frac{1}{2}\delta R}(x_i)\big\}_{i\in \mathbf{S}}$).

Set 
\begin{align}
\mathbf{I}= \Big\{i\in \mathbf{S}\big|\ \int_{B_{2\delta R}(x_i)} |f|\geq DK\cdot V\big(B_{2\delta R}(x_i)\big) \ or \ \int_{B_{2\delta R}(x_i)} |h|\geq D\hat{K}\cdot V\big(B_{2\delta R}(x_i)\big)\Big\} .\nonumber 
\end{align}

Let $\mathbf{J}= \mathbf{S}- \mathbf{I}$, then 
\begin{align}
V\big(\cup_{j\in \mathbf{J}} B_{\delta R}(x_j)\big)&\geq V\big(\cup_{i\in \mathbf{S}} B_{\delta R}(x_i)\big)- V\big(\cup_{i\in \mathbf{I}} B_{\delta R}(x_i)\big) \nonumber \\
&\geq V\big(B_{(1- 2\delta)R}(\tilde{p})\big)- V\big(\cup_{i\in \mathbf{I}} B_{\delta R}(x_i)\big) .\label{eq 1.19.1}
\end{align}

Now we have 
\begin{align}
V\big(\cup_{i\in \mathbf{I}} B_{\delta R}(x_i)\big)&\leq \sum_{i\in \mathbf{I}} V\big(B_{2\delta R}(x_i)\big)\leq \frac{1}{DK}\sum_{i\in \mathbf{I}} \int_{B_{2\delta R}(x_i)} |f|+ \frac{1}{D\hat{K}}\sum_{i\in \mathbf{I}} \int_{B_{2\delta R}(x_i)} |h| \nonumber \\
&\leq 10^n\cdot \Big[\frac{1}{DK}\int_{B_R(\tilde{p})} |f|+ \frac{1}{D\hat{K}}\int_{B_R(\tilde{p})} |h|\Big]\leq \frac{2\cdot 10^n}{D}V\big(B_R(\tilde{p})\big) .\label{eq 1.19.2}
\end{align}

From (\ref{eq 1.19.1}) and (\ref{eq 1.19.2}), we have 
\begin{align}
V\big(\cup_{j\in \mathbf{J}} B_{\delta R}(x_j)\big)&\geq \big[(1- 2\delta)^n- \frac{2\cdot 10^n}{D}\big] V\big(B_R(\tilde{p})\big) \nonumber \\
\fint_{B_{2\delta R}(x_i)} |f| &\leq DK\quad and \quad \fint_{B_{2\delta R}(x_i)} |h| \leq D\hat{K} \ , \quad \quad \quad \quad \forall i\in \mathbf{J} . \nonumber 
\end{align}

Then $\big\{B_{\frac{1}{2}\delta R}(x_j)\big\}_{j\in \mathbf{J}}$ is our choice satisfying all properties required.
}
\qed

\begin{lemma}\label{lem appr by harm func}
{Suppose $Rc(M^n)\geq 0$, and $p, q\in M$ with $d(p, q)> 2r$. For all $\epsilon> 0$, there exists finitely many balls $B_{\delta r}(x_i)\subset B_r(p), \forall i\in \mathbf{A}$, where $\mathbf{A}$ is a finite set, and harmonic functions $\mathbf{b}_i$ with $\delta= 2^{-150n}\epsilon^{18}$ such that 
\begin{align}
&V\big(\bigcup_{i\in \mathbf{A}} B_{\delta r}(x_i)\big)\geq (1- \epsilon)V\big(B_r(p)\big)\ , \quad \quad  \quad \sum_{i\in \mathbf{A}}V\big(B_{\delta r}(x_i)\big)\leq 2^{2n} V\big(B_r(p)\big) , \nonumber \\
&\fint_{B_{2\delta r}(x_i)}|\nabla (\mathbf{b}_i- b^+)| \leq \epsilon \ , \quad \quad  \quad
\fint_{B_{2\delta r}(x_i)} |\nabla^2 \mathbf{b}_i|\leq \frac{\epsilon}{\delta r} ,\nonumber 
\end{align}
where $b^{+}(\cdot)= d(q, \cdot)- d(q, p): M^n\rightarrow \mathbb{R}$.
}
\end{lemma}

\pf
{By scaling the metric $g$ to the new metric $r^{-2}g$, we only need to prove the conclusion for $r= 1$.

From Lemma \ref{lem estimate of Lap b}, we have $\fint_{B_1(p)} |\Delta b^+|\leq 3n$. Combining this with Lemma \ref{lem good choice of balls}, for $D> 1$, and $\frac{1}{2}> \delta_1> 0$ (to be determined later), we can choose finitely many balls $\Big\{B_{\delta_1}(y_j)\Big\}_{j\in \mathbf{J}}$ and $B_{2\delta_1}(y_j)\subset B_1(p)$ such that 
\begin{equation}\nonumber 
\left\{
\begin{array}{rl}
B_{\frac{1}{2}\delta_1}(y_j)\cap B_{\frac{1}{2}\delta_1}(y_i)&= \emptyset \ , \quad \quad \quad \quad i\neq j \\
V\big(\cup_{j\in \mathbf{J}}B_{\delta_1}(y_j)\big)&\geq \big[(1- 2\delta_1)^n- \frac{2\cdot 10^n}{D}\big]V\big(B_1(p)\big)  \\
\fint_{B_{2\delta_1}(y_j)} |\Delta b^+| &\leq D\cdot (3n) \ , \quad \quad \quad \quad \forall j\in \mathbf{J}  .
\end{array} \right.
\end{equation}

For each $j\in \mathbf{J}$, let $\mathbf{b}_j$ satisfy 
\begin{equation}\nonumber 
\left\{
\begin{array}{rl}
\Delta \mathbf{b}_j&= 0 \ , \quad \quad \quad \quad on\ B_{2\delta_1}(y_j) \\
\mathbf{b}_j&= b^+ \ , \quad \quad \quad \quad on \ \partial B_{2\delta_1}(y_j)  . 
\end{array} \right.
\end{equation}

Apply Lemma \ref{lem harmonic b} to $B_{2\delta_1}(y_j)$, we have 
\begin{align}
\fint_{B_{\delta_1}(y_j)} |\nabla (\mathbf{b}_j- b^+)|^2 &\leq \frac{V\big(B_{2\delta_1}(y_j)\big)}{V\big(B_{\delta_1}(y_j)\big)} \fint_{B_{2\delta_1}(y_j)} |\nabla (\mathbf{b}_j- b^+)|^2  \nonumber \\
&\leq 2^{n+ 2}\delta_1\cdot \fint_{B_{2\delta_1}(y_j)} |\Delta b^+| \leq D\cdot 2^{n+ 4}n\delta_1 ,\label{eq 1.23.1} \\
\fint_{B_{\delta_1}(y_j)} |\nabla^2 \mathbf{b}_j|^2  &\leq \frac{10^8n^3}{(\delta_1/2)^2}\Big[1+ \frac{d_p(y_j)}{\delta_1/2}\Big]^2\leq \frac{10^{10}n^3}{\delta_1^4} . \label{eq 1.23.2}
\end{align}

For any $j\in \mathbf{J}$, we can now apply Lemma \ref{lem good choice of balls} to the functions $|\nabla(\mathbf{b}_j- b^+)|^2$ and $|\nabla^2 \mathbf{b}_j|^2$ on $B_{\delta_1}(y_j)$, $j\in \mathbf{J}$. We get balls $\Big\{B_{\delta_2\delta_1}(x_i^j)\Big\}_{i, j}$, where $B_{2\delta_2\delta_1}(x_i^j)\subset B_{\delta_1}(y_j)$ and let $\delta= \delta_1 \delta_2$, such that 
\begin{align}
&\fint_{B_{2\delta}(x_i^j)} |\nabla (\mathbf{b}_j- b^+)|^2 \leq D\cdot \fint_{B_{\delta_1}(y_j)} |\nabla (\mathbf{b}_j- b^+)|^2\leq  D^2\cdot 2^{n+ 4}n\delta_1 ,\label{eq need 1} \\
&\fint_{B_{2\delta}(x_i^j)} |\nabla^2 \mathbf{b}_j|^2 \leq D\cdot \fint_{B_{\delta_1}(x_j)} |\nabla^2 \mathbf{b}_j|^2 \leq 
\frac{D\cdot 10^{10}n^3}{\delta_1^4}= \frac{D\cdot 10^{10}n^3(\delta_2 \delta_1^{-1})^2}{\delta^2} .\label{eq need 2} 
\end{align}
Furthermore, we set $\Omega_j= B_{\delta_1}(y_j)- \bigcup_i B_{\delta}(x_i^j)$, then 
\begin{align}
V(\Omega_j)&= V\big(B_{\delta_1}(y_j)\big)- V\big(\bigcup_i B_{\delta}(x_i^j)\big)\leq V\big(B_{\delta_1}(y_j)\big)- \Big[(1- 2\delta_2)^n- \frac{2\cdot 10^n}{D}\Big]V\big(B_{\delta_1}(y_j)\big) \nonumber \\
&= \Big[1+ \frac{2\cdot 10^n}{D}- (1- 2\delta_2)^n\Big]V\big(B_{\delta_1}(y_j)\big) .\nonumber 
\end{align}
And we get 
\begin{align}
&\quad \frac{V\big(\bigcup_{i, j} B_{\delta}(x_i^j)\big)}{V\big(B_1(p)\big)}= \frac{V\Big(\bigcup_j \big[B_{\delta_1}(y_j)- \Omega_j \big]\Big)}{V\big(B_1(p)\big)} \geq \frac{V\Big(\bigcup_j B_{\delta_1}(y_j)\Big)- \Sigma_j V(\Omega_j)}{V\big(B_1(p)\big)} \nonumber \\
&\geq \Big[(1- 2\delta_1)^n- \frac{2\cdot 10^n}{D}\Big]- \Big[1+ \frac{2\cdot 10^n}{D}- (1- 2\delta_2)^n\Big]\sum_j \frac{V\big(B_{\delta_1}(y_j)\big)}{V\big(B_1(p)\big)} \nonumber \\
&\geq \Big[(1- 2\delta_1)^n- \frac{2\cdot 10^n}{D}\Big]- 2^n\cdot\Big[1+ \frac{2\cdot 10^n}{D}- (1- 2\delta_2)^n\Big]\sum_j \frac{V\big(B_{\frac{1}{2}\delta_1}(y_j)\big)}{V\big(B_1(p)\big)} \nonumber \\
&\geq \Big\{(1- 2\delta_1)^n- \frac{2\cdot 10^n}{D}- 2^n\cdot\big[1+ \frac{2\cdot 10^n}{D}- (1- 2\delta_2)^n\big]\Big\}  .\label{eq need 3} 
\end{align}

If we choose $\delta_2= \delta_1^2$, $D= \delta_1^{-\frac{1}{3}}\cdot 10^n$. Using $(1- 2\delta_1)^n\geq 1- 2n\delta_1$, then for suitable $C(n)$, when $\delta_1= 2^{-50n}\epsilon^6$, $\delta= \delta_2\delta_1= 2^{-150n}\epsilon^{18}$, we get 
\begin{align}
D^2\cdot 2^{n+ 4}n\delta_1&\leq \epsilon^2 ,\nonumber \\
D\cdot 10^{10}n^3(\delta_2\delta_1^{-1})^2&\leq \epsilon^2 ,\nonumber \\
(1- 2\delta_1)^n- \frac{2\cdot 10^n}{D}- 2^n\cdot\big[1+ \frac{2\cdot 10^n}{D}- (1- 2\delta_2)^n\big]&\geq 1- \epsilon .\nonumber 
\end{align}
From (\ref{eq need 1}), (\ref{eq need 2}) and (\ref{eq need 3}), we have 
\begin{align}
\fint_{B_{2\delta r}(x_i)} |\nabla (\mathbf{b}_i- b^+)|^2 \leq \epsilon^2 \ , \quad \quad  \quad
\fint_{B_{2\delta r}(x_i)} |\nabla^2 \mathbf{b}_i|^2 \leq \frac{\epsilon^2}{\delta^2r^2} .\nonumber 
\end{align}
Using Cauchy-Schwartz inequality, the conclusion follows.
}
\qed

Let $SM^n$ be the unit tangent bundle of $M^n$, if $\pi: SM^n\rightarrow M^n$ is the projection map, for any $\Omega\subset SM^n$, the \textbf{Liouville measure} of $\Omega$, denoted by $V(\Omega)$, is defined by $V(\Omega)= \mu\big(\pi(\Omega)\big)\cdot V(\mathbb{S}^{n- 1})$, where $V(\mathbb{S}^{n- 1})$ is the volume of the conical Euclidean $(n- 1)$-sphere, and $\mu$ is the volume measure of $(M^n, g)$ determined by the metric $g$.

\begin{definition}\label{def geodesic flow}
{For $\nu\in S_xM^n$, let $\gamma_{\nu}(\cdot)$ be the geodesic starting from $x$ with $\gamma_{\nu}'(0)= \nu$, the geodesic flow $\mathfrak{g}^t(x, \nu): [0, \infty)\times SM^n\rightarrow SM^n$ is defined by
\begin{align}
\mathfrak{g}^t(x, \nu)= \big(\gamma_{\nu}(t), \gamma_{\nu}'(t)\big)\ , \quad \quad \quad \quad \forall t\geq 0 .\nonumber 
\end{align}
}
\end{definition}

\begin{theorem}[Liouville's Theorem]\label{thm Liouville's thm}
{For any region $D\subset SM^n$ we have $V(\mathfrak{g}^tD)= V(D)$, where $\mathfrak{g}^t: SM^n\rightarrow SM^n$ is the geodesic flow on $M^n$, and the measure on $D, \mathfrak{g}^tD$ is the Liouville measure.
}
\end{theorem}

\pf
{\cite{Arnold}.
}
\qed

\begin{lemma}\label{lem approximate by functions}
{Let $x\in M^n$, $l, r> 0$, suppose $f\in C^{\infty}\big(B_{r+ l}(x)\big)$, and $g$ is a Lipschitz function on $M^n$, then for any $0\leq t\leq l$,  
\begin{align}
\fint_{SB_r(x)} \Big|\big(g\circ \gamma_v\big)'(t)- \frac{g\big(\gamma_v(l)\big)- g\big(\gamma_v(0)\big)}{l}\Big|\leq \frac{2}{V\big(B_r(x)\big)} \int_{B_{r+ l}(x)}\Big[ l|\nabla^2 f|+ \big|\nabla (g- f)\big|\Big] .\nonumber 
\end{align}
}
\end{lemma}

\pf
{Let $h(\tau)= f\circ \gamma_v(\tau)$, then for any $0\leq t\leq l$
\begin{align}
\big|h'(t)- \frac{h(l)- h(0)}{l}\big|&= \Big|\int_0^t h''(\tau) d\tau+ h'(0)- \frac{\int_0^l \big[\int_0^s h''(\tau) d\tau+ h'(0)\big] ds}{l}\Big| \nonumber \\
&\leq \int_0^t |h''|d\tau + \frac{\int_0^l \int_0^s |h''| d\tau ds}{l}\leq 2\int_0^l |h''| d\tau .\nonumber 
\end{align}

%We define 
%\begin{align}
%(I)= \fint_{SB_r(x)} \Big|h'(t)- \frac{h(l)- h(0)}{l}\Big| \nonumber 
%\end{align}
Then
\begin{align}
\fint_{SB_r(x)} \Big|h'(t)- \frac{h(l)- h(0)}{l}\Big|&\leq \frac{2}{V\big(SB_r(x)\big)}\int_{SB_r(x)} \Big(\int_0^l \big|\frac{\partial^2}{\partial \tau^2}(f\circ \gamma_v)\big| d\tau\Big) \nonumber \\
&= \frac{2}{V\big(B_r(x)\big)V(\mathbb{S}^{n- 1})}\int_0^l \Big(\int_{SB_r(x)}\big|\frac{\partial^2}{\partial \tau^2}(f\circ \gamma_v)\big| \Big) d\tau \nonumber \\
&\leq \frac{2l}{V\big(B_r(x)\big)} \int_{B_{r+ l}(x)} |\nabla^2 f| .\label{eq 1.16.2}
\end{align}

For any $0\leq t\leq l$, applying Theorem \ref{thm Liouville's thm} and (\ref{eq 1.16.2}),
\begin{align}
&\fint_{SB_r(x)} \Big|\big(g\circ \gamma_v\big)'(t)- \frac{g\big(\gamma_v(l)\big)- g\big(\gamma_v(0)\big)}{l}\Big|\leq \fint_{SB_r(x)} \Big|\big(g\circ \gamma_v\big)'(t)- \big(f\circ \gamma_v\big)'(t)\Big| \nonumber \\
&\quad \quad + \fint_{SB_r(x)} \Big|h'(t)- \frac{h(l)- h(0)}{l}\Big|+ \fint_{SB_r(x)} \frac{\int_0^l \Big|\big((f- g)\circ \gamma_v\big)'\Big|ds}{l} \nonumber \\
&\quad  \leq \frac{2l}{V\big(B_r(x)\big)} \int_{B_{r+ l}(x)} |\nabla^2 f|+ \frac{1}{V\big(SB_r(x)\big)}\int_{SB_{r+ l}(x)} \Big|\big(g\circ \gamma_v\big)'(0)- \big(f\circ \gamma_v\big)'(0)\Big| \nonumber \\
&\quad \quad + \frac{1}{V\big(SB_r(x)\big)}\cdot \frac{1}{l}\int_0^l \Big[\int_{SB_r(x)} \Big|\big((f- g)\circ \gamma_v\big)'(s)\Big| \Big] ds \nonumber \\
&\quad \leq \frac{2l}{V\big(B_r(x)\big)} \int_{B_{r+ l}(x)} |\nabla^2 f|+ \frac{2}{V\big(B_r(x)\big)} \int_{B_{r+ l}(x)} \big|\nabla (g- f)\big| . \nonumber 
\end{align}
}
\qed

Now we prove Colding's integral Toponogov theorem in quantitative form.

\begin{theorem}\label{thm Colding-1.28}
{Suppose $Rc(M^n)\geq 0$, $p, q\in M$ with $d(p, q)> 2r$. For all $1> \epsilon> 0$, there exists $\delta= 2^{-240n}\epsilon^{18}$ such that for all $0\leq t\leq \delta$, 
\begin{align}
\fint_{SB_{r}(p)} \Big|(b^+\circ \gamma_v)'(tr)- \frac{(b^+\circ \gamma_v)(\delta r)- (b^+\circ \gamma_v)(0)}{\delta r}\Big|\leq \epsilon  ,\label{eq 2.20.1}
\end{align}
where $b^+(x)= d(x, q)- d(p, q)$.
}
\end{theorem}

\pf
{From Lemma \ref{lem appr by harm func}, for $\epsilon_1> 0$ (to be determined later), we can find finitely many balls $B_{\delta r}(x_i)\subset B_r(p)$ with $\delta= 2^{-150n}\epsilon_1^{18}$ such that 
\begin{align}
&V\big(\bigcup_{i\in \mathbf{A}} B_{\delta r}(x_i)\big)\geq (1- \epsilon_1)V\big(B_r(p)\big)\ , \quad \quad \quad \quad \sum_{i}V\big(B_{\delta r}(x_i)\big)\leq 2^{2n} V\big(B_r(p)\big) , \nonumber \\
&\fint_{B_{2\delta r}(x_i)} |\nabla (\mathbf{b}_i- b^+)|\leq \epsilon_1 \ , \quad \quad \quad \quad
\fint_{B_{2\delta r}(x_i)} |\nabla^2 \mathbf{b}_i| \leq \frac{\epsilon_1}{\delta r} .\nonumber  
\end{align}

Let $h(\nu)= (b^+\circ \gamma_v)'(tr)- \frac{(b^+\circ \gamma_v)(\delta r)- (b^+\circ \gamma_v)(0)}{\delta r}$, apply Lemma \ref{lem approximate by functions}, we get 
\begin{align}
&\quad \fint_{SB_{\delta r}(x_i)} \big|h(\nu)\big| \leq \frac{2\delta r}{V\big(B_{\delta r}(x_i)\big)}\int_{B_{2\delta r}(x_i)} |\nabla^2 \mathbf{b}_i|+ \frac{2}{V\big(B_{\delta r}(x_i)\big)}\int_{B_{2\delta r}(x_i)} |\nabla (\mathbf{b}_i- b^+)| 
\nonumber \\
& \leq 2\delta r\frac{V\big(B_{2\delta r}(x_i)\big)}{V\big(B_{\delta r}(x_i)\big)}\fint_{B_{2\delta r}(x_i)} |\nabla^2 \mathbf{b}_i|+ \frac{2V\big(B_{2\delta r}(x_i)\big)}{V\big(B_{\delta r}(x_i)\big)} \fint_{B_{2\delta r}(x_i)} |\nabla (\mathbf{b}_i- b^+)| \nonumber \\
&\leq 2^{n+ 2} \epsilon_1 .\nonumber 
\end{align}

Note $|\nabla b^+|\leq 1$, then we have 
\begin{align}
\int_{SB_r(p)} \big|h(\nu)\big|& \leq \int_{\bigcup_i SB_{\delta r}(x_i)} \big|h(\nu)\big|+ 2V\Big(SB_r(p)- \bigcup_i SB_{\delta r}(x_i)\Big) \nonumber \\
&\leq \sum_i \int_{SB_{\delta r}(x_i)} \big|h(\nu)\big|+ 2V(\mathbb{S}^{n- 1})V\Big(B_r(p)- \bigcup_i B_{\delta r}(x_i)\Big)  \nonumber \\
&\leq 2^{n+ 2}\epsilon_1 \sum_{i} V\big(SB_{\delta r}(x_i)\big)+ 2V(\mathbb{S}^{n- 1})\cdot \epsilon_1 V\big(B_r(p)\big) \nonumber \\
&=  2^{n+ 2}\epsilon_1 \cdot V(\mathbb{S}^{n- 1})\sum_{i} V\big(B_{\delta r}(x_i)\big)+ 2V(\mathbb{S}^{n- 1})\cdot \epsilon_1 V\big(B_r(p)\big)\nonumber \\
&\leq (2^{3n+ 2}+ 2)\epsilon_1\cdot V(\mathbb{S}^{n- 1}) V\big(B_{r}(p)\big)= (2^{3n+ 2}+ 2)\epsilon_1\cdot V\big(SB_{r}(p)\big) .\nonumber 
\end{align}

Let $\epsilon_1= 2^{-5n}\epsilon$, then we have the conclusion for $\delta= 2^{-240n}\epsilon^{18}$.
}
\qed

\begin{prop}\label{prop approximation component is almost orthogonal}
{Given $0< \epsilon< 1, 0\leq k\leq n, B_{10r}(q)\subset (M^n, g)$ with $Rc(g)\geq 0$, any $\delta\leq n^{-1250n}\cdot \epsilon^{100n}\big(\frac{r_1}{r}\big)$ and $\frac{r_1}{r}\leq n^{-110n}\epsilon^{10n}$. Assume there is an $(\delta r)$-Gromov-Hausdorff approximation $\Phi: B_{10r}(0, \hat{q})\rightarrow B_{10r}(q)$, where $B_{10r}(0, \hat{q})\subset \mathbb{R}^k\times \mathbf{X}_k$, and there are $\hat{q}_0, \hat{q}_1^+\in B_{3r}(\hat{q})\subset \mathbf{X}_k$ satisfying $d(\hat{q}_0, \hat{q}_1^+)= r_0\in \big[\frac{1}{16}r, 2r\big]$. Then we have    
\begin{align}
\fint_{B_{r_1}(q_1)}\sum_{i, j= 1}^{k+ 1} \big|\langle \nabla b_i^+, \nabla b_j^+\rangle- \delta_{ij}\big|\leq \epsilon ,\label{almost o.n.}
\end{align}
where $\big\{\mathbf{e}_i\big\}_{i= 1}^k$ is the standard basis for $\mathbb{R}^k, q_1= \Phi(0, \hat{q}_0)$ and 
\begin{align}
b_i^+(\cdot)&= d(\cdot, p_i^+)- d(q_1, p_i^+)\ , \quad \quad  \quad \quad p_i^+= \Phi(r_0\cdot \mathbf{e}_i, \hat{q}_0)\ , \quad \quad \quad \quad 1\leq i\leq k \nonumber \\
b_{k+ 1}^+(\cdot)&= d(\cdot, p_{k+ 1}^+)- d(q_1, p_{k+ 1}^+)\ ,  \quad \quad \quad \quad p_{k+ 1}^+= \Phi(0, \hat{q}_1^+) .\nonumber 
\end{align}
}
\end{prop}

\pf
{\textbf{Step (1)}. Firstly we have 
\begin{align}
d\big((r_0\mathbf{e}_i, \hat{q}_0), (0, \hat{q})\big)\leq r_0^2+ d(\hat{q}_0, \hat{q})^2\leq (2r)^2+ (3r)^2< (10r)^2 ,\nonumber 
\end{align}
which implies $(r_0\mathbf{e}_i, \hat{q}_0)\in B_{10r}(0, \hat{q})$, and $p_i^+$ is well-defined. Also we can easily see $B_{r_1}(q_1)\subset B_{10r}(q)$.

In this step we always assume that $x, y\in B_{r_1}(q_1)$. Because $\Phi$ is an $(\delta r)$-Gromov-Hausdorff approximation, there exists $(\tilde{x}, \hat{x}), (\tilde{y}, \hat{y})\in B_{10r}(0, \hat{q})$ such that 
\begin{align}
d\big(\Phi(\tilde{x}, \hat{x}), x\big)< \delta r \quad \quad  and \quad \quad d\big(\Phi(\tilde{y}, \hat{y}), y\big)< \delta r .\nonumber 
\end{align}

Let $d_0= d\big((\tilde{x}, \hat{x}), (\tilde{y}, \hat{y})\big)$, then 
\begin{align}
d_0&\leq \Big|d\big(\Phi(\tilde{x}, \hat{x}), \Phi(\tilde{y}, \hat{y})\big)- d_0\Big|+ d\big(x, \Phi(\tilde{x}, \hat{x})\big)+ d\big(y, \Phi(\tilde{y}, \hat{y})\big)+ d(x, y) \nonumber \\
&\leq 3\delta r+ d(x, y) .\label{mid 1.2}
\end{align}

Assume $r_1= \zeta r$, we have 
\begin{align}
|\tilde{x}|^2+ d(\hat{x}, \hat{q}_0)^2&= d\big((\tilde{x}, \hat{x})), (0, \hat{q}_0)\big)^2\leq \Big[d\big(\Phi(\tilde{x}, \hat{x}), \Phi(0, \hat{q}_0)\big)+ \delta r\Big]^2 \nonumber \\
&\leq\Big[ d(x, q_1)+ 2\delta r\Big]^2 \leq (r_1+ 2\delta r)^2= (\zeta+ 2\delta)^2r^2 , \nonumber 
\end{align}
which implies 
\begin{align}
|\tilde{x}|\leq (\zeta+ 2\delta)r\quad \quad and \quad \quad d(\hat{x}, \hat{q}_0)\leq (\zeta+ 2\delta)r .\label{tilde 3.3}
\end{align}
Similarly we have 
\begin{align}
|\tilde{y}|\leq (\zeta+ 2\delta)r\quad \quad and \quad \quad d(\hat{y}, \hat{q}_0)\leq (\zeta+ 2\delta)r . \label{tilde 3.3-y}
\end{align}

We define $\tilde{x}= (\tilde{x}_1, \cdots, \tilde{x}_k), \tilde{y}= (\tilde{y}_1, \cdots, \tilde{y}_k)$ and 
\begin{align}
\mathcal{I}_i&= \Big|d\big((\tilde{x}, \hat{x}), (r_0\mathbf{e}_i, \hat{q}_0)\big)- d\big((\tilde{y}, \hat{y}), (r_0\mathbf{e}_i, \hat{q}_0)\big)\Big| \nonumber \\
\mathcal{J}&= \Big|d\big((\tilde{x}, \hat{x}), (0, \hat{q}_1^+)\big)- d\big((\tilde{y}, \hat{y}), (0, \hat{q}_1^+)\big)\Big| .\nonumber 
\end{align}
From (\ref{mid 1.2}), (\ref{tilde 3.3}) and (\ref{tilde 3.3-y}),
\begin{align}
\mathcal{I}_i&= \Big|\sqrt{d(\tilde{x}, r_0\mathbf{e}_i)^2+ d(\hat{x}, \hat{q}_0)^2}- \sqrt{d(\tilde{y}, r_0\mathbf{e}_i)^2+ d(\hat{y}, \hat{q}_0)^2}\Big| \nonumber \\
& \leq \frac{|\sum_{j= 1}^k \tilde{x}_j^2- \sum_{j= 1}^k \tilde{y}_j^2|+ 2r_0|\tilde{x}_i- \tilde{y}_i|+ \big|d(\hat{x}, \hat{q}_0)^2- d(\hat{y}, \hat{q}_0)^2\big|}{\sqrt{d(\tilde{x}, r_0\mathbf{e}_i)^2+ d(\hat{x}, \hat{q}_0)^2}+ \sqrt{d(\tilde{y}, r_0\mathbf{e}_i)^2+ d(\hat{y}, \hat{q}_0)^2}} \nonumber \\
&\leq \frac{2(\zeta+ 2\delta)r\cdot \big[\sum_{j= 1}^k|\tilde{x}_j- \tilde{y}_j|+ d(\hat{x}, \hat{y})\big]+ 2r_0\cdot |\tilde{x}_i- \tilde{y}_i|}{|r_0- \tilde{x}_i|+ |r_0- \tilde{y}_i|} \nonumber \\
&\leq \frac{2(\zeta+ 2\delta)r}{2r_0- 2(\zeta+ 2\delta)r}\sqrt{k+ 1}d_0+ \frac{2r_0}{2r_0- 2(\zeta+ 2\delta)r}|\tilde{x}_i- \tilde{y}_i| \nonumber \\
&\leq \frac{8n(\zeta+ 2\delta)}{1- 16(\zeta+ 2\delta)}\big[d(x, y)+ 3\delta r\big]+ \frac{1}{1- 16(\zeta+ 2\delta)}|\tilde{x}_i- \tilde{y}_i| ,\nonumber 
\end{align}
in the last inequality above we used the assumption $r_0\geq \frac{1}{16}r$. Similarly, we have 
\begin{align}
\mathcal{J}&= \Big|\sqrt{|\tilde{x}|^2+ d(\hat{x}, \hat{q}_1^+)^2}- \sqrt{|\tilde{y}|^2+ d(\hat{y}, \hat{q}_1^+)^2} \Big| \nonumber \\
&\leq \frac{|\sum_{j= 1}^k \tilde{x}_j^2- \sum_{j= 1}^k \tilde{y}_j^2|+ \big|d(\hat{x}, \hat{q}_1^+)^2- d(\hat{y}, \hat{q}_1^+)^2\big|}{\sqrt{|\tilde{x}|^2+ d(\hat{x}, \hat{q}_1^+)^2}+ \sqrt{|\tilde{y}|^2+ d(\hat{y}, \hat{q}_1^+)^2}} \nonumber \\
&\leq \frac{2(\zeta+ 2\delta)r\sum_{j= 1}^k|\tilde{x}_j- \tilde{y}_j|+ 2(\zeta r+ 2\delta r+ r_0)\cdot d(\hat{x}, \hat{y})}{2r_0- 2(\zeta+ 2\delta)r} \nonumber \\
&\leq \frac{2n(\zeta+ 2\delta)rd_0}{2r_0- 2(\zeta+ 2\delta)r}+ \frac{1}{1- (\zeta+ 2\delta)\frac{r}{r_0}}d(\hat{x}, \hat{y}) \nonumber \\
&\leq \frac{8n(\zeta+ 2\delta)}{1- 16(\zeta+ 2\delta)}\big[d(x, y)+ 3\delta r\big]+ \frac{1}{1- 16(\zeta+ 2\delta)}d(\hat{x}, \hat{y}) . \nonumber 
\end{align}
Then we get 
\begin{align}
\sqrt{\sum_{i= 1}^k \mathcal{I}_i^2+ \mathcal{J}^2}&\leq \frac{1}{1- 16(\zeta+ 2\delta)}d_0+ \frac{8n^2(\zeta+ 2\delta)}{1- 16(\zeta+ 2\delta)}\big[d(x, y)+ 3\delta r\big] \nonumber \\
&\quad + \frac{8n(\zeta+ 2\delta)}{1- 16(\zeta+ 2\delta)}\sqrt{d(x, y)+ 3\delta r}\sqrt{d_0} \nonumber \\
&\leq \Big(\frac{1+ 16n^2(\zeta+ 2\delta)}{1- 16(\zeta+ 2\delta)}\Big)\big[d(x, y)+ 3\delta r\big] .\label{I and J}
\end{align}

For $i\leq k$, we have 
\begin{align}
\big|b_i^+(x)- b_i^+(y)\big|&= \big|d(x, p_i^+)- d(y, p_i^+)\big|\nonumber \\
&\leq \Big|d\big(\Phi(\tilde{x}, \hat{x}), \Phi(r_0\cdot \mathbf{e}_i, \hat{q}_0)\big)- d\big((\tilde{y}, \hat{y}), \Phi(r_0\cdot \mathbf{e}_i, \hat{q}_0)\big)\Big|+ 2\delta r \nonumber \\
&\leq 4\delta r+ \Big|d\big((\tilde{x}, \hat{x}), (r_0\cdot \mathbf{e}_i, \hat{q}_0)\big)- d\big((\tilde{y}, \hat{y}), (r_0\cdot \mathbf{e}_i, \hat{q}_0)\big)\Big| \nonumber \\
&\leq 4\delta r+ \mathcal{I}_i\label{eq 2.27.0} \\
\big|b_{k+ 1}^+(x)- b_{k+ 1}^+(y)\big|&= \big|d(x, p_{k+ 1}^+)- d(y, p_{k+ 1}^+)\big|\nonumber \\
&\leq \Big|d\big(\Phi(\tilde{x}, \hat{x}), \Phi(0, \hat{q}_1^+)\big)- d\big(\Phi(\tilde{y}, \hat{y}), \Phi(0, \hat{q}_1^+)\big)\Big|+ 2\delta r \nonumber \\
&\leq 4\delta r+ \Big|d\big((\tilde{x}, \hat{x}), (0, \hat{q}_1^+)\big)- d\big((\tilde{y}, \hat{y}), (0, \hat{q}_1^+)\big)\Big|  \nonumber \\
&\leq 4\delta r+ \mathcal{J} .\label{eq 2.27.0-k+1}
\end{align}

Let $\Psi= (b_1^+, \cdots, b_k^+, b_{k+ 1}^+): M^n \rightarrow \mathbb{R}^{k+ 1}$, from (\ref{mid 1.2}), (\ref{eq 2.27.0}), (\ref{eq 2.27.0-k+1}) and (\ref{I and J}), 
\begin{align}
&\quad d\big(\Psi(x), \Psi(y)\big)= \sqrt{\sum_{i= 1}^{k+ 1} \big|b_i^+(x)- b_i^+(y)\big|^2} \nonumber \\
&\leq \Big\{\sum_{i= 1}^k \Big(16\delta^2r^2+ 8\delta r\cdot d_0+ \mathcal{I}_i^2 \Big)  + 16\delta^2r^2+ 8\delta r\cdot d_0+ \mathcal{J}^2\Big\}^{\frac{1}{2}} \nonumber \\
&\leq 4n\delta r+ 4n\delta \sqrt{rd_0}+ \sqrt{\sum_{i= 1}^k\mathcal{I}_i^2+ \mathcal{J}^2}  \nonumber \\
&\leq 20n\delta \cdot r+ 4n\delta\cdot d(x, y)+ \Big(\frac{1+ 16n^2(\zeta+ 2\delta)}{1- 16(\zeta+ 2\delta)}\Big)\big[d(x, y)+ 3\delta r\big] . \nonumber 
\end{align}

If we assume $d(x, y)= l r_1$, where $0< l< 1$ is to be determined later, then
\begin{align}
\frac{d\big(\Psi(x), \Psi(y)\big)}{d(x, y)}&\leq 4n\delta+ \frac{1+ 16n^2(\zeta+ 2\delta)}{1- 16(\zeta+ 2\delta)}+ \Big(20n+ 3\cdot \frac{1+ 16n^2(\zeta+ 2\delta)}{1- 16(\zeta+ 2\delta)}\Big)\frac{\delta \cdot r}{d(x, y)} \nonumber \\
&\leq \frac{1+ 16n^2(\zeta+ 2\delta)}{1- 16(\zeta+ 2\delta)}+ \Big\{4n\delta+ \Big(20n+ 3\cdot \frac{1+ 16n^2(\zeta+ 2\delta)}{1- 16(\zeta+ 2\delta)}\Big)\frac{\delta}{l\zeta}\Big\} . \label{need 2.21.1}
\end{align}

For $0< \epsilon_1< 1$ (to be determined later), if the following holds: 
\begin{align}
\zeta\leq \frac{(\frac{\epsilon_1}{4})}{400n^2} \quad \quad and \quad \quad \delta\leq \frac{(\frac{\epsilon_1}{4})l}{800n^2}\cdot \zeta , \label{def of R and delta}
\end{align}
from (\ref{need 2.21.1}) and (\ref{def of R and delta}), we obtain 
\begin{align}
\frac{d\big(\Psi(x), \Psi(y)\big)}{d(x, y)}\leq 1+ \frac{\epsilon_1}{4} \ , \quad \quad \quad \quad if\ d(x, y)= l r_1\ and\  x, y\in B_{r_1}(q_1) . \label{Lip dist Phi is close to 1} 
\end{align}

From (\ref{Lip dist Phi is close to 1}), for all $\nu\in SB_{r_1}(q_1)$ except a zero-measure set, we have  
\begin{align}
\frac{d\big(\Psi\circ \gamma_{\nu}(lr_1), \Psi\circ \gamma_{\nu}(0)\big)}{lr_1}\leq 1+ \frac{\epsilon_1}{4} , \label{need ha-1}
\end{align}
where $\gamma_{\nu}$ is the geodesic satisfying $\gamma_{\nu}'(0)= \nu$.

\textbf{Step (2)}. From Theorem \ref{thm Colding-1.28}, let 
\begin{align}
l= 2^{-240n}\epsilon_1^{18} , \label{def of l}
\end{align}
then we have  
\begin{align}
\frac{1}{V\big(SB_{r_1}(q_1)\big)} \int_{SB_{r_1}(q_1)} \Big|\langle \nabla b_i^+, \nu \rangle- f(\nu)\Big|\leq \epsilon_1 \ , \quad \quad \quad \quad i= 1, \cdots, k+1 , \label{need 2.11}
\end{align}
where $f(\nu)\vcentcolon= \frac{(b_i^+\circ \gamma_{\nu})(lr_1)- (b_i^+\circ \gamma_{\nu})(0)}{lr_1}$ for $\nu\in SB_{r_1}(q_1)$.

For fixed $i$, and some $\theta\in (0, \frac{\pi}{2})$ (to be determined later), set 
\begin{align}
C_{\theta}= \Big\{\nu\in SB_{r_1}(q_1)|\ \angle(\nu, \nabla b_i^+)\leq  \theta \Big\}= \Big\{\nu\in SB_{r_1}(q_1)|\ \langle \nu, \nabla b_i^+\rangle\geq  \cos\theta \Big\} ,\nonumber 
\end{align}
then for any $\nu\in C_{\theta}$, 
\begin{align}
|\nu- \nabla b_i^+|^2= |\nu|^2+ |\nabla b_i^+|^2- 2\langle \nu, \nabla b_i^+\rangle \leq 2- 2\cos\theta .\label{need length difference}
\end{align}

From (\ref{need 2.11}),
\begin{align}
&\quad \frac{1}{V\big(SB_{r_1}(q_1)\big)}\int_{C_{\theta}} \Big|\big|f(\nu)\big|- 1\Big|  \nonumber \\
&\leq \frac{1}{V\big(SB_{r_1}(q_1)\big)}\int_{C_{\theta}} \Big|\big|f(\nu)\big|- \langle \nabla b_i^+, \nu \rangle\Big| +  \frac{1}{V\big(SB_{r_1}(q_1)\big)}\int_{C_{\theta}} \big|1- \langle \nabla b_i^+, \nu \rangle\big|  \nonumber \\
&\leq \epsilon_1+ (1- \cos\theta)\cdot \frac{V(C_{\theta})}{V\big(SB_{r_1}(q_1)\big)}\leq \epsilon_1+ \theta \frac{V(C_{\theta})}{V\big(SB_{r_1}(q_1)\big)} .\label{need lem 4.4.1}
\end{align}
%then 
%\begin{align}
%&\frac{1}{V\big(SB_1(p)\big)}\int_{C_{\theta}} \Big|\frac{(b_i^+\circ \gamma)(l)- (b_i^+\circ \gamma_{\nu})(0)}{l}\Big|^2\geq \Big(\frac{1}{V\big(SB_1(p)\big)}\int_{C_{\theta}} \Big|\frac{(b_i^+\circ \gamma)(l)- (b_i^+\circ \gamma_{\nu})(0)}{l}\Big|\Big)^2 \nonumber \\
%&\quad \geq \Big(1- C(n)\epsilon- (1- \cos\theta)\cdot \frac{V(C_{\theta})}{V\big(SB_1(p)\big)}\Big)^2
%\end{align}

Now from (\ref{need 2.11}), (\ref{need lem 4.4.1}) and (\ref{need ha-1}), for $j\neq i$, we have 
\begin{align}
&\quad \frac{1}{V\big(SB_{r_1}(q_1)\big)} \int_{C_{\theta}} \big|\langle \nabla b_j^+, \nu \rangle\big|\nonumber \\
&\leq \frac{1}{V\big(SB_{r_1}(q_1)\big)} \int_{SB_{r_1}(q_1)} \Big|\langle \nabla b_j^+, \nu \rangle- f(\nu)\Big|+ \frac{1}{V\big(SB_{r_1}(q_1)\big)} \int_{C_{\theta}} \big|f(\nu)\big| \nonumber \\
& \leq \epsilon_1 + \frac{V(C_{\theta})^{\frac{1}{2}}}{V\big(SB_{r_1}(q_1)\big)}\Big(\int_{C_{\theta}} \Big|\frac{d\big(\Psi\circ \gamma_{\nu}(lr_1), \Psi\circ \gamma_{\nu}(0)\big)}{lr_1}\Big|^2- \big|f(\nu)\big|^2\Big)^{\frac{1}{2}} \nonumber \\
&\leq \epsilon_1+ \frac{V(C_{\theta})^{\frac{1}{2}}}{V\big(SB_{r_1}(q_1)\big)}\Big\{\epsilon_1V(C_{\theta})^{\frac{1}{2}}+ \Big(\int_{C_{\theta}} \Big|1- \big|f(\nu)\big|^2\Big|\Big)^{\frac{1}{2}}\Big\} \nonumber \\
&\leq \epsilon_1 + \frac{V\big(C_{\theta}\big)}{V\big(SB_{r_1}(q_1)\big)} \epsilon_1+ \sqrt{2}\frac{V(C_{\theta})^{\frac{1}{2}}}{V\big(SB_{r_1}(q_1)\big)}\cdot \Big(\int_{C_{\theta}} \Big|\big|f(\nu)\big|- 1\Big| \Big)^{\frac{1}{2}} \nonumber \\
&\leq 2\epsilon_1+ \sqrt{2}\Big(\frac{V(C_{\theta})}{V\big(SB_{r_1}(q_1)\big)}\Big)^{\frac{1}{2}}\cdot \Big(\epsilon_1+ \theta\cdot \frac{V(C_{\theta})}{V\big(SB_{r_1}(q_1)\big)}\Big)^{\frac{1}{2}}\leq 4\sqrt{\epsilon_1}+ 2\sqrt{\theta}\frac{V(C_{\theta})}{V\big(SB_{r_1}(q_1)\big)} . \nonumber 
\end{align}

Note (\ref{need length difference}), then 
\begin{align}
\int_{C_{\theta}} \big|\langle \nabla b_i^+, \nabla b_j^+\rangle\big|&\leq \int_{C_{\theta}} \big|\langle \nabla b_{j}^+, v- \nabla b_i^+\rangle\big|+ \int_{C_{\theta}} \big|\langle \nabla b_j^+, \nu\rangle \big| \nonumber \\
&\leq \int_{C_{\theta}} \big|v- \nabla b_i^+\big|+ \int_{C_{\theta}} \big|\langle \nabla b_j^+, \nu\rangle \big|\leq \theta \cdot V(C_{\theta})+ \int_{C_{\theta}} \big|\langle \nabla b_j^+, \nu\rangle \big| .\nonumber 
\end{align}

Now note $\langle \nabla b_i^+, \nabla b_j^+\rangle$ is constant on $T_xM^n$ for any fixed $x\in M^n$, we have 
\begin{align}
\fint_{B_{r_1}(q_1)} \big|\langle \nabla b_i^+, \nabla b_j^+\rangle \big|&= \frac{1}{V\big(C_{\theta}\big)}\int_{C_{\theta}} \big|\langle \nabla b_i^+, \nabla b_j^+\rangle\big| \nonumber \\
&\leq \theta+ \frac{V\big(SB_{r_1}(q_1)\big)}{V(C_{\theta})}\cdot \Big\{4\sqrt{\epsilon_1}+ 2\sqrt{\theta}\frac{V(C_{\theta})}{V\big(SB_{r_1}(q_1)\big)}\Big\} \nonumber \\
&\leq 3\sqrt{\theta}+ 4\sqrt{\epsilon_1}\cdot \frac{V\big(SB_{r_1}(q_1)\big)}{V(C_{\theta})} .\nonumber 
\end{align}

Let $\epsilon_1\leq \Big(\frac{V(C_{\theta})}{4V\big(SB_{r_1}(q_1)\big)}\Big)^2\cdot \theta$ and $\theta= \frac{\epsilon^2}{16}\cdot \frac{1}{(n+ 1)^4}$, then $\fint_{B_{r_1}(q_1)} \big|\langle \nabla b_i^+, \nabla b_j^+\rangle \big|\leq 4\sqrt{\theta}= \frac{\epsilon}{(n+ 1)^2}$. So we let $\epsilon_1= 2n^{-50n}\epsilon^{4n+ 2}$, the above conclusion follows. Plug into (\ref{def of l}) and (\ref{def of R and delta}), the corresponding $\zeta, \delta$ can be determined. 
}
\qed

%\section{Approximate local Busemann function by harmonic function}

\section{Existence of almost orthonormal linear (A.O.L.) harmonic functions}

\begin{definition}\label{def excess est}
{For $q^+, q^-, p\in \mathbf{X}$, where $\mathbf{X}$ is a metric space, we say that $[q^+, q^-, p]$ is an \textbf{AG-triple on $\mathbf{X}$ with the excess $s$ and the scale $t$} if 
\begin{align}
\mathbf{E}(p)= s \quad\quad and \quad \quad 
\min\big\{d(p, q^+), d(p, q^-)\big\}= t ,\nonumber 
\end{align} 
where $\mathbf{E}(\cdot)= d(\cdot, q^+)+ d(\cdot, q^-)- d(q^+, q^-)$.
}
\end{definition}

For $p\in M^n, r_2\geq r_1\geq 0$, we define $A_{r_1, r_2}(p)$ as the following:
\begin{align}
A_{r_1, r_2}(p)= \{x\in M^n|\ r_1< d(x, p)< r_2\} . \nonumber 
\end{align}

We have the following Abresch-Gromoll lemma.
\begin{lemma}\label{lem general Abresch-Gromoll}
{On complete Riemannian manifold $(M^n, g)$ with $Rc\geq 0$, assume that $[q^+, q^-, p]$ is an AG-triple with the excess $\leq \frac{1}{n}\frac{r^2}{R}$ and the scale $\geq R$, furthermore assume $R\geq 2^{2n}r$, then $\sup_{B_{r}(p)} \mathbf{E}\leq 2^6\cdot \big(\frac{r}{R}\big)^{\frac{1}{n- 1}}r$, where $\mathbf{E}(\cdot)= d(\cdot, q^+)+ d(\cdot, q^-)- d(q^+, q^-)$.
}
\end{lemma}

\pf
{Define the function $\varphi(\rho)= \int_{\rho}^2 \int_t^2 \big(\frac{s}{t}\big)^{n- 1}ds dt$, and $\epsilon$ solves 
\begin{align}
\frac{\frac{R}{r}- 3}{4(n- 1)}= \frac{\varphi\big(\frac{\epsilon}{4}\big)}{\epsilon}\ , \quad \quad \quad \quad 0< \epsilon\leq 4 .\label{def epsilon in A-G}
\end{align}
From $R\geq 2^{2n}r> \big[3+ 4(n- 1)\varphi(1)\big]r$, it is easy to see that $\epsilon$ exists and is unique.

By $R\geq 2^{2n}r$ and (\ref{def epsilon in A-G}), we have 
\begin{align}
\varphi(\frac{\epsilon}{4})= \epsilon\cdot \frac{\frac{R}{r}- 3}{4(n- 1)}\geq \frac{\epsilon}{8n}\frac{R}{r} .\label{A-G quantity-1}
\end{align}

Note for $\rho\leq 1$, 
\begin{align}
\varphi(\rho)= \frac{1}{n}\Big[\frac{1}{2}\rho^2+ \frac{2^n}{n- 2}\rho^{2- n}+ 4\big(-\frac{1}{n- 2}- \frac{1}{2}\big)\Big] \leq 2^n\rho^{2- n} . \label{A-G quantity-2}
\end{align}

From (\ref{A-G quantity-1}) and (\ref{A-G quantity-2}), we obtain $2^n\big(\frac{\epsilon}{4}\big)^{2- n}\geq \frac{\epsilon}{8n}\frac{R}{r}$, which implies $\epsilon\leq 2^6\cdot \big(\frac{r}{R}\big)^{\frac{1}{n- 1}} $.

To prove the conclusion, we only need to prove that $\sup_{B_r(p)} \mathbf{E}\leq \epsilon r$. By contradiction, if $\sup_{B_r(p)} \mathbf{E}> \epsilon r$, then there exists $x_0\in B_r(p)$, such that 
\begin{align}
\mathbf{E}(x_0)> \epsilon r .\label{max point}
\end{align}
From $|\nabla \mathbf{E}|\leq 2$ and (\ref{def epsilon in A-G}), we have 
\begin{align}
d(p, x_0)&\geq \frac{\mathbf{E}(x_0)- \mathbf{E}(p)}{2}> \frac{\epsilon r- \frac{1}{n}\frac{r^2}{R}}{2} \geq \frac{1}{2}\Big(\epsilon r-  \frac{2(n- 1)}{\frac{R}{r}- 3}\varphi(1) r\Big) \nonumber \\
&\geq  \frac{1}{2}\Big(\epsilon r-  \frac{2(n- 1)}{\frac{R}{r}- 3}\varphi\big(\frac{\epsilon}{4}\big) r\Big)= \frac{1}{2}\big(\epsilon r- \frac{\epsilon}{2}r\big)= \frac{\epsilon}{4}r ,\nonumber 
\end{align}
which implies $p\in A_{\frac{\epsilon}{4}r, 2r}(x_0)$.

Define the function $h: A_{\frac{\epsilon}{4}r, 2r}(x_0)\rightarrow [0, \infty)$ by $\epsilon\leq 2^6\cdot \big(\frac{r}{R}\big)^{\frac{1}{n- 1}}$, and it is easy to check $(\mathbf{E}- h)\geq 0$ on $\partial B_{2r}(x_0)$. 

For any $x\in \partial B_{\frac{\epsilon}{4}r}(x_0)$, using $|\nabla \mathbf{E}|\leq 2$, (\ref{max point}) and (\ref{def epsilon in A-G}), we have 
\begin{align}
\mathbf{E}(x)- h(x)\geq \mathbf{E}(x_0)- 2\cdot \frac{\epsilon}{4}r- \frac{2(n- 1)r}{\frac{R}{r}- 3}\varphi\big(\frac{\epsilon}{4}\big)> \frac{\epsilon}{2}r- \frac{2(n- 1)r}{\frac{R}{r}- 3}\varphi\big(\frac{\epsilon}{4}\big)= 0 ,\nonumber 
\end{align}
hence we obtain $\min_{x\in \partial A_{\frac{\epsilon}{4}r, 2r}(x_0)}(\mathbf{E}-h)(x)\geq 0$.

On the other hand, let $d= d(x, x_0)$, then $\varphi''(\frac{d}{r})+ \frac{n- 1}{d}r\cdot \varphi'(\frac{d}{r})= 1$. From $Rc\geq 0$ and Laplace Comparison Theorem $\Delta d\leq  \frac{n- 1}{d}$, note $\varphi'\leq 0$, for any $x\in A_{\frac{\epsilon}{4}r, 2r}(x_0)$, we have the following inequality in weak sense:
\begin{align}
\Delta h(x)&= \frac{2(n- 1)}{\frac{R}{r}- 3}r\big[r^{-2}\varphi''+ r^{-1}\varphi'\Delta d\big]\geq \frac{2(n- 1)}{R- 3r} \nonumber \\
\Delta \mathbf{E}(x)&\leq \frac{n- 1}{d(x, q^+)}+ \frac{n- 1}{d(x, q^-)}\leq \frac{2(n- 1)}{R- 3r}  .\nonumber
\end{align}
Then we get $\min_{x\in A_{\frac{\epsilon}{4}r, 2r}(x_0)}\Delta (\mathbf{E}- h)(x)\leq 0$.

From Weak Maximum Principle for weak superharmonic function (c.f. \cite[Theorem $8.1$]{GT}), note $p\in A_{\frac{\epsilon}{4}r, 2r}(x_0)$, we have $(\mathbf{E}- h)(p)\geq \min_{\partial A_{\frac{\epsilon}{4}r, 2r}(x_0)}(\mathbf{E}- h)\geq 0$. And note $d(p, x_0)< r$, we obtain
\begin{align}
\mathbf{E}(p)\geq h(p)= \frac{2(n- 1)}{\frac{R}{r}- 3}\varphi\Big(\frac{d(p, x_0)}{r}\Big)\cdot r> \frac{2(n- 1)}{\frac{R}{r}- 3}\varphi(1)r> \frac{1}{n}\frac{r^2}{R} ,\nonumber 
\end{align}
which is contradicting the assumption on the excess, the conclusion is proved.
}
\qed

The following lemma provides the existence of good cut-off function on manifolds with $Rc\geq 0$, which will be used later.

\begin{lemma}\label{lem good cut-off function}
{If $Rc(M^n)\geq 0$ and $p\in M^n$, for any $\tau\in (0, 1)$, there is a nonnegative smooth function $\phi: M^n\rightarrow [0, 1]$ \begin{equation}\nonumber 
\phi(x)= \left\{
\begin{array}{rl}
&1 \quad \quad \quad \quad \quad \quad \quad \quad x\in B_{\tau r}(p) \\
&0 \quad \quad \quad \quad \quad \quad \quad \quad x\notin B_{r}(p) \nonumber  
\end{array} \right.
\end{equation}
satisfying $\sup\limits_{x\in B_r(p)} |\Delta \phi(x)|\leq \frac{10^{15}n^{10}}{\tau^{2n}(1- \tau)^8}r^{-2}$.
}
\end{lemma}

\pf
{Scaling $(M^n, g)$ to $(M^n, r^{-2}g)$, then we only to construct $\phi(x)$ for $r= 1$. We define $h_1(\rho)= \int_{\rho}^{1} \int_t^{1} (\frac{s}{t})^{n- 1} ds dt$ and 
\begin{align}
W_1(x)&= \frac{h_1\big(d(p, x)\big)}{h_1(\tau )} \ , \quad \quad \quad \quad h_2(\rho)= \int_0^{\rho} \int_0^t (\frac{s}{t})^{n- 1} ds dt . \nonumber
\end{align}

It is easy to see 
\begin{align}
\frac{1}{2n}(1- \tau)^2\leq h_1(\tau)\leq \frac{\tau^{2- n}}{n(n- 2)} . \label{h_1 bound}
\end{align}

Using Laplace Comparison Theorem, it is straightforward to get $\Delta W_1\geq \frac{1}{h_1(\tau )}$. From the theory of elliptic equations of second order (c.f. \cite{GT}), we can define the function $W: \overline{B_{1}(p)}- B_{\tau }(p)\rightarrow \mathbb{R}$ satisfying
\begin{equation}\nonumber 
\left\{
\begin{array}{rl}
&\Delta W= \frac{1}{h_1(\tau )} \quad \quad \quad \quad \quad on \ B_{1}(p)- \overline{B_{\tau }(p)} \\
&W= 1 \quad \quad \quad \quad \quad \quad \quad on \ \partial B_{\tau }(p) \nonumber  \\
&W= 0 \quad \quad \quad \quad \quad \quad \quad on \ \partial B_{1}(p) \nonumber  \\
\end{array} \right.
\end{equation}

Apply the Maximum Principle to $W- W_1$ on $\overline{B_{1}(p)}- B_{\tau }(p)$, we get 
\begin{align}
W\geq W_1 \quad \quad \quad \quad \quad on\ \overline{B_{1}(p)}- B_{\tau }(p) , \label{W-lower-bound}
\end{align}
which implies that $W\geq 0$.

For any $x_0\in \partial B_{\frac{1}{2}(1+ \tau)}(p)$, we define $W_2(x)= \frac{h_2\big(d(x_0, x)\big)}{h_1(\tau )}$. Note $B_{\frac{1}{2}(1- \tau)}(x_0)\subset \Big(\overline{B_{1}(p)}- B_{\tau }(p)\Big)$, from Laplace Comparison Theorem again, 
\begin{align}
\Delta (W- W_2)\geq 0 \quad \quad \quad \quad \quad on\ B_{\frac{1}{2}(1- \tau)}(x_0) .\nonumber 
\end{align}
Then Maximum Principle yields
\begin{align}
W(x_0)- W_2(x_0)\leq \max_{\partial B_{\frac{1}{2}(1- \tau)}(x_0)} (W- W_2) . \nonumber 
\end{align}

From the definition of $W$ and Maximum Principle, we know 
\begin{align}
W\leq 1 \quad \quad \quad \quad \quad on\ \overline{B_{1}(p)}- B_{\tau }(p) .\nonumber 
\end{align}

On $\partial B_{\frac{1}{2}(1- \tau)}(x_0)$, we have $W_2(x)= \frac{h_2(\frac{1}{2}(1- \tau))}{h_1(\tau )}$. Also note $W_2(x_0)= 0$, then we get 
\begin{align}
W(x_0)\leq 1- \frac{h_2(\frac{1}{2}(1- \tau))}{h_1(\tau )}\ , \quad \quad \quad \quad \forall x_0\in \partial B_{\frac{1}{2}(1+ \tau)}(p) .\nonumber 
\end{align}

By $\Delta W(x)\geq 0$, apply Maximum Principle to $W$ on $\overline{B_{1}(p)}- B_{\frac{1}{2}(1+ \tau)}(p)$, we have 
\begin{align}
W\leq 1- \frac{h_2(\frac{1}{2}(1- \tau))}{h_1(\tau )} \ , \quad \quad \quad \quad on\ \overline{B_{1}(p)}- B_{\frac{1}{2}(1+ \tau)}(p) .\label{W-upper-bound}
\end{align}

We define $\delta_0= \frac{(1- \tau)^2\tau}{16n^2}$, then we can get  
\begin{align}
\frac{h_1(\tau + \delta_0)}{h_1(\tau )}- \Big(1-\frac{h_2(\frac{1}{2}(1- \tau))}{h_1(\tau )}\Big)\geq \frac{(1- \tau)^2}{32n\cdot h_1(\tau)} . \nonumber 
\end{align}
We can find a smooth function $f: [0, 1]\rightarrow [0, 1]$ as the following:
\begin{equation}\nonumber 
f(s)= \left\{
\begin{array}{rl}
&0 \quad \quad \quad \quad \quad 0\leq s\leq 1-\frac{h_2(\frac{1}{2}(1- \tau))}{h_1(\tau )}  \\
&1 \quad \quad \quad \quad \quad \quad \quad \frac{h_1(\tau + \delta_0)}{h_1(\tau )}\leq s\leq 1 ,\nonumber 
\end{array} \right.
\end{equation}
which satisfies  
\begin{align}
\big|f'(s)\big|+ \big|f''(s)\big|\leq \frac{10^4n^2}{(1- \tau)^4}\big[h_1(\tau)^2+ 1\big] .\label{upper bound of f's derivatives}
\end{align}

Note when $\tau < d(x, p)\leq \tau + \delta_0$, from (\ref{W-lower-bound}), we have $W(x)\geq W_1(x)\geq \frac{h_1(\tau + \delta_0)}{h_1(\tau )}$, then $f\big(W(x)\big)= 1$. And when $\frac{1}{2}(1+ \tau)\leq d(x, p)< 1$, from (\ref{W-upper-bound}), we get $W(x)\leq 1- \frac{h_2(\frac{1}{2}(1- \tau))}{h_1(\tau )}$, hence $f\big(W(x)\big)= 0$.

Now we can define smooth function $\phi$ as the following,
\begin{equation}\nonumber 
\phi(x)= \left\{
\begin{array}{rl}
&1 \quad \quad \quad \quad \quad \quad \quad \quad x\in B_{\tau }(p) \\
&0 \quad \quad \quad \quad \quad \quad \quad \quad x\notin B_{1}(p) \nonumber  \\
&f\big(W(x)\big) \quad \quad \quad \quad \quad x\in \overline{B_{1}(p)}- B_{\tau }(p) \nonumber
\end{array} \right.
\end{equation}

From $0\leq W\leq 1$ and Theorem \ref{thm Cheng-Yau's lemma}, for any $x$ satisfying $\tau + \delta_0\leq d(x, p)\leq \frac{1}{2}(1+ \tau)$,  
\begin{align}
\big|\nabla W(x)\big|\leq \frac{200n(\frac{1}{2}\delta_0+ 1)}{\frac{1}{2}\delta_0}\Big[\sup_{B_{\delta_0}(x)}|W|+ \frac{1}{h_1(\tau )}\Big] \leq \frac{800n\big(1+ h_1(\tau)\big)}{\delta_0\cdot h_1(\tau)} . \nonumber 
\end{align}
Then from (\ref{upper bound of f's derivatives}) and (\ref{h_1 bound}), on $M^n$, 
\begin{align}
|\Delta \phi|\leq |f''|\cdot |\nabla W|^2+ |f'|\cdot |\Delta W|\leq 10^{10}n^4\cdot \frac{\big(1+ h_1(\tau)\big)^4}{\delta_0^2\cdot h_1(\tau)^2} \leq \frac{10^{15}n^{10}}{\tau^{2n}(1- \tau)^8} .\nonumber 
\end{align}
}
\qed

\begin{lemma}\label{lem gradient est imply Hessian est}
{On complete Riemannian manifold $M^n$ with $Rc\geq 0$, for  harmonic function $\mathbf{b}$ defined on $B_r(p)$ and any $\tau\in (0, 1)$, we have  
\begin{align}
\fint_{B_{\tau r}(p)} \big|\nabla ^2 \mathbf{b}\big| \leq \frac{10^8n^5r^{-1}}{\tau^{\frac{3n}{2}}(1- \tau)^4}\cdot \sqrt{\sup_{B_r(p)}|\nabla \mathbf{b}|+ 1}\cdot \sqrt{\fint_{B_r(p)} \big||\nabla \mathbf{b}|- 1\big|} .\nonumber 
\end{align}
}
\end{lemma}

\pf
{From Bochner formula and $\Delta \mathbf{b}= 0$, we have 
\begin{align}
\frac{1}{2}\Delta |\nabla \mathbf{b}|^2= |\nabla^2 \mathbf{b}|^2+ Rc(\nabla \mathbf{b}, \nabla \mathbf{b})\geq |\nabla^2 \mathbf{b}|^2 . \nonumber 
\end{align}

From Lemma \ref{lem good cut-off function}, one can choose a nonnegative cut-off function $\phi$ such that 
\begin{equation}\nonumber 
\phi(x)= \left\{
\begin{array}{rl}
&1 \quad \quad \quad \quad \quad \quad \quad \quad x\in B_{\tau r}(p) \\
&0 \quad \quad \quad \quad \quad \quad \quad \quad x\notin B_{r}(p) \nonumber  
\end{array} \right.
\end{equation}
satisfying $\sup_{x\in B_r(p)} |\Delta \phi(x)|\leq \frac{10^{15}n^{10}}{\tau^{2n}(1- \tau)^8}r^{-2}$. Now we have  
\begin{align}
\int_{B_{\tau r}(p)} |\nabla^2 \mathbf{b}|^2&\leq \int_{B_{r}(p)} |\nabla^2 \mathbf{b}|^2\cdot \phi \leq \frac{1}{2}\int_{B_{r}(p)} \Delta \big(|\nabla \mathbf{b}|^2\big)\cdot \phi  \nonumber \\
&= \frac{1}{2} \int_{B_{r}(p)}  |\nabla \mathbf{b}|^2 \cdot \Delta\phi= \frac{1}{2} \int_{B_{r}(p)}  \big(|\nabla \mathbf{b}|^2- 1\big) \cdot \Delta\phi \nonumber \\
&\leq \frac{1}{2}\sup_{B_r(p)} |\Delta \phi| \cdot \int_{B_{r}(p)}   \Big|\big|\nabla \mathbf{b}\big|^2- 1\Big|\nonumber \\
&\leq \frac{\tau^{-n}}{2}\sup_{B_r(p)} |\Delta \phi| \cdot \big(\sup_{B_r(p)}|\nabla \mathbf{b}|+ 1\big)V\big(B_{\tau r}(p)\big)\fint_{B_r(p)} \big||\nabla \mathbf{b}|- 1\big| . \nonumber 
\end{align}

From Cauchy-Schwartz inequality and the above inequality, we have 
\begin{align}
\fint_{B_{\tau r}(p)} \big|\nabla ^2 \mathbf{b}\big| \leq \Big(\fint_{B_{\tau r}(p)} \big|\nabla ^2 \mathbf{b}\big|^2\Big)^{\frac{1}{2}}\leq \frac{10^8n^5r^{-1}}{\tau^{\frac{3n}{2}}(1- \tau)^4}\sqrt{\big(\sup_{B_r(p)}|\nabla \mathbf{b}|+ 1\big)\fint_{B_r(p)} \big||\nabla \mathbf{b}|- 1\big|} .\nonumber 
\end{align}
}
\qed

\begin{theorem}[Li-Yau]\label{thm Li-Yau harnack}
{Let $(M^n, g)$ be a complete Riemannian manifold with $Rc\geq 0$, then the heat kernel $H(x, y, t)$ satisfies
\begin{align}
H(x, y, t)\leq (100n)^{2n+ 2}\frac{1}{V\big(B_{\sqrt{t}}(y)\big)}\exp\Big\{-\frac{d^2(x, y)}{100t}\Big\} .\nonumber 
\end{align}
}
\end{theorem}

\pf
{By choosing suitable $\epsilon$ in \cite[Corollary $3.1$]{LY}, the conclusion follows. 
}
\qed

Now we have the following existence result of almost linear harmonic function $\mathbf{b}$ with respect to the local Busemann function $b^+$.

\begin{lemma}\label{lem existence of harmonic function-r}
{On complete Riemannian manifold $M^n$ with $Rc\geq 0$, assume that $[q^+, q^-, p]$ is an AG-triple with the excess $\leq \frac{4}{n}\frac{r^2}{R}$ and the scale $\geq R$, also assume $R\geq 2^{2n+ 1}r$. Then there exists harmonic function $\mathbf{b}$ defined on $B_{2r}(p)$ such that 
\begin{align}
\sup_{B_{r}(p)} |\nabla \mathbf{b}|\leq 1+ 2^{51n^2}\big(\frac{r}{R}\big)^{\frac{1}{4(n- 1)}} \quad \quad and \quad \quad 
\fint_{B_{r}(p)} \big|\nabla (\mathbf{b}- b^+)\big| \leq 2^{4n}\big(\frac{r}{R}\big)^{\frac{1}{2(n- 1)}} , \nonumber 
\end{align}
where $b^{+}(x)= d(x, q^{+})- d(p, q^+)$.
}
\end{lemma}

\begin{remark}
{The bound of $|\nabla \mathbf{b}|$ was only a uniform bound $C(n)$ in \cite{CC-Ann}, it was observed in \cite{ChN} that $|\nabla \mathbf{b}|$ can have the improved bound $1+ \epsilon$. The argument to get $|\nabla \mathbf{b}|\leq 1+ 2^{51n^2}\big(\frac{r}{R}\big)^{\frac{1}{4(n- 1)}}$, comes from $(3.42)- (3.45)$ of \cite{ChN}, which was suggested to us by A. Naber.
}
\end{remark}

\pf
{\textbf{Step (1)}. We define the function $\mathbf{b}$ by
\begin{equation}\nonumber 
\left\{
\begin{array}{rl}
&\Delta \mathbf{b}= 0 \quad \quad \quad \quad \quad on \ B_{2r}(p) \\
&\mathbf{b}= b^+ \quad \quad \quad \quad \quad on \ \partial B_{2r}(p) . \nonumber  
\end{array} \right.
\end{equation}

Define $\psi(\rho)= \int_{\rho}^4 \int_t^4 \big(\frac{s}{t}\big)^{n- 1} ds dt$, choose $q\in B_{6r}(p)- \overline{B_{4r}(p)}$, and define the function $\tilde{h}(x): A_{2r, 8r}(q)\rightarrow [0, \infty)$ by $\tilde{h}(x)= \frac{2(n- 1)r}{\frac{R}{2r}- 7}\psi\Big(\frac{d(x, q)}{2r}\Big)$. From Laplace Comparison Theorem, it is straightforward to get 
\begin{align}
\min_{x\in A_{2r, 8r}(q)}\Delta (\tilde{h}- b^+)(x)\geq 0 .\nonumber 
\end{align}
Note $B_{2r}(p)\subset A_{2r, 8r}(q)$, hence $\min_{x\in B_{2r}(p)}\Delta (\mathbf{b}- b^++ \tilde{h})(x)\geq 0$.

From Maximum principle, for any $x\in B_{2r}(p)$, we get 
\begin{align}
\big(\mathbf{b}- b^+\big)(x)&\leq \big(\mathbf{b}- b^++ \tilde{h}\big)(x)\leq \max_{\partial B_{2r}(p)} \big(\mathbf{b}- b^++ \tilde{h}\big)= \max_{\partial B_{2r}(p)} \tilde{h}\leq \frac{n- 1}{\frac{R}{2r}- 7}\psi(1)\cdot 2r \nonumber \\
&\leq 4^n\cdot \frac{2r}{R}\cdot 2r\leq 2^{2n+ 2}\frac{r^2}{R} .\label{gag 0.1}
\end{align}

Let $b^{-}(x)= d(x, q^-)- d(p, q^-)$, then note $\mathbf{E}(x)= b^+(x)+ b^-(x)+ \mathbf{E}(p)$, then from $R\geq 2^{2n+ 1}r$ and 
Lemma \ref{lem general Abresch-Gromoll}, $\sup_{B_{2r}(p)}\mathbf{E}(x)\leq 2^8 \big(\frac{r}{R}\big)^{\frac{1}{n- 1}}r$, and we have 
\begin{align}
\mathbf{b}- b^+= \mathbf{b}+ b^-- \mathbf{E}+ \mathbf{E}(p)\geq \mathbf{b}+ b^-- 2^8 \big(\frac{r}{R}\big)^{\frac{1}{n- 1}}r\ , \quad \quad \quad \quad on\ B_{2r}(p) . \label{gag 1}
\end{align}

From Laplace Comparison Theorem again, 
\begin{align}
\Delta \Big(\mathbf{b}+ b^-- 2^8 \big(\frac{r}{R}\big)^{\frac{1}{n- 1}}r- \tilde{h}\Big)= \Delta \big(b^-- \tilde{h}\big)\leq \frac{n- 1}{d(x, q^-)}- \frac{n- 1}{R- 14r}\leq 0 \ ,  \quad \quad on \ B_{2r}(p) .\nonumber 
\end{align}

By Maximum Principle, note $\mathbf{E}(x)\geq 0$ and the assumption on the excess, we get 
\begin{align}
&\quad \min_{x\in B_{2r}(p)}(\mathbf{b}+ b^-- 2^8 \big(\frac{r}{R}\big)^{\frac{1}{n- 1}}r- \tilde{h})(x) \nonumber \\
&\geq \min_{\partial B_{2r}(p)} (\mathbf{b}+ b^-- 2^8 \big(\frac{r}{R}\big)^{\frac{1}{n- 1}}r- \tilde{h}) \nonumber \\
&\geq \min_{\partial B_{2r}(p)} (b^++ b^-)- 2^8 \big(\frac{r}{R}\big)^{\frac{1}{n- 1}}r- \frac{n- 1}{\frac{R}{2r}- 7}\psi(1)\cdot 2r \nonumber \\
&= \min_{\partial B_{2r}(p)} \big(\mathbf{E}(x)- \mathbf{E}(p)\big)- 2^8 \big(\frac{r}{R}\big)^{\frac{1}{n- 1}}r- \frac{n- 1}{\frac{R}{2r}- 7}\psi(1)\cdot 2r \nonumber \\
&\geq - \frac{1}{n}\frac{4r^2}{R}- 2^8 \big(\frac{r}{R}\big)^{\frac{1}{n- 1}}r- \frac{n- 1}{\frac{R}{2r}- 7}\psi(1)\cdot 2r  \nonumber  \\
&\geq -2^{3n}\big(\frac{r}{R}\big)^{\frac{1}{n- 1}}r .\nonumber 
\end{align}

From (\ref{gag 1}) and $\tilde{h}\geq 0$, we obtain
\begin{align}
\min_{x\in B_{2r}(p)}(\mathbf{b}- b^+)\geq \mathbf{b}+ b^-- 2^8 \big(\frac{r}{R}\big)^{\frac{1}{n- 1}}r- \tilde{h} \geq -2^{3n}\big(\frac{r}{R}\big)^{\frac{1}{n- 1}}r  .\label{gag 0.2}
\end{align}

By (\ref{gag 0.1}) and (\ref{gag 0.2}), it yields 
\begin{align}
\sup_{B_{2r}(p)}\big|\mathbf{b}- b^+\big|\leq 2^{3n}\big(\frac{r}{R}\big)^{\frac{1}{n- 1}}r ,\label{require 2.8.0} 
\end{align}
which implies
\begin{align}
\sup_{B_{2r}(p)}|\mathbf{b}|\leq 2^{3n}\big(\frac{r}{R}\big)^{\frac{1}{n- 1}}r+ 2r .\label{gag c1} 
\end{align}

From Theorem \ref{thm Cheng-Yau's lemma} and (\ref{gag c1}), note $\frac{R}{2r}\geq 1$, we get 
\begin{align}
\sup_{B_{\frac{15}{8}r}(p)}|\nabla \mathbf{b}|\leq \frac{60n}{\frac{1}{16}r}\sup_{B_{2r}(p)} |\mathbf{b}|\leq 2^{4n+ 12} .\label{gradient of b has a bound} 
\end{align}

Note $d(p, q^+)\geq R> 4r$, then from Lemma \ref{lem estimate of Lap b},  $\fint_{B_{2r}(p)} |\Delta b^{+}|\leq \frac{3n}{2r}$. Now do integration by parts, from (\ref{require 2.8.0}), we get 
\begin{align}
\fint_{B_{2r}(p)} \big|\nabla (\mathbf{b}- b^+)\big|^2& = \fint_{B_{2r}(p)} \Delta b^+\cdot (\mathbf{b}- b^+) \leq \sup_{B_{2r}(p)} \big|\mathbf{b}- b^+\big|\cdot \fint_{B_{2r}(p)} |\Delta b^+|  \nonumber \\
& \leq 2^{3n}\big(\frac{r}{R}\big)^{\frac{1}{n- 1}}\cdot (3n)\leq 2^{7n}\big(\frac{r}{R}\big)^{\frac{1}{n- 1}} .\label{2r ball aver integral}
\end{align}

From (\ref{2r ball aver integral}) and the Bishop-Gromov Comparison Theorem, we have
\begin{align}
\fint_{B_{r}(p)} \big|\nabla (\mathbf{b}- b^+)\big|&\leq \Big(\fint_{B_{r}(p)} \big|\nabla (\mathbf{b}- b^+)\big|^2\Big)^{\frac{1}{2}}\leq \Big[2^n\fint_{B_{2r}(p)} \big|\nabla (\mathbf{b}- b^+)\big|^2\Big]^{\frac{1}{2}}\leq 2^{4n}\big(\frac{r}{R}\big)^{\frac{1}{2(n- 1)}} .\nonumber 
\end{align} 

\textbf{Step (2)}. Also (\ref{2r ball aver integral}) implies the following
\begin{align}
\fint_{B_{\frac{7}{4}r}(p)} \big||\nabla \mathbf{b}|- 1\big|&\leq \fint_{B_{\frac{7}{4}r}(p)} \big|\nabla \mathbf{b}- \nabla b^+\big|\leq (\frac{8}{7})^n \fint_{B_{2r}(p)} \big|\nabla \mathbf{b}- \nabla b^+\big| \nonumber \\
&\leq (\frac{8}{7})^n\big(\fint_{B_{2r}(p)} \big|\nabla (\mathbf{b}- b^+)\big|^2\big)^{\frac{1}{2}}\leq 2^{7n}\big(\frac{r}{R}\big)^{\frac{1}{2(n- 1)}} \label{almost linear ineq-1} \\
\fint_{B_{\frac{15}{8}r}(p)} \big||\nabla \mathbf{b}|- 1\big|&\leq  2^{7n}\big(\frac{r}{R}\big)^{\frac{1}{2(n- 1)}} .\label{almost linear ineq-2}
\end{align}

From (\ref{gradient of b has a bound}), (\ref{almost linear ineq-2}) and apply Lemma \ref{lem gradient est imply Hessian est} for $\tau= \frac{14}{15}$ there, we have 
\begin{align}
\fint_{B_{\frac{7}{4} r}(p)}|\nabla^2 \mathbf{b}|& \leq 2^{7n}10^8n^5r^{-1}\sqrt{\sup_{B_{\frac{15}{8}r}(p)}\big(|\nabla \mathbf{b}|+ 1\big)\cdot \fint_{B_{\frac{15}{8}r}(p)} \big||\nabla \mathbf{b}|- 1\big|} \nonumber \\
&\leq 2^{27n}r^{-1}\cdot \sqrt{2^{20n}\cdot \big(\frac{r}{R}\big)^{\frac{1}{2(n- 1)}}}\leq 2^{37n}r^{-1}\big(\frac{r}{R}\big)^{\frac{1}{4(n- 1)}} . \label{hessian of b upper bound} 
\end{align}

From the Bochner formula and $\Delta \mathbf{b}= 0$, using $Rc\geq 0$, we have $\Delta |\nabla \mathbf{b}|^2\geq 2|\nabla^2 \mathbf{b}|^2$, combining
\begin{align}
|\nabla^2 \mathbf{b}|\geq \big|\nabla |\nabla \mathbf{b}|\big| , \label{2nd deri bound cont 1st deri}
\end{align}
we get  
\begin{align}
\Delta |\nabla \mathbf{b}|\geq \frac{|\nabla^2 \mathbf{b}|^2- \big|\nabla |\nabla \mathbf{b}|\big|^2}{|\nabla \mathbf{b}|}\geq 0 .\label{gradient Lap is nonneg}
\end{align}

Let $\phi\in C^{\infty}(M^n)$ be a nonnegative cut-off function such that 
\begin{equation}\nonumber 
\phi(x)= \left\{
\begin{array}{rl}
&1 \quad \quad \quad \quad \quad x\in B_{\frac{3}{2}r}(p) \\
&0 \quad \quad \quad \quad \quad x\notin B_{\frac{7}{4}
r}(p)  ,
\end{array} \right.
\end{equation}
and 
\begin{align}
-\frac{8}{r}\leq \phi '\leq 0\ , \quad \quad \quad \quad |\phi ''|\leq \frac{64}{r^2} .\label{eq bound of cut-off}
\end{align}

Now for any $y\in B_r(p)$, from (\ref{2nd deri bound cont 1st deri}), (\ref{gradient Lap is nonneg}) and (\ref{eq bound of cut-off}), we have 
\begin{align}
&\quad \frac{d}{dt}\int_{B_{2r}(p)} \big(|\nabla \mathbf{b}|- 1\big)\phi\cdot H(x, y, t)dx= \int_{B_{2r}(p)} \big(|\nabla \mathbf{b}|- 1\big)\phi\cdot \Delta_x H(x, y, t)dx \nonumber \\
&= \int_{B_{2r}(p)} \Big\{ \phi\cdot \Delta\big(|\nabla \mathbf{b}|- 1\big)+ \big(|\nabla \mathbf{b}|- 1\big)\cdot \Delta\phi+ 2\nabla |\nabla \mathbf{b}|\cdot \nabla \phi \Big\} H(x, y, t)dx \nonumber \\
&\geq -\int_{A_{\frac{3r}{2}, \frac{7}{4}r}(p)} \Big[\big||\nabla \mathbf{b}|- 1\big|\cdot \frac{64}{r^2}+ 2|\nabla^2 \mathbf{b}|\cdot \frac{8}{r}\Big]H(x, y, t)dx . \label{evolution lower bound} 
\end{align}

From Theorem \ref{thm Li-Yau harnack}, for $x\in A_{\frac{3r}{2}, \frac{7}{4}r}(p), y\in B_r(p)$, we have 
\begin{align}
H(x, y, t)&\leq (100n)^{2n+ 2}\frac{1}{V\big(B_{\sqrt{t}}(y)\big)}\exp\Big\{-\frac{d^2(x, y)}{100t}\Big\} \nonumber \\
&\leq (100n)^{2n+ 2}\frac{1}{V\big(B_{3r}(y)\big)}\cdot \Big(\frac{3r}{\sqrt{t}}\Big)^n\exp\Big\{-\frac{d^2(x, y)}{100t}\Big\} \nonumber \\
&\leq (300n)^{2n+ 2}\frac{1}{V\big(B_{2r}(p)\big)}\Big(\frac{r}{\sqrt{t}}\Big)^n\exp\Big\{-\frac{r^2}{400t}\Big\} \nonumber \\
&\leq (300n)^{2n+ 2}\frac{1}{V\big(B_{2r}(p)\big)}\cdot (200n)^{\frac{n}{2}}e^{- \frac{n}{2}} \nonumber \\
&\leq (300n)^{4n}\frac{1}{V\big(B_{2r}(p)\big)} .\label{heat kernel upper bound}
\end{align}

From (\ref{evolution lower bound}), (\ref{heat kernel upper bound}), (\ref{almost linear ineq-1}) and (\ref{hessian of b upper bound}), we get 
\begin{align}
&\quad \frac{d}{dt}\int_{B_{2r}(p)} \big(|\nabla \mathbf{b}|- 1\big)\phi\cdot H(x, y, t)dx \nonumber \\
&\geq -(400n)^{4n}\fint_{B_{\frac{7}{4} r}(p)} \Big[\big||\nabla \mathbf{b}|- 1\big|\cdot \frac{1}{r^2}+ |\nabla^2 \mathbf{b}|\cdot \frac{1}{r}\Big] \nonumber \\
&\geq -(400n)^{4n}\Big[2^{7n}\big(\frac{r}{R}\big)^{\frac{1}{2(n- 1)}}\frac{1}{r^2}+ 2^{37n}r^{-2}\big(\frac{r}{R}\big)^{\frac{1}{4(n- 1)}}\Big] \nonumber \\
&\geq -\frac{2^{50n^2}}{r^2}\cdot \big(\frac{r}{R}\big)^{\frac{1}{4(n- 1)}} . \label{diff ineq-1}
\end{align}

Take the integration of (\ref{diff ineq-1}), from (\ref{almost linear ineq-2}) and (\ref{heat kernel upper bound}), we get 
\begin{align}
|\nabla \mathbf{b}|(y)- 1&\leq \int_{B_{2r}(p)} \big(|\nabla \mathbf{b}|- 1\big)\phi\cdot H(x, y, r^2)dx+ 2^{50n^2}\cdot \big(\frac{r}{R}\big)^{\frac{1}{4(n- 1)}} \nonumber \\
&\leq (300n)^{4n}\fint_{B_{\frac{7r}{4}}(p)} \big||\nabla \mathbf{b}|- 1\big|+ 2^{50n^2}\cdot \big(\frac{r}{R}\big)^{\frac{1}{4(n- 1)}} \nonumber \\
&\leq 2^{51n^2}\big(\frac{r}{R}\big)^{\frac{1}{4(n- 1)}} .\nonumber 
\end{align}
}
\qed

On Riemannian manifolds, if there is a segment $\gamma_{p, q}$ between two points $p, q$, we can choose the middle point of the segment $\gamma_{p, q}$, denoted as $z$. Then $[p, q, z]$ is an AG-triple with the excess $0$ and the scale $\frac{1}{2}d(p, q)$. 

For a metric space $\mathbf{X}_k$, generally we can not find the middle point as in Riemannian manifolds. However, if there exists a suitable Gromov-Hausdorff approximation from $\mathbb{R}^k\times \mathbf{X}_k$ to manifold $M^n$ locally, the following lemma provides the existence of almost middle point and AG-triple in metric space $\mathbf{X}_k$.

\begin{lemma}\label{lem almost mid pt}
{For $0\leq k\leq n, 0< \delta< \frac{1}{3}$, $B_{10r}(0, \hat{q})\subset \mathbb{R}^k\times \mathbf{X}_k$ and $B_{10r}(q)\subset (M^n, g)$, if there is an $(\delta r)$-Gromov-Hausdorff approximation
\begin{align}
\Phi: B_{10r}(0, \hat{q})\rightarrow B_{10r}(q) , \nonumber 
\end{align}
then for $\hat{q}_1^+, \hat{q}_1^-\in B_r(\hat{q})\subset \mathbf{X}_k$, there exists $\hat{q}_0 \in B_{3r}(\hat{q})$ such that 
\begin{align}
\big|d(\hat{q}_0, \hat{q}_1^+)- \frac{1}{2}d(\hat{q}_1^+, \hat{q}_1^-)\big|\leq 8\sqrt{\delta}r \quad \quad and \quad \quad \big|d(\hat{q}_0, \hat{q}_1^-)- \frac{1}{2}d(\hat{q}_1^+, \hat{q}_1^-)\big|\leq 8\sqrt{\delta}r . \label{almost mid ineq conclusion} 
\end{align}
And if $\sqrt{\delta}r\leq \frac{d(\hat{q}_1^+, \hat{q}_1^-)}{100}$, then $\big[p_{k+1}^+, p_{k+ 1}^-, q_1\big]$ is an AG-triple with the excess $\leq 16\sqrt{\delta}r$ and the scale $\geq \frac{1}{4}d(\hat{q}_1^+, \hat{q}_1^-)$, where  
\begin{align}
&q_1= \Phi(0, \hat{q}_0), \quad \quad \quad p_{k+ 1}^+= \Phi(0, \hat{q}_1^+), \quad \quad \quad p_{k+ 1}^-= \Phi(0, \hat{q}_1^-) . \nonumber 
\end{align}
}
\end{lemma}

\pf
{From the assumption that $\Phi$ is an $(\delta r)$-Gromov-Hausdorff approximation,
\begin{align}
\big|d(p_{k+ 1}^+, p_{k+ 1}^-)- d(\hat{q}_1^+, \hat{q}_1^-)\big|= \big|d\big(\Phi(0, \hat{q}_1^+), \Phi(0, \hat{q}_1^-)\big)- d\big((0, \hat{q}_1^+), (0, \hat{q}_1^-)\big)\big|< \delta r  .\label{length of the seg on subspace}
\end{align}

It is easy to see that $p_{k+ 1}^+, p_{k+ 1}^-\in B_{r+ \delta r}(q)$. In fact the segment $\gamma_{p_{k+ 1}^+, p_{k+ 1}^-}\subset B_{3r}(q)$, otherwise, note $d(\hat{q}_1^+, \hat{q}_1^-)\leq 2r$, we will get 
\begin{align}
d(p_{k+ 1}^+, p_{k+ 1}^-)\geq 2\big[3r- (r+ \delta r)\big]> d(\hat{q}_1^+, \hat{q}_1^-)+ \delta r ,\nonumber 
\end{align}
which is contradicting (\ref{length of the seg on subspace}).
 
Then we can choose the middle point of $\gamma_{p_{k+ 1}^+, p_{k+ 1}^-}$, denoted as $z_0\in B_{3r}(q)$, then there exists $(s_0, \hat{q}_0)\in B_{10r}(0, \hat{q})$ such that 
\begin{align}
d\big(\Phi(s_0, \hat{q}_0), z_0\big)< \delta r .\label{middle pt is appro by image}
\end{align}

Now using $\frac{1}{2}d(p_{k+1}^+, p_{k+ 1}^-)= d(p_{k+1}^+, z_0)$ and (\ref{middle pt is appro by image}), we have 
\begin{align}
&\quad \Big|d\big((s_0, \hat{q}_0), (0, \hat{q}_1^+)\big)- \frac{1}{2}d(p_{k+1}^+, p_{k+ 1}^-)\Big| \nonumber \\
&\leq \Big|d\big((s_0, \hat{q}_0), (0, \hat{q}_1^+)\big)- d\big(\Phi(s_0, \hat{q}_0), p_{k+ 1}^+\big)\Big|\nonumber \\
&\quad + \Big|d\big(\Phi(s_0, \hat{q}_0), p_{k+ 1}^+\big)- d(p_{k+1}^+, z_0)\Big| \nonumber \\
&< \delta r+ d\big(\Phi(s_0, \hat{q}_0), z_0\big) < 2\delta r .\label{scale bound-1}
\end{align}
Similarly, by using $\frac{1}{2}d(p_{k+1}^+, p_{k+ 1}^-)= d(p_{k+1}^-, z_0)$, we have 
\begin{align}
\Big|d\big((s_0, \hat{q}_0), (0, \hat{q}_1^-)\big)- \frac{1}{2}d(p_{k+1}^+, p_{k+ 1}^-)\Big| < 2\delta r .\label{scale bound-2}
\end{align}

From (\ref{length of the seg on subspace}), (\ref{scale bound-1}) and (\ref{scale bound-2}), 
\begin{align}
d(\hat{q}_0, \hat{q}_1^+)^2+ s_0^2&< \big(\frac{1}{2}d(p_{k+1}^+, p_{k+ 1}^-)+ 2\delta r\big)^2< \big(\frac{1}{2}d(\hat{q}_1^+, \hat{q}_1^-)+ \frac{5}{2}\delta r \big)^2 \nonumber \\
d(\hat{q}_0, \hat{q}_1^-)^2+ s_0^2&< \big(\frac{1}{2}d(\hat{q}_1^+, \hat{q}_1^-)+ \frac{5}{2}\delta r \big)^2 .\nonumber
\end{align}
Take the sum of the above two inequalities, we have 
\begin{align}
2\cdot \big(\frac{1}{2}d(\hat{q}_1^+, \hat{q}_1^-)+ \frac{5}{2}\delta r \big)^2&> 2s_0^2+ d(\hat{q}_0, \hat{q}_1^+)^2+ d(\hat{q}_0, \hat{q}_1^-)^2 \nonumber \\
&\geq 2s_0^2+ \frac{1}{2}\big[d(\hat{q}_0, \hat{q}_1^+)+ d(\hat{q}_0, \hat{q}_1^-)\big]^2 \nonumber \\
&\geq 2s_0^2+ \frac{1}{2}d(\hat{q}_1^+, \hat{q}_1^-)^2= \frac{1}{2}d(\hat{q}_1^+, \hat{q}_1^-)^2+ 2s_0^2 . \nonumber 
\end{align}
Simplify the above inequality, using $d(\hat{q}_1^+, \hat{q}_1^-)\leq 2r$, we get 
\begin{align}
|s_0|\leq \frac{5}{2}\delta r+ \sqrt{5}\sqrt{\delta}r\leq 5\sqrt{\delta}r . \label{s_0 bound}
\end{align}

From (\ref{scale bound-1}), (\ref{scale bound-2}) and (\ref{s_0 bound}), we have 
\begin{align}
\Big|d\big((0, \hat{q}_0), (0, \hat{q}_1^+)\big)- \frac{1}{2}d(p_{k+1}^+, p_{k+ 1}^-)\Big| < 7\sqrt{\delta} r \label{new scale-1} \\
\Big|d\big((0, \hat{q}_0), (0, \hat{q}_1^-)\big)- \frac{1}{2}d(p_{k+1}^+, p_{k+ 1}^-)\Big| < 7\sqrt{\delta} r .\label{new scale-2}
\end{align}

From (\ref{length of the seg on subspace}), (\ref{new scale-1}) and (\ref{new scale-2}),  we get (\ref{almost mid ineq conclusion}) and the following
\begin{align}
&\quad d(q_1, p_{k+1}^+)+ d(q_1, p_{k+1}^-)- d(p_{k+1}^+, p_{k+1}^-) \nonumber \\
&\leq d\big((0, \hat{q}_0), (0, \hat{q}_1^+)\big)+ d\big((0, \hat{q}_0), (0, \hat{q}_1^-)\big)- d(p_{k+1}^+, p_{k+ 1}^-)+ 2\delta r \nonumber \\
&< 14\sqrt{\delta}r+ 2\delta r\leq 16\sqrt{\delta}r \nonumber \\
&d(q_1, p_{k+1}^+)\geq d\big((0, \hat{q}_0), (0, \hat{q}_1^+)\big)- \delta r> \frac{1}{2}d(p_{k+1}^+, p_{k+ 1}^-)- 7\sqrt{\delta} r- \delta r \nonumber \\
&\quad\quad   \geq  \frac{1}{2}(d(\hat{q}_1^+, \hat{q}_1^-)- \delta r)- 8\sqrt{\delta}r\geq \frac{1}{4}d(\hat{q}_1^+, \hat{q}_1^-) . \nonumber 
\end{align}
Similarly we have $d(q_1, p_{k+1}^-)\geq \frac{1}{4}d(\hat{q}_1^+, \hat{q}_1^-)$. Hence the scale is $\geq \frac{1}{4}d(\hat{q}_1^+, \hat{q}_1^-)$.

Finally from (\ref{almost mid ineq conclusion}), we get 
\begin{align}
d(\hat{q}_0, \hat{q})\leq d(\hat{q}_0, \hat{q}_1^+)+ d(\hat{q}, \hat{q}_1^+)\leq \frac{1}{2}d(\hat{q}_1^+, \hat{q}_1^-)+ 8\sqrt{\delta}r+ r< 3r , \nonumber 
\end{align}
which implies $\hat{q}_0\in B_{3r}(\hat{q})$.
}
\qed

\begin{prop}\label{thm AG-triple imply one more splitting-pre}
{For $B_{10r}(q)\subset (M^n, g)$ with $Rc(g)\geq 0$ and any $\epsilon\in (0, 1)$, any $\frac{r_1}{r}\leq n^{-300n^3}\epsilon^{10n}$ and $\delta\leq n^{-2000n^3}\epsilon^{70n}\big(\frac{r_1}{r}\big)^4$. If there is an $(\delta r)$-Gromov-Hausdorff approximation for $0\leq k\leq n$ and $B_{10r}(0, \hat{q})\subset \mathbb{R}^k\times \mathbf{X}_k$, 
\begin{align}
\Phi: B_{10r}(0, \hat{q})\rightarrow B_{10r}(q) , \nonumber 
\end{align}
where $\mathrm{diam}\big(B_r(\hat{q})\big)= r_0\geq \frac{1}{4}r$. Then there are harmonic functions $\big\{\mathbf{b}_i\big\}_{i= 1}^{k+ 1}$ defined on some geodesic ball $B_{r_1}(q_1)\subset B_{10r}(q)$, such that 
\begin{align}
\sup_{B_{r_1}(q_1)\atop i= 1, \cdots, k+ 1} |\nabla \mathbf{b}_i|\leq 1+ \epsilon \quad \quad and \quad \quad \fint_{B_{r_1}(q_1)}\sum_{i, j= 1}^{k+ 1} \big|\langle \nabla \mathbf{b}_i, \nabla \mathbf{b}_j\rangle- \delta_{ij}\big|\leq \epsilon  .\nonumber 
\end{align}
}
\end{prop}

\pf
{We firstly assume $\delta\leq (400)^{-2}$, then from assumption $\mathrm{diam}\big(B_r(\hat{q})\big)= r_0$ and Lemma \ref{lem almost mid pt}, there are $\hat{q}_1^+, \hat{q}_1^-\in B_{r}(\hat{q})$ and $\hat{q}_0\in B_{3r}(\hat{q})$ such that 
\begin{align}
&d(\hat{q}_1^+, \hat{q}_1^-)= r_0 \nonumber \\
&\big|d(\hat{q}_0, \hat{q}_1^+)- \frac{1}{2}r_0\big|\leq 8\sqrt{\delta}r \quad \quad and \quad \quad \big|d(\hat{q}_0, \hat{q}_1^-)- \frac{1}{2}r_0\big|\leq 8\sqrt{\delta}r .\nonumber 
\end{align}

Then we have 
\begin{align}
\tilde{r}_0&\vcentcolon= d(\hat{q}_0, \hat{q}_1^+)\geq \frac{1}{2}r_0- 8\sqrt{\delta}r\geq \frac{1}{4}r_0\geq \frac{1}{16}r ,\label{lower bound as assp of prop}\\
\tilde{r}_0&\leq \frac{1}{2}r_0+ 8\sqrt{\delta}r\leq 2r .\label{upper bound as assp of prop}
\end{align}

Let $\big\{\mathbf{e}_i\big\}_{i= 1}^k$ be the standard basis for $\mathbb{R}^k$, put $q_1= \Phi(0, \hat{q}_0)$ and
\begin{align}
& p_{k+ 1}^+= \Phi(0, \hat{q}_1^+), \quad \quad \quad b_{k+ 1}^+(\cdot)= d(\cdot, p_{k+ 1}^+)- d(q_1, p_{k+ 1}^+) ,\nonumber \\
& p_i^+= \Phi(\tilde{r}_0\cdot \mathbf{e}_i, \hat{q}_0)\ , \quad \quad \quad b_i^+(\cdot)= d(\cdot, p_i^+)- d(q_1, p_i^+) .\nonumber 
\end{align}

From (\ref{lower bound as assp of prop}), (\ref{upper bound as assp of prop}), let $r_1= \zeta\cdot r$, if 
\begin{align}
\zeta&\leq n^{-110n}\cdot \epsilon_1^{10n} \label{zeta need}\\
\delta&\leq n^{-1250n}\epsilon_1^{100n}\zeta , \label{delta upper bound 1}
\end{align}
then we can apply Proposition \ref{prop approximation component is almost orthogonal} to get  
\begin{align}
\fint_{B_{r_1}(q_1)}\sum_{i, j= 1}^{k+ 1} \big|\langle \nabla b_i^+, \nabla b_j^+\rangle- \delta_{ij}\big|\leq \epsilon_1 .\label{almost o.n.-need-1}
\end{align}

Define $p_{k+ 1}^-= \Phi(0, \hat{q}_1^-),  p_i^-= \Phi(-\tilde{r}_0\cdot \mathbf{e}_i, \hat{q}_0)$, then from Lemma \ref{lem almost mid pt}, $\big[p_{k+1}^+, p_{k+ 1}^-, q_1\big]$ is an AG-triple with the excess $\leq 16\sqrt{\delta}r$ and the scale $\geq \frac{1}{4}r_0$. Because $\Phi$ is an $(\delta r)$-Gromov-Hausdorff approximation, it is easy to show that $[p_i^+, p_i^-, q_1]$ is AG-triple with the excess $\leq 16\sqrt{\delta}r$ and the scale $\geq \frac{1}{4}r_0$ for $i= 1, \cdots, k$. 

If we assume 
\begin{align}
\delta\leq \frac{\zeta^4}{n^2} ,\label{delta upper bound 2}
\end{align}
then $16\sqrt{\delta}r\leq \frac{4}{n}\frac{r_1^2}{\frac{1}{4}r_0}$. Also from $\frac{1}{4}r_0\geq \frac{1}{16}r$ and (\ref{zeta need}), we have $\frac{1}{4}r_0\geq 2^{2n+ 1}r_1$.

Now we can apply Lemma \ref{lem existence of harmonic function-r} to obtain harmonic functions $\big\{\mathbf{b}_i\big\}_{i= 1}^{k+ 1}$ satisfying
\begin{align}
&\sup_{B_{r_1}(q_1)\atop i= 1, \cdots, k+ 1} |\nabla \mathbf{b}_i|\leq 1+ 2^{51n^2}\Big(\frac{r_1}{\frac{1}{4}r_0}\Big)^{\frac{1}{4(n- 1)}}\leq 1+ 2^{51n^2}\epsilon_1^2  \label{gradient bound of b need} \\
&\sup_{i= 1, \cdots, k+ 1}\fint_{B_{r_1}(q_1)} \big|\nabla (\mathbf{b}_i- b_i^+)\big|\leq 2^{4n}\Big(\frac{r_1}{\frac{1}{4}r_0}\Big)^{\frac{1}{2(n- 1)}}\leq 2^{4n}\epsilon_1^5 .\label{harmonic is almost linear}
\end{align}

From (\ref{gradient bound of b need}) and (\ref{harmonic is almost linear}), we get 
\begin{align}
&\quad \fint_{B_{r_1}(q_1)} \big|\langle \nabla \mathbf{b}_i, \nabla\mathbf{b}_j\rangle- \delta_{ij}\big| \nonumber \\
&\leq \fint_{B_{r_1}(q_1)} \big|\nabla (\mathbf{b}_i- b_i^+)\big|\cdot |\nabla \mathbf{b}_j|+ 
 \Big|\big\langle \nabla b_i^+, \nabla(\mathbf{b}_j- b_j^+)\big\rangle\Big|+ \big|\langle \nabla b_i^+, \nabla b_j^+ \rangle- \delta_{ij}\big| \nonumber \\
 &\leq (2+ 2^{51n^2}\epsilon_1^2)\cdot 2^{4n}\epsilon_1^5+ \fint_{B_{r_1}(q_1)} \big|\langle \nabla b_i^+, \nabla b_j^+ \rangle- \delta_{ij}\big| .\nonumber 
\end{align}

From (\ref{almost o.n.-need-1}) and the above inequality, we have
\begin{align}
\fint_{B_{r_1}(q_1)} \sum_{i, j= 1}^{k+ 1}\big|\langle \nabla \mathbf{b}_i, \nabla\mathbf{b}_j\rangle- \delta_{ij}\big|\leq (n+ 1)^2\cdot (2+ 2^{51n^2}\epsilon_1^2)\cdot 2^{4n}\epsilon_1^5+ \epsilon_1 .\nonumber 
\end{align}
To get the conclusion, it is easy to check that $\epsilon_1= 2^{-25n^2}\epsilon$ is enough for the need.

From (\ref{zeta need}), (\ref{delta upper bound 1}) and (\ref{delta upper bound 2}), we get the conclusion.

}
\qed

\begin{lemma}\label{lem weak type 1-1 ineq-first}
{Suppose $(M^n, g)$ has $Rc\geq 0$, for $z\in M^n$, let $f: B_R(z)\rightarrow \mathbb{R}$ be a nonnegative function and $f\in L^1\big(B_R(z)\big)$, then there exists $p\in B_{\frac{R}{2}}(z)$ such that 
\begin{align}
\sup_{t\leq \frac{R}{2}} \fint_{B_t(p)} f\leq 15^n\fint_{B_R(z)} f .\label{small ball integral average}
\end{align}
}
\end{lemma}

\pf
{Let $c= 15^n\fint_{B_R(z)} f$, define
\begin{align}
\mathcal{S}_c= \Big\{x\in B_R(z): \sup_{B_r(x)\subset B_R(z)}\fint_{B_r(x)}f> c\Big\} ,\nonumber 
\end{align}
and $\mathcal{I}(f)(x)= \sup_{x\in B\subset B_R(z)} \fint_B f$, where the supremum is taken over all geodesic balls $B$ containing $x$ and $B\subset B_R(z)$.

Let $T_c= \Big\{x\in B_R(z): \mathcal{I}(f)(x)> c\Big\}$, then $\mathcal{S}_c\subset T_c$. Choose any $S\subset \subset T_c$, for any $x\in S$, there exists geodesic ball $B_x\subset B_R(z)$ such that $x\in B_x$ and 
\begin{align}
V(B_x)< \frac{1}{c} \int_{B_x} f .\label{use it to control volume-1}
\end{align}

By compactness of $S$, we can select a finite collection of such balls $\{B_\alpha\}_{\alpha\in \mathtt{F}}$ that cover $S$. Set $R_0= \sup\big\{r_{\alpha}|\ \alpha\in \mathtt{F}\big\}$, and $\mathrm{F}_j= \big\{\alpha\in \mathtt{F}|\ \frac{R_0}{2^j}< r_\alpha\leq \frac{R_0}{2^{j- 1}}\big\}$, $j=1, 2, \cdots$.

We define $\mathrm{G}_j\subset \mathrm{F}_j$ as follows:
\begin{enumerate}
\item[(1)] Let $\mathrm{G}_1$ be any maximal collection of $\alpha$ in $\mathrm{F}_1$, such that $\{B_{r_\alpha}(x_{\alpha})\}_{\alpha\in \mathrm{G}_1}$ are disjoint.
\item[(2)] Assume $\mathrm{G}_1, \cdots, \mathrm{G}_{k- 1}$ have been selected, choose $\mathrm{G}_k$ to be any maximal subcollection of
\begin{align}
\big\{\alpha\in \mathrm{F}_k|\ B_{r_\alpha}(x_\alpha)\cap B_{r_{\alpha'}}(x_{\alpha'})= \emptyset\ for\ any \ \alpha'\in \bigcup_{j= 1}^{k- 1} \mathrm{G}_j\big\}, \nonumber 
\end{align} 
which also satisfies $B_{r_\alpha}(x_\alpha)\cap B_{r_\beta}(x_\beta)= \emptyset$ for any $\alpha, \beta\in G_k, \alpha\neq \beta$.
\end{enumerate}

Now we define $\mathtt{G}= \bigcup_{k= 1}^{\infty} \mathrm{G}_k\subset \mathtt{F}$. For any $\alpha\in \mathtt{F}$, there is $j\in \mathbb{N}$, such that $\alpha\in \mathrm{F}_j$. If $\alpha\in \mathrm{G}_j$, we get $B_{r_\alpha}(x_\alpha)\subset \bigcup_{\alpha ' \in \mathtt{G}} B_{5r_{\alpha '}}(x_{\alpha'})$.

Otherwise $\alpha\notin \mathrm{G}_j$, by the definition of $\mathrm{G}_j$, there is $\alpha'\in \bigcup_{i= 1}^j \mathrm{G}_i$ such that $B_{r_\alpha}(x_\alpha)\cap B_{r_{\alpha'}}(x_{\alpha'})\neq \emptyset$. Note $r_{\alpha}\leq \frac{R_0}{2^{j- 1}}< 2r_{\alpha'}$, hence $B_{r_{\alpha}}(x_{\alpha})\subset B_{r_{\alpha'}+ 2r_{\alpha}}(x_{\alpha'})\subset \bigcup_{\beta \in \mathtt{G}} B_{5r_{\beta}}(x_{\beta})$.

Then we find $\mathtt{G}\subset \mathtt{F}$, such that for any $\alpha'\neq \beta'\in \mathtt{G}$, $B_{r_{\alpha'}}(x_{\alpha'})\cap B_{r_{\beta'}}(x_{\beta'})= \emptyset$ and 
\begin{align}
\bigcup_{\alpha\in \mathtt{F}} B_{r_{\alpha}}(x_{\alpha}) \subset \bigcup_{\alpha ' \in \mathtt{G}} B_{5r_{\alpha '}}(x_{\alpha'}) .\nonumber
\end{align}

From Bishop-Gromov Comparison Theorem and (\ref{use it to control volume-1}), we get 
\begin{align}
\mu(S)&\leq \mu\big(\bigcup_{\alpha\in \mathtt{F}} B_\alpha\big)\leq \mu\big(\bigcup_{\alpha ' \in \mathtt{G}} B_{5r_{\alpha '}}(x_{\alpha'})\big) \leq 5^n\sum_{\alpha ' \in \mathtt{G}} \mu(B_{r_{\alpha '}}(x_{\alpha '})) \nonumber \\
&< 5^n\cdot \sum_{\alpha ' \in \mathtt{G}} \Big(\frac{1}{c}\int_{B_{r_{\alpha '}}(x_{\alpha '})} f\Big) \leq \frac{5^n}{c}\int_{B_R(z)} f . \nonumber 
\end{align}

If one takes the supremum over all such $S\subset \subset T_c$, we have 
\begin{align}
\mu(\mathcal{S}_c)\leq \mu(T_c)\leq \frac{5^n}{c}\int_{B_R(z)} f\leq 3^{-n}V\big(B_R(z)\big) .\nonumber 
\end{align}

If (\ref{small ball integral average}) does not hold for any point in $B_{\frac{R}{2}}(z)$, then we have $B_{\frac{R}{2}}(z)\subset \mathcal{S}_{c}$, from Bishop-Gromov Volume Comparison Theorem, $\frac{\mu(\mathcal{S}_{c})}{\mu\big(B_R(z)\big)}\geq \frac{\mu\big(B_{\frac{R}{2}}(z)\big)}{\mu\big(B_R(z)\big)}\geq 2^{-n}$, it is the contradiction, the conclusion follows.
}
\qed

The following lemma is well known so we omit its proof here.
\begin{lemma}\label{lem one side G-H map imply G-H dist}
{Let $(\mathbf{X}, d_\mathbf{X}, x_0)$ and $(\mathbf{Y}, d_\mathbf{Y}, y_0)$ be two pointed metric spaces, if there is a pointed $\epsilon$-Gromov-Hausdorff approximation $f: \big(\mathbf{X}, x_0\big)\rightarrow (\mathbf{Y}, y_0)$, then there exists a pointed $(3\epsilon)$-Gromov-Hausdorff approximation $h: \big(\mathbf{Y}, y_0\big)\rightarrow (\mathbf{X}, x_0)$. 
}
\end{lemma}\qed

\begin{theorem}\label{thm AG-triple imply one more splitting}
{For $B_{10r}(q)\subset (M^n, g)$ with $Rc(g)\geq 0$, any $0< \epsilon< 1, \delta= n^{-3400n^3}\epsilon^{110n}$ and integer $0\leq k\leq n$, assume there is an $(\delta r)$-Gromov-Hausdorff approximation ,
\begin{align}
f: B_{10r}(q)\rightarrow B_{10r}(0, \hat{q})\subset \mathbb{R}^k\times \mathbf{X}_k ,\nonumber 
\end{align}
and $\mathrm{diam}\big(B_r(\hat{q})\big)= r_0\geq \frac{1}{4}r$. Then there are harmonic functions $\big\{\mathbf{b}_i\big\}_{i= 1}^{k+ 1}$ defined on some geodesic ball $B_{s}(p)\subset B_{10r}(q)$, where $s= n^{-320n^3}\epsilon^{10n}r$, such that 
\begin{align}
\sup_{B_{s}(p)\atop i= 1, \cdots, k+ 1} |\nabla \mathbf{b}_i|\leq 1+ \epsilon \quad \quad and \quad \quad \sup_{t\leq s}\fint_{B_{t}(p)}\sum_{i, j= 1}^{k+ 1} \big|\langle \nabla \mathbf{b}_i, \nabla \mathbf{b}_j\rangle- \delta_{ij}\big|\leq \epsilon .\nonumber 
\end{align}
}
\end{theorem}

\pf
{From Lemma \ref{lem one side G-H map imply G-H dist}, there is an $(3\delta r)$-Gromov-Hausdorff approximation,
\begin{align}
\Phi: B_{10r}(0, \hat{q})\rightarrow B_{10r}(q)  .\nonumber 
\end{align}

Let $r_1= n^{-300n^3}\epsilon_1^{10n}r$, assume $\delta\leq \frac{1}{3}n^{-3200n^3}\epsilon_1^{110n}$, where $\epsilon_1> 0$ is to be determined later. Apply Proposition \ref{thm AG-triple imply one more splitting-pre}, we obtain harmonic functions $\big\{\mathbf{b}_i\big\}_{i= 1}^{k+ 1}$ defined on $B_{r_1}(q_1)\subset B_{10r}(q)$, such that 
\begin{align}
\sup_{B_{r_1}(q_1)\atop i= 1, \cdots, k+ 1} |\nabla \mathbf{b}_i|\leq 1+ \epsilon_1 \quad \quad and \quad \quad \fint_{B_{r_1}(q_1)}\sum_{i, j= 1}^{k+ 1} \big|\langle \nabla \mathbf{b}_i, \nabla \mathbf{b}_j\rangle- \delta_{ij}\big|\leq \epsilon_1 . \nonumber 
\end{align}

Apply Lemma \ref{lem weak type 1-1 ineq-first}, we get $B_{s}(p)\subset B_{r_1}(q_1)$, where $s= \frac{r_1}{2}= \frac{1}{2} n^{-300n^3}\epsilon_1^{10n}r$, such that 
\begin{align}
\sup_{B_{s}(p)\atop i= 1, \cdots, k+ 1} |\nabla \mathbf{b}_i|\leq 1+ \epsilon_1 \quad \quad and \quad \quad \sup_{t\leq s}\fint_{B_{t}(p)}\sum_{i, j= 1}^{k+ 1} \big|\langle \nabla \mathbf{b}_i, \nabla \mathbf{b}_j\rangle- \delta_{ij}\big|\leq 15^n\epsilon_1 .\nonumber 
\end{align}

Choose $\epsilon_1= 15^{-n}\epsilon$, let $\delta= n^{-3400n^3}\epsilon^{110n}$ and $s= n^{-320n^3}\epsilon^{10n}r$, the conclusion follows.
}
\qed

\section*{Part II. A.O.L. harmonic functions produce G-H approximation}

In Part II of this paper, we will prove the following quantitative version of almost splitting theorem, which was implied in the argument of Cheeger and Colding in a series of papers, \cite{Colding-shape}, \cite{Colding-large}, \cite{Colding-volume}, \cite{CC-Ann} and \cite{Cheeger-GAFA}. 

For our application, we need the quantitative estimate, which relate the Gromov-Hausdorff distance to the average integral bound of almost orthonormal linear harmonic functions. Although we believe that many results in Part II are well-known to some experts in this field, but we can not find the reference providing those quantitative estimates exactly. So we elaborate the concise argument of Cheeger-Colding to present the proof in all the details for self-contained reason.
  
\begin{theorem}\label{thm existence of harmonic map components imply GH-dist is small}
{For $\epsilon> 0$ and $1\leq k\leq n$, there is $\delta= n^{-700n^4}\epsilon^{18n^4}$ such that for complete Riemannian manifold $(M^n, g)$ with $Rc\geq 0$, if there exist harmonic functions $\big\{\mathbf{b}_i\big\}_{i= 1}^k$, defined on $B_r(p)$, satisfying $\mathbf{b}_i(p)= 0$ and 
\begin{align}
\sup_{B_r(p)}|\nabla \mathbf{b}_i|\leq 2  \ , \quad \quad  \quad \quad \fint_{B_r(p)} \sum_{i, j= 1}^k \big|\langle \nabla \mathbf{b}_i, \nabla \mathbf{b}_j\rangle- \delta_{ij}\big|\leq \delta ,\label{assump on harmonic linear function}
\end{align}
then we can find a metric space $\mathbf{X}_k$ and an $(\epsilon s)$-Gromov-Hausdorff approximation $f_k= (\mathbf{b}_1, \cdots, \mathbf{b}_k, \mathcal{P}_k): B_{s}(p)\rightarrow B_{s}(0, \hat{p})\subset \mathbb{R}^k\times \mathbf{X}_k$, where $s= \frac{1}{1280} r$.
}
\end{theorem}

\begin{remark}
{Colding and Minicozzi \cite{CM} gave a characterization of Gromov-Hausdorff distance through the integral estimate of Hessian of harmonic functions among other things.
}
\end{remark}

\section{Existence of almost linear function and Hessian estimate}\label{SEC existence of linear function}

When there is a harmonic function $\mathbf{b}$ defined locally on manifold $M^n$, with bounded gradient and the average integral of $\big||\nabla \mathbf{b}|- 1\big|$ is small enough, we will show the existence of an almost linear function, which is a generalization of linear function in $\mathbb{R}^n$. The proof of Proposition \ref{prop existence of almost linear function} has close relationship with the argument in \cite{Cheeger-GAFA}.

\begin{definition}
{For $\mathbf{X}\subset M^n, t\in \mathbb{R}$, the function $\rho(x)= d(x, \mathbf{X})+ t: M^n\rightarrow \mathbb{R}$ is called \textbf{almost linear function}, which is the generalization of function $f(x)= d(x, \mathbb{R}^{n- 1})+ t$ defined on $\mathbb{R}^n$. 
}
\end{definition}

\begin{definition}\label{def Lip constant function}
{For Lipschitz function $f$ defined on metric space $\mathbf{M}$, the \textbf{pointwise Lipschitz constant function} $\mathscr{L}(f)$ is defined by 
\begin{align}
\mathscr{L}(f)(z)= \varlimsup_{d(z, z')\rightarrow 0} \frac{\big|f(z)- f(z')\big|}{d(z, z')}\ , \quad \quad \quad \quad z, z'\in \mathbf{M} ,\nonumber 
\end{align}
and the \textbf{Lipschtiz constant} $\mathbf{L}(f)$ is defined by $\mathbf{L}(f)= \sup_{z\in \mathbf{M}} \big\{\mathscr{L}(f)(z)\big\}$.
}
\end{definition}

\begin{remark}\label{rem Lip implies almost differentiable}
{From the classical Rademacher theorem, the Lipschitz function is almost differentiable on manifolds, for the general argument on metric measure spaces see \cite{Cheeger-GAFA}. Hence when the pointwise Lipschitz constant function $\mathscr{L}(f)(z)$ appears as the integrand function in an integral, we can replace it by $|\nabla f|(z)$, and we will use this fact freely in the following argument.
}
\end{remark}

\begin{prop}\label{prop existence of almost linear function}
{For any $\epsilon> 0$, there is $\delta= 2^{-100n^2}\epsilon^{2n^2}$ such that for any complete Riemannian manifold $(M^n, g)$ with $Rc\geq 0$, if there exists one harmonic function $\mathbf{b}$ defined on $B_r(p)$ satisfying $\sup_{B_r(p)}|\nabla \mathbf{b}|\leq 2$ and $\fint_{B_r(p)} \big||\nabla \mathbf{b}|- 1\big|\leq \delta$. Then one can find $t_0, t_1\in \mathbb{R}$ and two functions $\rho, \tilde{\rho}$ defined on $B_{\frac{r}{10}}(p)$, such that
\begin{align}
&\rho(x)= d\big(x, \rho^{-1}(t_0)\big)+ t_0 \ ,  \quad \quad \tilde{\rho}(x)= t_1- d\big(x, \tilde{\rho}^{-1}(t_1)\big)\ , \quad  \forall \ x\in B_{\frac{r}{160}}(p) \nonumber \\
&\frac{1}{320}r\leq d\big(x, \rho^{-1}(t_0)\big)\leq \frac{3}{320}r\ , \quad \quad \quad \quad \forall x\in B_{\frac{r}{320}}(p) \label{need 3.6.0} \\
&\sup_{B_{\frac{1}{20}r}(p)}\big|\mathbf{b}- \rho\big|\leq \epsilon\cdot r \ , \quad \quad \quad \quad \fint_{B_{\frac{1}{20}r}(p)} \big|\nabla (\mathbf{b}- \rho)\big|\leq \epsilon \nonumber \\
&\rho(x)\leq \tilde{\rho}(x)+ \frac{\epsilon}{2} r \ , \quad \quad \quad \quad \forall x\in B_{\frac{1}{10}r}(p) .\nonumber  
\end{align}
}
\end{prop}

\pf
{\textbf{Step (1)}. For $1> \theta> 0$ (to be determined later), we define 
\begin{align}
\mathbf{A}= \Big\{x\in B_r(p)\big| \ \big||\nabla \mathbf{b}(x)|- 1\big|< \frac{\delta}{\theta}\Big\} ,\nonumber 
\end{align}
then from assumption, we get 
\begin{align}
V\big(B_r(p)- \mathbf{A}\big)\leq  \theta\cdot V\big(B_r(p)\big) . \label{small bad set}
\end{align}

We define $\mathbf{b}_*$ as the following:
\begin{equation}\nonumber 
\mathbf{b}_*(x)= \left\{
\begin{array}{rl}
&\mathbf{b}(x) \quad \quad \quad \quad \quad \quad \quad \quad \quad x\in \mathbf{A} \\
&\sup\limits_{z'\in \mathbf{A}}\big[\mathbf{b}(z')- (1+ \frac{\delta}{\theta})\cdot d(x, z')\big] \quad \quad \quad \quad \quad x\in B_r(p)- \mathbf{A} .\nonumber  
\end{array} \right.
\end{equation}
It is easy to get $\sup_{x\in B_r(p)}\big|\mathscr{L}(\mathbf{b}_*)(x)\big|\leq 1+ \frac{\delta}{\theta}$.

Put $h(x)= \frac{\theta}{\delta+ \theta}\mathbf{b}_*(x)+ \frac{\delta}{\delta+ \theta}\mathbf{b}_*(p)$, then 
\begin{align}
\sup_{x\in B_r(p)}\big|\mathscr{L}(h)(x)\big|= \frac{\theta}{\delta+ \theta}\sup_{x\in B_r(p)}\big|\mathscr{L}(\mathbf{b}_*)(x)\big|\leq 1 .   \nonumber 
\end{align}

Note $h(p)= \mathbf{b}_*(p)$ and $\sup\limits_{B_r(p)}\big|\mathscr{L} (h- \mathbf{b}_*)\big|\leq \frac{\delta}{\theta}$, we get $\sup_{B_r(p)} |h- \mathbf{b}_*|\leq \frac{\delta}{\theta} r$. 

For any $x\in \overline{B_{\frac{r}{5}}(p)}$, from (\ref{small bad set}), there exists $x'\in B_{\frac{7}{5}\theta^{\frac{1}{n}} r}(x)$ such that $x'\in \mathbf{A}$. Note $\mathbf{b}(x')- \mathbf{b}_*(x')= 0$ and 
\begin{align}
\sup_{B_r(p)}\big|\mathscr{L}(\mathbf{b}- \mathbf{b}_*)\big| \leq \sup_{B_r(p)} |\nabla \mathbf{b}|+ \sup_{B_r(p)} |\mathscr{L}(\mathbf{b}_*)|\leq 3+ \frac{\delta}{\theta}. \nonumber
\end{align}

Then for any $x\in \overline{B_{\frac{r}{5}}(p)}$, we have 
\begin{align}
\big|\mathbf{b}(x)- \mathbf{b}_*(x)\big|\leq \big|\mathbf{b}(x')- \mathbf{b}_*(x')\big|+ (3+ \frac{\delta}{\theta})\cdot \frac{7}{5}\theta^{\frac{1}{n}} r = (3+ \frac{\delta}{\theta})\cdot \frac{7}{5}\theta^{\frac{1}{n}} r .\nonumber 
\end{align}

From above we have 
\begin{align}
\sup_{\overline{B_{\frac{r}{5}}(p)}}|h- \mathbf{b}|&\leq 5\frac{\delta}{\theta} r+ (3+ \frac{\delta}{\theta})\cdot \frac{7}{5}\theta^{\frac{1}{n}} r\leq \big[6\theta^{\frac{1}{n}}+ 7\frac{\delta}{\theta}\big]r .\label{h- f}
\end{align}

From (\ref{small bad set}), note $h= \frac{\theta}{\delta+ \theta}\mathbf{b}+ \frac{\delta}{\delta+ \theta}\mathbf{b}_*(p)$ on $\mathbf{A}$, we have 
\begin{align}
&\fint_{B_{\frac{r}{5}}(p)} \big(1- \langle \nabla \mathbf{b}, \nabla h\rangle\big)\leq \frac{1}{V\big(B_{\frac{r}{5}}(p)\big)} \Big[\int_{B_{\frac{r}{5}}(p)- \mathbf{A}} 3 + \int_{\mathbf{A}\cap B_{\frac{r}{5}}(p)} \big(1- \langle \nabla \mathbf{b}, \nabla h\rangle\big) \Big] \nonumber \\
&\quad\quad\quad\quad \leq 3\frac{V\big(B_{\frac{r}{5}}(p)- \mathbf{A}\big)}{V\big(B_{\frac{r}{5}}(p)\big)}+  \frac{1}{V\big(B_{\frac{r}{5}}(p)\big)} \int_{\mathbf{A}\cap B_{\frac{r}{5}}(p)} \Big\{1- \langle \nabla \mathbf{b}, \frac{\theta}{\delta+ \theta} \nabla \mathbf{b}\rangle \Big\} \nonumber \\
&\quad\quad\quad\quad  \leq 3\theta \frac{V\big(B_{r}(p)\big)}{V\big(B_{\frac{r}{5}}(p)\big)}+ \frac{1}{V\big(B_{\frac{r}{5}}(p)\big)} \int_{\mathbf{A}\cap B_{\frac{r}{5}}(p)} \Big\{\frac{\delta}{\delta+ \theta}+ \frac{\theta}{\delta+ \theta}\big(1- |\nabla \mathbf{b}|^2\big)\Big\} \nonumber \\
&\quad\quad \quad \quad \leq 3\cdot 5^n\theta+ \frac{4\delta}{\theta+ \delta}\leq 4\big(5^n\theta+ \frac{\delta}{\theta}\big)  .\label{gradient of f- h}
\end{align}

Now from (\ref{h- f}) and (\ref{gradient of f- h}), we can choose $\theta= 2^{-12n^2}\epsilon_1^n$, where $\epsilon_1> 0$ is to be determined later, and 
\begin{align}
\delta\leq 2^{-20n^2}\epsilon_1^{n+ 1} .\label{delta need 1}
\end{align}
Then we have 
\begin{align}
\sup_{\overline{B_{\frac{r}{5}}(p)}}|h- \mathbf{b}|\leq \epsilon_1 r \quad \quad and \quad \quad  \fint_{B_{\frac{r}{5}}(p)} \big(1- \langle \nabla \mathbf{b}, \nabla h\rangle\big)\leq \epsilon_1 .\label{step 1 ineq}
\end{align}

\textbf{Step (2)}. Now for $x\in \overline{B_{\frac{1}{5}r}(p)}$, we consider 
\begin{align}
\tilde{\rho}(x)&= \sup_{z'\in \partial B_{\frac{1}{5}r}(p)} \big[h(z')- d(x, z')\big]\ ,  \nonumber \\
\rho(x)&= \inf_{z'\in \partial B_{\frac{1}{5}r}(p)} \big[h(z')+ d(x, z')\big]. \nonumber
\end{align}

For $x\in B_{\frac{r}{10}}(p)$, assume 
\begin{align}
\rho(x)= h(x^*)+ d(x^*, x), \nonumber 
\end{align}
where $x^*\in \partial B_{\frac{r}{5}}(p)$. Let $\gamma_{x, x^*}(s)$ be the minimizing geodesic from $x$ to $x^*$, parametrized by arc-length. From $\mathbf{L}(\rho)\leq 1$, we get
\begin{align}
\rho\big(\gamma_{x, x^*}(s)\big)\geq \rho(x)- s= h(x^*)+ d(x, x^*)- s= h(x^*)+ d\big(\gamma_{x, x^*}(s), x^*\big) \geq \rho\big(\gamma_{x, x^*}(s)\big) ,\nonumber 
\end{align}
which implies
\begin{align}
\rho\big(\gamma_{x, x^*}(s)\big)= \rho(x)- s \ , \quad \quad \quad \quad \forall s\in [0, d(x, x^*)] .\nonumber 
\end{align}

But $d(x, x^*)\geq \frac{r}{10}$, we have 
\begin{align}
\rho\big(\gamma_{x, x^*}(s)\big)= \rho(x)- s \ , \quad \quad \quad \quad \forall s\in \big[0, \frac{r}{10}\big] .\label{barrier along geodesic-1}
\end{align}

Let $t_0= \rho(p)- \frac{r}{160}$, note $0\leq \rho(x)- t_0\leq \frac{r}{10}$ for $x\in B_{\frac{r}{160}}(p)$, then from (\ref{barrier along geodesic-1}),
\begin{align}
t_0= \rho(x)- \big[\rho(x)- t_0\big]= \rho\Big(\gamma_{x, x^*}\big(\rho(x)- t_0\big)\Big) .\nonumber 
\end{align}

Now we have 
\begin{align}
\rho(x)- t_0\leq d\big(x, \rho^{-1}(t_0)\big)\leq d\Big(x, \gamma_{x, x^*}\big(\rho(x)- t_0\big)\Big)= \rho(x)- t_0 .\nonumber 
\end{align}

By the above, for $x\in B_{\frac{r}{160}}(p)$, we get 
\begin{align}
\rho(x)= t_0+ d\big(x, \rho^{-1}(t_0)\big) .\label{barrier linear-1} 
\end{align}

Similarly, let $t_1= \tilde{\rho}(p)+ \frac{r}{160}$, for $x\in B_{\frac{r}{160}}(p)$, we get 
\begin{align}
\tilde{\rho}(x)= t_1- d\big(x, \tilde{\rho}^{-1}(t_1)\big) . \nonumber 
\end{align}

For any $x\in B_{\frac{r}{320}}(p)$, using $\mathbf{L}(\rho)\leq 1$, we have 
\begin{align}
d\big(x, \rho^{-1}(t_0)\big)= \rho(x)- t_0= \rho(x)- \rho(p)+ \frac{r}{160}\in \big[\frac{r}{320}, \frac{3r}{320}\big] .\nonumber 
\end{align}

\textbf{Step (3)}. Let $\epsilon_2= 2^{11}\epsilon_1^{\frac{1}{n}}$, we will prove that 
\begin{align}
\rho(w)\leq \tilde{\rho}(w)+ \epsilon_2 r \ , \quad \quad \quad \quad \forall w\in B_{\frac{1}{10}r}(p) \label{upper and lower barrier ineq}
\end{align}

By contradiction, if there exists $w\in B_{\frac{1}{10}r}(p)$ such that 
\begin{align}
\tilde{\rho}(w)+ 2\theta_1 r\leq \rho(w) ,\nonumber 
\end{align}
where $\theta_1= 2^{10}\epsilon_1^{\frac{1}{n}}\leq \frac{1}{20}$. Let $K= \frac{1}{2}\big[\tilde{\rho}(w)+ \rho(w) \big]$, since $\big| \mathbf{L}(\rho)\big|\leq 1, \big|\mathbf{L}(\tilde{\rho})\big|\leq 1$ and $\tilde{\rho}\leq \rho$, also note $B_{\theta r}(w)\subset B_{\frac{1}{5}r}(p)$, we have 
\begin{align}
\tilde{\rho}(x)\leq K\leq \rho(x)\ , \quad \quad \quad \quad x\in B_{\theta_1 r}(w) .\nonumber 
\end{align}

Now we define 
\begin{equation}\nonumber 
\check{h}(x)=\left\{
\begin{array}{rl}
&\rho(x)\ , \quad\quad \quad \quad if\ \rho(x)\leq K  \\
&K \ , \quad\quad \quad \quad if\ \tilde{\rho}(x)\leq K\leq \rho(x)  \\
&\tilde{\rho}(x)\ , \quad\quad \quad \quad if\  K\leq \tilde{\rho}(x) ,
\end{array} \right.
\end{equation}
then $\big|\mathbf{L}(\check{h})\big|\leq 1$, also note $\check{h}= \tilde{\rho}= \rho= h$ on $\partial B_{\frac{1}{5}r}(p)$, from (\ref{step 1 ineq}), we get
\begin{align}
\sup_{\partial B_{\frac{1}{5}r}(p)} |\check{h}- \mathbf{b}|=\sup_{\partial B_{\frac{1}{5}r}(p)} |h- \mathbf{b}|\leq \epsilon_1 r .\nonumber 
\end{align}

We define the Lipschitz function $\phi(x): \overline{B_{\frac{1}{5}r}(p)}\rightarrow \mathbb{R}$ satisfying 
\begin{align}
0\leq \phi(x)\leq 1\quad \quad and \quad  \quad \mathbf{L}(\phi)\leq (\epsilon_1 r)^{-1} , \nonumber 
\end{align}
also 
\begin{equation}\nonumber 
\phi(x)= \left\{
\begin{array}{rl}
&0 \quad \quad \quad \quad x\in \partial B_r(p) \\
&1 \quad \quad \quad \quad x\in  B_{(\frac{1}{5}- \epsilon_1)r}(p)  .
\end{array} \right.
\end{equation}

Now define $\tilde{\mathbf{b}}= \check{h}+ (1- \phi)(\mathbf{b}- \check{h})$, we have 
\begin{align}
\Big|\mathbf{L} \big[(1- \phi)\cdot (\mathbf{b}- \check{h})\big]\Big| &\leq  \sup_{B_{\frac{r}{5}}(p)} (1- \phi) \cdot \big|\mathscr{L} (\mathbf{b}- \check{h})\big| + \sup_{B_{\frac{r}{5}}(p)\backslash B_{(\frac{1}{5}- \epsilon_1)r}(p)} |\mathbf{b}- \check{h}|\cdot \big|\mathscr{L} (1- \phi)\big| \nonumber \\
&\leq  \big[|\nabla\mathbf{b}|+ |\mathbf{L}(\check{h})|\big]+ \big[\sup_{\partial B_{\frac{r}{5}}(p)} |\mathbf{b}- \check{h}|+ (\epsilon_1 r)\cdot \sup_{B_{\frac{r}{5}}(p)}|\mathscr{L}(\mathbf{b}- \check{h})|\big]\cdot (\epsilon_1 r)^{-1} \nonumber \\
&\leq 3+ \big[\epsilon_1r+ 3\epsilon_1r\big](\epsilon_1r)^{-1}= 7 . \nonumber 
\end{align}

Let $\chi_E$ denote the characteristic function of a set $E$, then from above and the Bishop-Gromov Comparison Theorem we get 
\begin{align}
\fint_{B_{\frac{r}{5}}(p)} |\nabla\tilde{\mathbf{b}}|^2&\leq \fint_{B_{\frac{r}{5}}(p)} \Big(|\nabla\check{h}|+ \big|\nabla \big[(1- \phi)\cdot (\mathbf{b}- \check{h})\big]\big|\Big)^2\leq \fint_{B_{\frac{r}{5}}(p)} \Big(|\nabla\check{h}|+ 7\chi_{B_{\frac{r}{5}}\backslash B_{(\frac{1}{5}- \epsilon_1)r}}\Big)^2 \nonumber \\
&\leq \fint_{B_{\frac{r}{5}}(p)} |\nabla\check{h}|^2+ 70\Big(1- \frac{V\big(B_{(\frac{1}{5}- \epsilon_1)r}(p)\big)}{V\big(B_{\frac{r}{5}}(p)\big)}\Big) \leq \fint_{B_{\frac{r}{5}}(p)} |\nabla\check{h}|^2+ 70\big(1- (1- 5\epsilon_1)^n\big) \nonumber \\
&\leq 350n\cdot \epsilon_1+ \fint_{B_{\frac{r}{5}}(p)} |\nabla\check{h}|^2 . \nonumber 
\end{align}

From the fact that $\tilde{\mathbf{b}}|_{\partial B_{\frac{r}{5}}(p)}= \mathbf{b}|_{\partial B_{\frac{r}{5}}(p)}$ and $\mathbf{b}$ is harmonic, note the harmonic function has the smallest energy, we have 
\begin{align}
\fint_{B_{\frac{r}{5}}(p)} |\nabla \mathbf{b}|^2\leq \fint_{B_{\frac{r}{5}}(p)} |\nabla \tilde{\mathbf{b}}|^2\leq 350n\epsilon_1+ \fint_{B_{\frac{r}{5}}(p)} |\nabla\check{h}|^2  ,\nonumber 
\end{align}
which implies
\begin{align}
\fint_{B_{\frac{r}{5}}(p)} \Big[1- |\nabla\check{h}|^2 \Big]&\leq 350n\epsilon_1+ \fint_{B_{\frac{r}{5}}(p)} \Big[1- |\nabla \mathbf{b}|^2\Big] \leq 350n\epsilon_1+ 3\cdot 5^n\delta \leq 400n\epsilon_1 . \nonumber 
\end{align}
In the last inequality above, we used (\ref{delta need 1}).

On the other hand, note $\mathscr{L}(\check{h})= 0$ on $B_{\theta_1 r}(w)$, we have 
\begin{align}
400n\epsilon_1 &\geq \fint_{B_{\frac{r}{5}}(p)} \Big[1- |\nabla\check{h}|^2 \Big]\geq \frac{1}{V\big(B_{\frac{r}{5}}(p)\big)} \int_{B_{\theta_1 r}(w)} \Big[1- |\nabla\check{h}|^2 \Big]= \frac{V\big(B_{\theta_1 r}(w)\big)}{V\big(B_{\frac{r}{5}}(p)\big)}\geq \big(\frac{10\theta_1}{3}\big)^n ,\nonumber 
\end{align}
which implies $\theta_1< 2^{10}\epsilon_1^{\frac{1}{n}}$, it is the contradiction with  the choice of $\theta_1$.

\textbf{Step (4)}. Note $\tilde{\rho}(x)\leq h(x)\leq \rho(x)$ and (\ref{upper and lower barrier ineq}), then we get 
\begin{align}
\big|h(x)- \rho(x)\big|\leq \rho(x)- \tilde{\rho}(x)\leq \epsilon_2 r\ , \quad \quad \quad \forall x\in B_{\frac{r}{20}}(p). \label{the diff between h and rho}
\end{align}

From (\ref{step 1 ineq}) and (\ref{the diff between h and rho}) we have 
\begin{align}
\sup_{B_{\frac{1}{20}r}(p)} |\mathbf{b}- \rho|\leq (2^{11}+ 1)\epsilon_1^{\frac{1}{n}}r\leq 2^{12}\epsilon_1^{\frac{1}{n}}r  . \nonumber 
\end{align}

From (\ref{delta need 1}) and (\ref{step 1 ineq}),
\begin{align}
&\quad \fint_{B_{\frac{1}{20}r}(p)} \big|\nabla (\mathbf{b}- \rho)\big|^2 \nonumber \\
&\leq  4^n\fint_{B_{\frac{1}{5}r}(p)} \big|\nabla (\mathbf{b}- \rho)\big|^2= 4^n\fint_{B_{\frac{1}{5}r}(p)} |\nabla \mathbf{b}|^2+ |\nabla \rho|^2- 2\langle \nabla \mathbf{b}, \nabla \rho\rangle \nonumber \\
&\leq 4^n\fint_{B_{\frac{1}{5}r}(p)}  \big( |\nabla \mathbf{b}|^2- 1\big)+ \big(|\nabla \rho|^2- 1\big)+ 2\big(1- \langle \nabla \mathbf{b}, \nabla h\rangle\big)+ 2 \langle \nabla \mathbf{b}, \nabla (h- \rho)\rangle \nonumber \\
&\leq 4^n\big[3\cdot 5^n\delta+ 0+ 2\epsilon_1+ 0\big]=4^n\cdot (3\epsilon_1) , \nonumber 
\end{align}
then 
\begin{align}
\fint_{B_{\frac{1}{20}r}(p)} \big|\nabla (\mathbf{b}- \rho)\big|\leq \Big(\fint_{B_{\frac{1}{20}r}(p)} \big|\nabla (\mathbf{b}- \rho)\big|^2\Big)^{\frac{1}{2}}\leq 2^{n+ 1}\sqrt{\epsilon_1} .\nonumber  
\end{align}

If we choose $\epsilon_1\leq 2^{-23n}\epsilon^n$ and also $2^{10}\epsilon_1^{\frac{1}{n}}\leq \frac{1}{20}$, then we obtain the conclusion. From (\ref{delta need 1}), we only need to choose $\delta\leq 2^{-100n^2}\epsilon^{2n^2}$.
}
\qed

\section{Segment inequality and measure of `good' points}

In the proofs of this section, when the context is clear, for simplicity, we use $B_r$ instead of $B_r(p)$, similar for $B_{4r}$ etc.

\begin{lemma}[Segment Inequality]\label{lem segment ineq}
{Assume $(M^n, g)$ is a complete Riemannian manifold with $Rc\geq 0$, then for any nonnegative function $f$ defined on $B_{2r}(p)\subset M^n$, 
\begin{align}
\int_{B_r(p)\times B_r(p)} \Big(\int_0^{d(y_1, y_2)} f\big(\gamma_{y_1, y_2}(s)\big) ds\Big) dy_1 dy_2\leq 2^{n+ 1}r\cdot V\big(B_r(p)\big)\cdot \int_{B_{2r}(p)} f . \nonumber 
\end{align}
}
\end{lemma}

\pf
{In the proof, we assume $y_1, y_2\in B_r$, and set  
\begin{align}
& E(y_1, y_2)= \int_0^{d(y_1, y_2)} f\big(\gamma_{y_1, y_2}(s)\big) ds \nonumber \\
& E_1(y_1, y_2)= \int_{\frac{1}{2}d(y_1, y_2)}^{d(y_1, y_2)} f\big(\gamma_{y_1, y_2}(s)\big) ds\ , \quad \quad \quad \quad E_2(y_1, y_2)= \int_0^{\frac{1}{2}d(y_1, y_2)} f\big(\gamma_{y_1, y_2}(s)\big) ds . \nonumber 
\end{align}
Then $E= E_1+ E_2$. Along any geodesic $\gamma$ starting from $y_1$, write the volume element of $M^n$ in geodesic polar coordinate as $ds\wedge \mathcal{A}(s)$. Then $\mathcal{A}(s)= J(\theta, s)d\theta$, from Bishop-Gromov Comparison Theorem, 
\begin{align}
\mathcal{A}(s)\leq \big(\frac{s}{u}\big)^{n- 1}\mathcal{A}(u)\leq 2^{n- 1}\cdot \mathcal{A}(u) \ , \quad \quad \quad \quad \forall u\in [\frac{s}{2}, s] .\nonumber 
\end{align} 

For $y\in B_r, v\in S_{y}M^n$, we define 
\begin{align}
I(y, v)= \Big\{t\geq 0| \ \gamma(t)\in B_r, \gamma '(0)= v, \gamma|_{[0, t]} \ is \ minimal \Big\} .\nonumber 
\end{align}
Then we have 
\begin{align}
\sup_{y\in B_r, v\in S_{y}M^n} \big|I(y, v)\big|\leq 2r ,\label{inteval upper bound}
\end{align}
where $\big|I(y, v)\big|$ denotes the measure of $I(y, v)$. 

Assume $\gamma_{y_1}^{v_1}$ is the geodesic starting from $y_1$ with $(\gamma_{y_1}^{v_1})'(0)= v_1$, then for any $y_1\in B_r$, $s\in I(y_1, v_1)$,
\begin{align}
E_1\big(y_1, \gamma_{y_1}^{v_1}(s)\big)\mathcal{A}(s)&= \mathcal{A}(s)\int_{\frac{1}{2}s}^s f\big(\gamma_{y_1}^{v_1}(u)\big) du\leq 2^{n- 1}\cdot \int_{\frac{1}{2}s}^s f\big(\gamma_{y_1}^{v_1}(u)\big) \mathcal{A}(u) du \nonumber \\
&\leq 2^{n- 1}\cdot \int_0^{\mathcal{T}(y_1, v_1)} f\big(\gamma_{y_1}^{v_1}(t)\big) \mathcal{A}(t) dt , \nonumber 
\end{align}
where $\mathcal{T}(y_1, v_1)= \max\limits_{t\in I(y_1, v_1)} t$.

Thus from (\ref{inteval upper bound}), for any $y_1\in B_r, v_1\in S_{y_1}M^n$, 
\begin{align}
\int_{I(y_1, v_1)} E_1\big(y_1, \gamma_{y_1}^{v_1}(s)\big) \mathcal{A}(s) ds\leq 2^nr\cdot \int_0^{\mathcal{T}(y_1, v_1)} f\big(\gamma_{y_1}^{v_1}(t)\big)\mathcal{A}(t) dt .\label{seg ineq 1}
\end{align}

Integrating (\ref{seg ineq 1}) with respect to $v_1$ over the unit tangent space $S_{y_1}M^n$, and note $\Big(\bigcup\limits_{y_1\in B_{r}(p) \atop y_2\in B_{r}(p) } \gamma_{y_1, y_2}\Big)\subset B_{2r}$, where $\gamma_{y_1, y_2}$ is the minimal geodesic connecting $y_1$ with $y_2$ in $M^n$, we have 
\begin{align}
\int_{B_r} E_1(y_1, y_2) dy_2\leq 2^nr\cdot \int_{B_{2r}} f \ , \quad \quad \quad \quad \forall y_1\in B_r .\label{seg ineq 2}
\end{align}
And we integrate (\ref{seg ineq 2}) with respect to $y_1$ over $B_r$, 
\begin{align}
\int_{B_r\times B_r}\int_{\frac{1}{2}d(y_1, y_2)}^{d(y_1, y_2)} f\big(\gamma_{y_1, y_2}(s)\big) ds=\int_{B_r\times B_r} E_1(y_1, y_2) \leq 2^nr V(B_r)\cdot \int_{B_{2r}} f .\label{seg ineq 3}
\end{align}
Similarly, we get 
\begin{align}
\int_{B_r\times B_r}\int_{0}^{\frac{1}{2}d(y_1, y_2)} f\big(\gamma_{y_1, y_2}(s)\big) ds=\int_{B_r\times B_r}  E_2(y_1, y_2) \leq 2^nr V(B_r)\cdot \int_{B_{2r}} f .\label{seg ineq 4}
\end{align}

Take the sum of (\ref{seg ineq 3}) and (\ref{seg ineq 4}), the conclusion follows.
}
\qed

For $x\in B_{2r}(p)\subset M^n$ and a closed subset $\mathbf{X}\subseteq M^n$, we define 
\begin{align}
\rho(x)= d\big(x, \mathbf{X}\big)+ t_0\ , \quad \quad \quad \quad \hat{\rho}(x)= d\big(x, \mathbf{X}\big) ,\nonumber 
\end{align}
where $t_0\in \mathbb{R}$ is some constant, we define $\mathfrak{P}(x)\in \mathbf{X}$ by $d\big(x, \mathbf{X}\big)= d\big(x, \mathfrak{P}(x)\big)$ (if there are two points $y_1, y_2$ satisfying $d\big(x, \mathbf{X}\big)= d(x, y_1)= d(x, y_2)$, then define $\mathfrak{P}(x)= y_1$ or $y_2$ freely). We assume 
\begin{align}
0< r\leq \hat{\rho}(x)\leq 3r\ , \quad \quad \quad \quad \forall x\in B_r(p) .\label{hat of rho assume}
\end{align}
For $0< \eta< \frac{1}{2}$, we have $\big|\frac{\hat{\rho}(y)- \hat{\rho}(x)}{\hat{\rho}(x)- \eta r}\big|\leq 4$. And we also define 
\begin{align}
\mathfrak{G}_{\rho}(x)= \big(\hat{\rho}(x), \mathfrak{P}(x)\big): B_r(p)\rightarrow \mathbb{R}\times \mathbf{X} .\nonumber 
\end{align}

For $x, y\in B_r(p)$, define 
\begin{align}
&\sigma_x(s)= \gamma_{\mathfrak{P}(x), x}(s+ \eta r)\ , \quad \quad \quad \quad \tilde{\sigma}_y(s)= \sigma_y\Big(\frac{\hat{\rho}(y)- \eta r}{\hat{\rho}(x)- \eta r}s\Big)= \gamma_{\mathfrak{P}(y), y}\Big(\frac{\hat{\rho}(y)- \eta r}{\hat{\rho}(x)- \eta r}s+ \eta r\Big)\nonumber \\
&\tau_s= \gamma_{\sigma_x(s), \tilde{\sigma}_y(s)}\ , \quad \quad \quad \quad l_s= d\big(\sigma_x(s), \tilde{\sigma}_y(s)\big) .\nonumber 
\end{align}

\begin{definition}\label{def local integral has good control}
{For $0< \eta< \frac{1}{2}$, we define 
\begin{align}
Q_{\eta, \mathbf{b}}^{r, \rho}&= \Big\{x\in B_r(p): \int_
{0}^{\hat{\rho}(x)- \eta r} |\nabla \mathbf{b}- \nabla \rho|\big(\sigma_x(s)\big)ds\leq \eta r\Big\} \nonumber \\
T_{\eta, \mathbf{b}}^{r, \rho}&= \Big\{x\in B_r(p): \fint_{B_r(p)} dy\Big(\int_
0^{\hat{\rho}(x)- \eta r} \big(\int_{\gamma_{\sigma_x(s), \tilde{\sigma}_y(s)}} |\nabla^2 \mathbf{b}|\big) ds\Big)\leq \eta r\Big\} \ , \nonumber \\
T_{\eta, \mathbf{b}}^{r, \rho}(x)&= \Big\{y\in B_r(p): \int_
0^{\hat{\rho}(x)- \eta r} \big(\int_{\gamma_{\sigma_x(s), \tilde{\sigma}_y(s)}} |\nabla^2 \mathbf{b}|\big) ds\leq \sqrt{\eta} r\Big\}\ , \quad \quad \quad \quad \forall x\in T_{\eta, \mathbf{b}}^{r, \rho} .\nonumber 
\end{align}
}
\end{definition}

For non-negative function $f$ defined on $B_r(p)$, we define 
\begin{align}
\mathfrak{Q}_{\eta, f}^{r, \rho}\vcentcolon = \Big\{x\in B_r(p): \int_
{0}^{\hat{\rho}(x)- \eta r} f\big(\sigma_x(s)\big)ds\leq \eta r\Big\} . \nonumber 
\end{align}

\begin{lemma}\label{lem lower bound of good local integra set-general}
{Assume (\ref{hat of rho assume}), then for any non-negative function $f$ satisfying $\fint_{B_{4r}(p)} f\leq \delta$, we have $\frac{V(\mathfrak{Q}_{\eta, f}^{r, \rho})}{V\big(B_r(p)\big)}\geq 1- 3^n\eta^{-n}\delta$.
}
\end{lemma}

\pf
{Firstly we have 
\begin{align}
\int_{B_r\backslash \mathfrak{Q}_{\eta, f}^{r, \rho}} dx \int_0^{\hat{\rho}(x)- \eta r} f\big(\sigma_x(s)\big) ds \geq \eta r\cdot V\big(B_r\backslash \mathfrak{Q}_{\eta, f}^{r, \rho}\big) .\nonumber 
\end{align}

Assume $\theta_s(x)$ is the gradient flow of $\rho(\cdot)$ starting from $\mathfrak{P}(x)$ at time $s+ \eta r$, using Co-Area formula and Bishop-Gromov volume comparison theorem, 
\begin{align}
&\quad \int_{B_r\backslash \mathfrak{Q}_{\eta, f}^{r, \rho}} dx \int_0^{\hat{\rho}(x)- \eta r} f\big(\sigma_x(s)\big) ds\leq \int_{B_r} dx\int_0^{\hat{\rho}(x)- \eta r} f\big(\sigma_x(s)\big) ds \nonumber \\
&= \int_{r}^{3r} dt \int_{\hat{\rho}^{-1}(t)\cap B_r} dx \int_0^{t- \eta r} f\big(\theta_s(x)\big) ds= \int_{r}^{3r} dt \int_0^{t- \eta r} ds \int_{\hat{\rho}^{-1}(t)\cap B_r}  f\big(\theta_s(x)\big) dx \nonumber \\
&\leq \big(\frac{3r}{\eta r}\big)^{n- 1} \int_{r}^{3r} dt \int_0^{t- \eta r} ds \int_{\theta_s\big(\hat{\rho}^{-1}(t)\cap B_r\big)}   f\big(\tilde{x}\big) d\tilde{x} \leq \big(\frac{3}{\eta}\big)^{n- 1} \int_{r}^{3r} dt  \int_{B_{4r}}   f\big(\tilde{x}\big) d\tilde{x} \nonumber \\
&\leq \big(\frac{3}{\eta}\big)^{n- 1}\cdot 2r\cdot \fint_{B_{4r}}   f\cdot V\big(B_{4r}\big)\leq 3^n\eta^{1- n}\delta r \cdot V\big(B_r\big) . \nonumber 
\end{align}

From the above, we obtain $\frac{V(\mathfrak{Q}_{\eta, f}^{r, \rho})}{V\big(B_r\big)}\geq 1- 3^n\eta^{-n}\delta$.
}
\qed

\begin{lemma}\label{lem lower bound of good local integra set}
{Assume (\ref{hat of rho assume}), $\mathbf{b}$ is harmonic function on $B_{16r}(p)$ satisfying $\sup_{B_{16r}(p)}|\nabla \mathbf{b}|\leq 2$ and $\fint_{B_{16r}(p)} |\nabla \mathbf{b}- \nabla \rho|\leq \delta$. Then we have 
\begin{align}
\frac{V(Q_{\eta, \mathbf{b}}^{r, \rho})}{V\big(B_r(p)\big)}&\geq 1- (12\eta^{-1})^n\delta\quad \quad and \quad \quad 
\frac{V(T_{\eta, \mathbf{b}}^{r, \rho})}{V\big(B_r(p)\big)}\geq 1- n^{20n}\eta^{-n}\sqrt{\delta}\ , \nonumber \\
\frac{V\big(T_{\eta, \mathbf{b}}^{r, \rho}(x)\big)}{V\big(B_r(p)\big)}&\geq 1- \sqrt{\eta}\ , \quad \quad \quad \quad \forall x\in T_{\eta, \mathbf{b}}^{r, \rho} .\nonumber 
\end{align}
}
\end{lemma}

\pf
{From the assumption and Bishop-Gromov Comparison Theorem, 
\begin{align}
\fint_{B_{4r}} |\nabla \mathbf{b}- \nabla \rho|\leq \frac{V\big(B_{16r}\big)}{V\big(B_{4r}\big)}\fint_{B_{16r}} |\nabla \mathbf{b}- \nabla \rho|\leq 4^n\delta  .\nonumber 
\end{align}
From the above inequality, apply Lemma \ref{lem lower bound of good local integra set-general} to $|\nabla \mathbf{b}- \nabla \rho|$, we get the first inequality of the conclusion. To prove the $3$rd inequality of the conclusion, we note 
\begin{align}
\int_{B_r\backslash T_{\eta, \mathbf{b}}^{r, \rho}(x)} dy \int_0^{\hat{\rho}(x)- \eta r}\Big(\int_{\gamma_{\sigma_x(s), \tilde{\sigma}_y(s)}} |\nabla^2 \mathbf{b}|\Big) ds\geq \sqrt{\eta} rV\big(B_r\backslash T_{\eta, \mathbf{b}}^{r, \rho}(x)\big) .\nonumber 
\end{align}

On the other hand, note $x\in T_{\eta, \mathbf{b}}^{r, \rho}$, we have 
\begin{align}
\int_{B_r\backslash T_{\eta, \mathbf{b}}^{r, \rho}(x)} dy \int_0^{\hat{\rho}(x)- \eta r}\Big(\int_{\gamma_{\sigma_x(s), \tilde{\sigma}_y(s)}} |\nabla^2 \mathbf{b}|\Big) ds &\leq V\big(B_r\big) \fint_{B_r} dy\int_0^{\hat{\rho}(x)- \eta r}\Big(\int_{\gamma_{\sigma_x(s), \tilde{\sigma}_y(s)}} |\nabla^2 \mathbf{b}|\Big) ds \nonumber \\
&\leq  V\big(B_r\big)\cdot \eta r .\nonumber 
\end{align}

Hence we obtain $\frac{V\big(T_{\eta, \mathbf{b}}^{r, \rho}(x)\big)}{V\big(B_r\big)}\geq 1- \sqrt{\eta}$. Finally we prove the $2$nd inequality. From assumption and Lemma \ref{lem gradient est imply Hessian est}, we get 
\begin{align}
\fint_{B_{8r}} |\nabla^2 \mathbf{b}| &\leq \frac{3\cdot 10^8n^5(16r)^{-1}}{2^{-\frac{3n}{2}}\cdot 2^{-4}}\cdot \Big(\fint_{B_{16r}} \big||\nabla \mathbf{b}|- 1\big|\Big)^{\frac{1}{2}} \nonumber \\
&\leq n^{14n}r^{-1} \Big(\fint_{B_{16r}} |\nabla \mathbf{b}- \nabla \rho|\Big)^{\frac{1}{2}}\leq n^{14n}r^{-1}\sqrt{\delta} .\label{L2 upper bound}
\end{align}

Note we have
\begin{align}
\int_{B_r\backslash T_{\eta, \mathbf{b}}^{r, \rho}} dx\fint_{B_r} dy \int_0^{\hat{\rho}(x)- \eta r}\Big(\int_{\gamma_{\sigma_x(s), \tilde{\sigma}_y(s)}} |\nabla^2 \mathbf{b}|\Big) ds \geq \eta r \cdot V\big(B_r\backslash T_{\eta, \mathbf{b}}^{r, \rho}\big) .\nonumber 
\end{align}

On the other hand, from  Lemma \ref{lem segment ineq}, we have 
\begin{align}
&\quad \int_{B_r\backslash T_{\eta, \mathbf{b}}^{r, \rho}} dx\fint_{B_r} dy \int_0^{\hat{\rho}(x)- \eta r}\Big(\int_{\gamma_{\sigma_x(s), \tilde{\sigma}_y(s)}} |\nabla^2 \mathbf{b}|\Big) ds \nonumber \\
& \leq \frac{1}{V\big(B_r\big)}\int_{B_r\times B_r} dxdy \int_0^{\hat{\rho}(x)- \eta r}\Big(\int_{\gamma_{\sigma_x(s), \tilde{\sigma}_y(s)}} |\nabla^2 \mathbf{b}|\Big) ds  \nonumber \\
&\leq \frac{1}{V\big(B_r\big)}\int_{r}^{3r} dt_1\int_{r}^{3r} dt_2 \int_{\big(\hat{\rho}^{-1}(t_1)\cap B_r\big)\times \big(\hat{\rho}^{-1}(t_2)\cap B_r\big)} dxdy \int_0^{t_1- \eta r}\Big(\int_{\gamma_{\sigma_x(s), \tilde{\sigma}_y(s)}} |\nabla^2 \mathbf{b}|\Big) ds \nonumber \\
&\leq \frac{(3\eta^{-1})^{n-1}}{V\big(B_r\big)}\int_{r}^{3r} dt_1\int_{r}^{3r} dt_2 \int_0^{t_1- \eta r}\Big(\int_{\theta_s\big(\hat{\rho}^{-1}(t_1)\cap B_r\big)\times \theta_{\frac{t_2- \eta r}{t_1- \eta r}\cdot s}\big(\hat{\rho}^{-1}(t_2)\cap B_r\big)} d\tilde{x}d\tilde{y} \Big(\int_{\gamma_{\tilde{x}, \tilde{y}}} |\nabla^2 \mathbf{b}|\Big) \Big) ds \nonumber \\
&\leq \frac{(3\eta^{-1})^{n-1}}{V\big(B_r\big)}\int_{r}^{3r} dt_1\int_0^{t_1- \eta r} ds \int_{\theta_s\big(\hat{\rho}^{-1}(t_1)\cap B_r\big)\times B_{4r}} d\tilde{x}d\tilde{y} \Big(\int_{\gamma_{\tilde{x}, \tilde{y}}} |\nabla^2 \mathbf{b}|\Big)  \nonumber \\
&\leq \frac{2(3\eta^{-1})^{n-1}r}{V\big(B_r\big)} \int_{B_{4r}\times B_{4r}} d\tilde{x}d\tilde{y} \Big(\int_{\gamma_{\tilde{x}, \tilde{y}}} |\nabla^2 \mathbf{b}|\Big) \nonumber \\
&\leq \frac{2(3\eta^{-1})^{n-1}r}{V\big(B_r\big)} 2^{n+ 1}\cdot (4r)\cdot V\big(B_{4r}\big) \int_{B_{8r}}|\nabla^2 \mathbf{b}| \nonumber \\
&\leq 2^{6n+ 4}(3\eta^{-1})^{n-1}r^2\cdot V\big(B_r\big)\fint_{B_{8r}}|\nabla^2 \mathbf{b}|  \nonumber \\
&\leq n^{20n}\eta^{1- n}r\sqrt{\delta} \cdot V\big(B_r\big) .\nonumber 
\end{align}
We used (\ref{L2 upper bound}) in the last inequality above.

From above, we get 
\begin{align}
\frac{V(T_{\eta, \mathbf{b}}^{r, \rho})}{V\big(B_r\big)}\geq 1- n^{20n}\eta^{-n} \sqrt{\delta} .\nonumber 
\end{align}
}
\qed

For $x, y\in B_r(p)$, for each $\mathbf{b}_i$ from Theorem \ref{thm existence of harmonic map components imply GH-dist is small} and the corresponding $\rho_i, \mathfrak{P}_i$, we define 
\begin{align}
\sigma_{x, i}(s)= \gamma_{\mathfrak{P}_i(x), x}(s+ \eta r) \ , \quad \quad \quad \quad i=1, \cdots, k .\nonumber 
\end{align}

\begin{definition}\label{def local integral has good control-2}
{For $0< \eta< \frac{1}{2}c_1$, we define 
\begin{align}
Q_{\eta, \mathbf{b}_i}^{r, \rho_i, \rho_j}&= \Big\{x\in B_r(p): \int_
{0}^{\hat{\rho}_j(x)- \eta r} \big|\nabla (\rho_i- \mathbf{b}_i)\big|\big(\sigma_{x, j}(s)\big)ds\leq \eta r\Big\} \nonumber \\
P_{\eta, \mathbf{b}_k, \mathbf{b}_l}^{r, \rho_j}&= \Big\{x\in B_r(p): \int_
{0}^{\hat{\rho}_j(x)- \eta r} \big|\langle\nabla \mathbf{b}_k, \nabla \mathbf{b}_l \rangle \big|\big(\sigma_{x, j}(s)\big)ds\leq \eta r\Big\} . \nonumber 
\end{align}
}
\end{definition}

\begin{lemma}\label{lem lower bound of good local integra set-2}
{Assume $0< r\leq \hat{\rho}_j(x)\leq 3r$ for any $x\in B_r(p)$, if 
\begin{align}
\fint_{B_{4r}(p)} \big|\langle \nabla \mathbf{b}_k, \nabla \mathbf{b}_l \rangle\big|\leq \delta \quad \quad and \quad \quad  \fint_{B_{4r}(p)}\big|\nabla (\mathbf{b}_i- \rho_i)\big| \leq \delta  ,\nonumber 
\end{align}
then we have 
\begin{align}
\frac{V(Q_{\eta, \mathbf{b}_i}^{r, \rho_i, \rho_j})}{V\big(B_r(p)\big)}\geq 1- (3\eta^{-1})^n\delta \quad \quad and \quad \quad 
\frac{V(P_{\eta, \mathbf{b}_k, \mathbf{b}_l}^{r, \rho_j})}{V\big(B_r(p)\big)}\geq 1- (3\eta^{-1})^n\delta .\nonumber 
\end{align}
}
\end{lemma}

\pf
{Apply Lemma \ref{lem lower bound of good local integra set-general} to $\big|\nabla (\mathbf{b}_i- \rho_i)\big|$ and $\Big|\big\langle\nabla \mathbf{b}_k, \nabla\mathbf{b}_l \big\rangle\Big|$ respectively, we get our conclusion.
}
\qed

\section{Quantitative almost splitting theorem}\label{SEC splitting on Ricci limit spaces}

The main results of this section were sort of implied in Cheeger-Colding's work (see \cite{CC-Ann} and \cite{CC1}), however we will not follow their argument there. Instead, we adapt the argument of Colding-Naber in \cite{CN} to prove the main result of this section. Although our argument has close relationship with \cite{CC-Ann} and \cite{CC1}, the main difference is that the angle between two segments is not involved into our argument, and the first variation formula is applied to the case of both end points are moving.

\begin{lemma}\label{lem property 2 of G-H appr}
{For any $0< \eta< \frac{1}{2}$, assume (\ref{hat of rho assume}), $\mathbf{b}$ is harmonic function on $B_{16r}(p)$ satisfying $\sup_{B_r(p)} |\nabla \mathbf{b}|\leq 2$ and $\sup_{B_r(p)}|\mathbf{b}- \rho|\leq \eta r $. Then for any $x\in T_{\eta, \mathbf{b}}^{r, \rho}\cap Q_{\eta, \mathbf{b}}^{r, \rho}$, $y\in T_{\eta, \mathbf{b}}^{r, \rho}(x)\cap Q_{\eta, \mathbf{b}}^{r, \rho}$, we have $\Big|d(x, y)- d\big(\mathfrak{G}_{\rho}(x), \mathfrak{G}_{\rho}(y)\big)\Big|\leq 5000\eta^{\frac{1}{8}}\cdot r$.
}
\end{lemma}

\pf
{Using the first variation formula for arc length, and note 
\begin{align}
\sigma_x'= \nabla \rho\ , \quad \quad and \quad \quad  \tilde{\sigma}_y'= \frac{\hat{\rho}(y)- \eta r}{\hat{\rho}(x)- \eta r}\cdot \nabla \rho , \nonumber 
\end{align} 
\begin{align}
l_t- l_0&= d\big(\sigma_x(t), \tilde{\sigma}_y(t)\big)- d\big(\sigma_x(0), \tilde{\sigma}_y(0)\big) \nonumber \\
&= \int_0^t \langle\tilde{\sigma}_y', \tau_s'\rangle\big(\tau_s(l_s)\big)- \langle\sigma_x', \tau_s'\rangle\big(\tau_s(0)\big) \nonumber \\
&= \int_0^t \langle \tilde{\sigma}_y'- \nabla \mathbf{b}, \tau_s'\rangle\big(\tau_s(l_s)\big)- \langle \sigma_x'- \nabla \mathbf{b}, \tau_s'\rangle\big(\tau_s(0)\big) \nonumber \\
&\quad + \int_0^t \langle\nabla \mathbf{b}, \tau_s'\rangle\big(\tau_s(l_s)\big)- \langle\nabla \mathbf{b}, \tau_s'\rangle\big(\tau_s(0)\big) ds \nonumber \\
&= \int_0^{t} \Big\langle \frac{\hat{\rho}(y)- \eta r}{\hat{\rho}(x)- \eta r}\nabla \rho- \nabla \mathbf{b}, \tau_s'\Big\rangle \big(\tau_s(l_s)\big)- \langle \nabla (\rho- \mathbf{b}), \tau_s'\rangle \big(\tau_s(0)\big) \nonumber \\
&\quad + \int_0^t\int_0^{l_s} \nabla^2 \mathbf{b} (\tau_s', \tau_s')\big(\tau_s(v)\big)dvds .
\label{important Hessian formula}
\end{align}

\textbf{Step (1)}. If 
\begin{align}
l_0\leq 40\eta^{\frac{1}{8}}r , \label{small assumption on l_0}
\end{align}
From (\ref{important Hessian formula}) and $x, y\in Q_{\eta, \mathbf{b}}^{r, \rho}$, we have 
\begin{align}
&\quad d(x, y)- l_0= l_{\hat{\rho}(x)- \eta r}- l_0
 \nonumber \\
&\leq \Big|\frac{\hat{\rho}(y)- \hat{\rho}(x)}{\hat{\rho}(x)- \eta r}\Big|\cdot \int_0^{\hat{\rho}(x)- \eta r} \Big\langle \nabla \rho, \tau_s'\Big\rangle \big(\tau_s(l_s)\big)ds \nonumber \\
&\quad + \int_0^{\hat{\rho}(x)- \eta r} \big|\nabla (\rho-  \mathbf{b})\big|  \big(\tilde{\sigma}_y(s)\big)+ \big|\nabla (\rho-  \mathbf{b})\big|  \big(\sigma_x(s)\big)ds \nonumber\\
&\quad + \int_0^{\hat{\rho}(x)- \eta r}\int_0^{l_s} \big|\nabla^2 \mathbf{b} \big|\big(\tau_s(v)\big)dvds \nonumber \\
&\leq \big|\hat{\rho}(y)- \hat{\rho}(x)\big|+ \frac{\hat{\rho}(x)- \eta r}{\hat{\rho}(y)- \eta r}\int_0^{\hat{\rho}(y)- \eta r} |\nabla \mathbf{b}- \nabla \rho|  \big(\sigma_y(s)\big)ds  \nonumber \\
&\quad + \int_0^{\hat{\rho}(x)- \eta r} |\nabla \mathbf{b}- \nabla \rho|  \big(\sigma_x(s)\big)ds+ \sqrt{\eta} r\nonumber \\
&\leq \big|\hat{\rho}(y)- \hat{\rho}(x)\big|+ 8\sqrt{\eta}r .\nonumber 
\end{align}
Hence 
\begin{align}
d(x, y)- \big|\hat{\rho}(y)- \hat{\rho}(x)\big|\leq l_0+ 8\sqrt{\eta}r . \label{upper bound of dist on manifold}
\end{align}

Note $d(x, y)\geq \big|\hat{\rho}(y)- \hat{\rho}(x)\big|$, by $d\big(\mathfrak{G}_{\rho}(x), \mathfrak{G}_{\rho}(y)\big)= \sqrt{d\big(\mathfrak{P}(x), \mathfrak{P}(y)\big)^2+ \big|\hat{\rho}(y)- \hat{\rho}(x)\big|^2}$, (\ref{small assumption on l_0}) and (\ref{upper bound of dist on manifold}), we have 
\begin{align}
&\quad \Big|d(x, y)- d\big(\mathfrak{G}_{\rho}(x), \mathfrak{G}_{\rho}(y)\big)\Big|
\nonumber \\
&\leq \Big|d(x, y)- \big|\hat{\rho}(y)- \hat{\rho}(x)\big| \Big|+ 
\Big|\big|\hat{\rho}(y)- \hat{\rho}(x)\big| - d\big(\mathfrak{G}_{\rho}(x), \mathfrak{G}_{\rho}(y)\big)\Big|  \nonumber \\
&\leq d(x, y)- \big|\hat{\rho}(y)- \hat{\rho}(x)\big|+ d\big(\mathfrak{P}(x), \mathfrak{P}(y)\big) \nonumber \\
&\leq l_0+ 8\sqrt{\eta}r+  l_0+ 2\sqrt{\eta}r\leq 90\eta^{\frac{1}{8}}r . \nonumber 
\end{align}

\textbf{Step (2)}. In the rest of the proof, we assume that 
\begin{align}
l_0> 40\eta^{\frac{1}{8}}r  .\label{big assumption on l_0}
\end{align}
From (\ref{important Hessian formula}) and $x, y\in Q_{\eta, \mathbf{b}}^{r, \rho}$, for $t\in [0, \hat{\rho}(x)- \eta r]$, we get 
\begin{align}
&\quad \Big|l_t- l_0- \Big(\frac{\hat{\rho}(y)- \hat{\rho}(x)}{\hat{\rho}(x)- \eta r}\Big)\cdot \int_0^{t} \langle \nabla \mathbf{b}, \tau_s'\rangle\big(\tau_s(l_s)\big)\Big| \nonumber \\
&\leq  \int_0^{t} \big|\nabla (\rho- \mathbf{b})\big|\big(\sigma_x(s)\big)+ \frac{\hat{\rho}(y)- \eta r}{\hat{\rho}(x)- \eta r} \big|\nabla (\rho- \mathbf{b})\big|\big(\tilde{\sigma}_y(s)\big) \nonumber \\
&\quad + \int_0^{t}\int_0^{l_s} \big|\nabla^2\mathbf{b}\big|\big(\tau_s(v)\big) dvds  \nonumber \\
&\leq \int_0^{t} |\nabla \mathbf{b}- \nabla \rho|\big(\sigma_x(s)\big)+\int_0^{\frac{\hat{\rho}(y)- \eta r}{\hat{\rho}(x)- \eta r}t} |\nabla \mathbf{b}- \nabla \rho|\big(\sigma_y(s)\big)ds+ \sqrt{\eta}r  \nonumber \\
&\leq 3\sqrt{\eta}r .\label{step 1 need}
\end{align}

Now we estimate $\langle \nabla \mathbf{b}, \tau_s'\rangle\big(\tau_s(l_s)\big)$. For any $0\leq t_1\leq l_s$,  we have 
\begin{align}
(\mathbf{b}\circ \tau_s)'(l_s)\leq (\mathbf{b}\circ \tau_s)'(t_1)+ \int_0^{l_s} \big|(\mathbf{b}\circ \tau_s)''(t)\big|dt ,\nonumber 
\end{align}
Take the integral of the above inequality from $0$ to $l_s$ with respect to $t_1$, we get 
\begin{align}
l_s\cdot (\mathbf{b}\circ \tau_s)'(l_s)&\leq \big[(\mathbf{b}\circ \tau_s)(l_s)- (\mathbf{b}\circ \tau_s)(0)\big]+ l_s\cdot \int_0^{l_s} |\nabla^2 \mathbf{b}|\big(\tau_s(t)\big)dt  \nonumber \\
&\leq 2\eta r+ \big[(\rho\circ \tau_s)(l_s)- (\rho\circ \tau_s)(0)\big]+ l_s\cdot \int_0^{l_s} |\nabla^2 \mathbf{b}|\big(\tau_s(t)\big)dt \nonumber \\
&= 2\eta r+ \Big(\frac{\hat{\rho}(y)- \hat{\rho}(x)}{\hat{\rho}(x)- \eta r}\Big)\cdot s+ l_s\cdot \int_0^{l_s} |\nabla^2 \mathbf{b}|\big(\tau_s(t)\big)dt .\nonumber
\end{align}
In the second inequality above we used the assumption $\sup_{B_r(p)}|\mathbf{b}- \rho|\leq \eta r$. Then 
\begin{align}
(\mathbf{b}\circ \tau_s)'(l_s)\leq \frac{s}{l_s}\cdot \Big(\frac{\hat{\rho}(y)- \hat{\rho}(x)}{\hat{\rho}(x)- \eta r}\Big)+ \frac{2\eta r}{l_s}+ \int_0^{l_s} |\nabla^2 \mathbf{b}|\big(\tau_s(t)\big)dt . \label{derivative at 0 ineq}
\end{align}

Similarly, we can also have 
\begin{align}
(\mathbf{b}\circ \tau_s)'(l_s)&\geq \frac{s}{l_s}\cdot \Big(\frac{\hat{\rho}(y)- \hat{\rho}(x)}{\hat{\rho}(x)- \eta r}\Big)- \frac{2\eta r}{l_s} - \int_0^{l_s} |\nabla^2 \mathbf{b}|\big(\tau_s(t)\big)dt .\label{derivative at 0 ineq-1}
\end{align}

Hence we have 
\begin{align}
\Big|(\mathbf{b}\circ \tau_s)'(l_s) - \frac{s}{l_s}\cdot \Big(\frac{\hat{\rho}(y)- \hat{\rho}(x)}{\hat{\rho}(x)- \eta r}\Big)\Big|\leq \frac{2\eta r}{l_s}+  \int_0^{l_s} |\nabla^2 \mathbf{b}|\big(\tau_s(t)\big)dt .\label{integrand has some Hessian bound}
\end{align}

\textbf{Step (3)}. We will show the uniform lower bound of $l_t$ when $t\in \big[0, \hat{\rho}(x)- \eta r\big]$. There are two cases to be discussed.

\begin{enumerate}
\item[\textbf{(3.A)}] If $\big|\hat{\rho}(y)- \hat{\rho}(x)\big|\leq \frac{1}{4} l_0$. From (\ref{step 1 need}) and  $|\nabla \mathbf{b}|\leq 2$, 
\begin{align}
l_t- l_0&\geq  \Big(\frac{\hat{\rho}(y)- \hat{\rho}(x)}{\hat{\rho}(x)- \eta r}\Big)\cdot \int_0^{t} \langle \nabla \mathbf{b}, \tau_s'\rangle\big(\tau_s(l_s)\big)ds- 3\sqrt{\eta}r \nonumber \\
&\geq - 2\big|\hat{\rho}(y)- \hat{\rho}(x)\big|- 3\sqrt{\eta}r
\geq -\frac{1}{2}l_0- 3\sqrt{\eta}r .\nonumber 
\end{align}
From (\ref{big assumption on l_0}), we get 
\begin{align}
l_t\geq \frac{1}{2}l_0- 3\sqrt{\eta}r\geq \eta^{\frac{1}{4}}r .\nonumber 
\end{align}

\item[\textbf{(3.B)}] If $\big|\hat{\rho}(y)- \hat{\rho}(x)\big|> \frac{1}{4} l_0$. Let $\Big(\frac{\hat{\rho}(y)- \hat{\rho}(x)}{\hat{\rho}(x)- \eta r}\Big)^2= \alpha_1$, from (\ref{derivative at 0 ineq}), (\ref{derivative at 0 ineq-1}) and (\ref{step 1 need}), we can get 
\begin{align}
l_t- l_0&\geq  \Big(\frac{\hat{\rho}(y)- \hat{\rho}(x)}{\hat{\rho}(x)- \eta r}\Big)\cdot \int_0^{t} \langle \nabla \mathbf{b}, \tau_s'\rangle\big(\tau_s(l_s)\big)ds- 3\sqrt{\eta}r\nonumber \\
&\geq \alpha_1\int_0^t \frac{s}{l_s}ds- \sqrt{\alpha_1}\cdot \int_0^{t} \Big(\frac{2\eta r}{l_s}+  \int_0^{l_s} |\nabla^2 \mathbf{b}|\big(\tau_s(t)\big)dt \Big)ds- 3\sqrt{\eta}r \nonumber \\
&\geq -\sqrt{\alpha_1}\cdot 2\eta r\int_0^{\frac{2\eta r}{\sqrt{\alpha_1}}} \frac{1}{l_s}ds- 16\int_0^{\hat{\rho}(x)- \eta r} \int_0^{l_s} |\nabla^2 \mathbf{b}|\big(\tau_s(t)\big)dt ds- 3\sqrt{\eta}r .\label{4B-1}
\end{align}
Note for any $0\leq s\leq \frac{2\eta r}{\sqrt{\alpha_1}}$, using $\sqrt{\alpha_1}\geq 2^{-2}\cdot \frac{l_0}{\hat{\rho}(x)- \eta r}\geq \frac{l_0}{12r}$ and (\ref{big assumption on l_0}), we have 
\begin{align}
l_s&\geq l_0- s- \frac{\hat{\rho}(y)- \eta r}{\hat{\rho}(x)- \eta r}s\geq l_0- \big(1+ \frac{3r- \eta r}{r- \eta r}\big)\frac{2\eta r}{\sqrt{\alpha_1}}\geq l_0- 7\cdot 2\eta r\cdot \frac{12r}{l_0} \nonumber \\
&\geq l_0- \frac{200\eta r^2}{l_0}\geq \sqrt{\eta} r .\label{4B-2}
\end{align}

From (\ref{4B-1}) and (\ref{4B-2}), we get 
\begin{align}
l_t- l_0\geq -\sqrt{\alpha_1}\cdot 2\eta r\cdot \frac{2\eta r}{\sqrt{\alpha_1}}\cdot \frac{1}{\sqrt{\eta} r}- 19\sqrt{\eta}r \geq  -25\sqrt{\eta}r .\nonumber 
\end{align}
Then by (\ref{big assumption on l_0}) again,
\begin{align}
l_t\geq l_0- 25\sqrt{\eta}r \geq \eta^{\frac{1}{4}}r .\nonumber 
\end{align}
\end{enumerate}

From above two cases, we always have
\begin{align}
l_t\geq \eta^{\frac{1}{4}}r \ , \quad \quad \quad \quad \forall t\in \big[0, \hat{\rho}(x)- \eta r\big] .\label{uniform lower bound of l_t}
\end{align}

\textbf{Step (4)}. From (\ref{integrand has some Hessian bound}), (\ref{uniform lower bound of l_t}) and $\eta\geq \delta$, we obtain
\begin{align}
&\quad \Big|\Big(\frac{\hat{\rho}(y)- \hat{\rho}(x)}{\hat{\rho}(x)- \eta r}\Big)\cdot \int_0^{t} \langle \nabla \mathbf{b}, \tau_s'\rangle\big(\tau_s(l_s)\big)- \alpha_1\int_0^t \frac{s}{l_s} ds\Big|\nonumber \\
&\leq \Big|\frac{\hat{\rho}(y)- \hat{\rho}(x)}{\hat{\rho}(x)- \eta r}\Big|\int_0^{t} \Big(\frac{2\eta r}{l_s}+  \int_0^{l_s} |\nabla^2 \mathbf{b}|\big(\tau_s(t)\big)dt \Big)ds\nonumber \\
&\leq 3r\cdot \frac{8\eta}{\eta^{\frac{1}{4}}}+ 4\int_0^{\hat{\rho}(x)- \eta r} \int_0^{l_s} |\nabla^2 \mathbf{b}|\big(\tau_s(t)\big)dt ds
\leq 28\sqrt{\eta}r .\label{step 2 need}
\end{align}

From (\ref{step 1 need}) and (\ref{step 2 need}), we get
\begin{align}
\Big|l_t- l_0- \alpha_1\int_0^t \frac{s}{l_s} ds\Big| \leq 31\sqrt{\eta}r . \nonumber 
\end{align}

Define $\mathfrak{L}(t)= \int_0^t \frac{s}{l_s}ds$, then $l_t= \frac{t}{\mathfrak{L}'(t)}$ and we have 
\begin{align}
\Big|\frac{t}{\mathfrak{L}'(t)}- \alpha_1\mathfrak{L}(t)- l_0\Big|\leq 31\sqrt{\eta}r .\label{ode ineq 1}
\end{align}

On the other hand, let $f(s)= \sqrt{\alpha_1s^2+ l_0^2}$, and define 
\begin{align}
F(t)= \int_0^t \frac{s}{f(s)}ds .\nonumber 
\end{align}
Then it is easy to get 
\begin{align}
\frac{t}{F(t)'}- \alpha_1 F(t)- l_0= 0 .\label{ode ineq 2}
\end{align}

From (\ref{ode ineq 1}) and (\ref{ode ineq 2}), we obtain
\begin{align}
\Big|\Big[\big(\frac{1}{2}\alpha_1\mathfrak{L}^2+ l_0\mathfrak{L}\big)- \big(\frac{1}{2}\alpha_1F^2+ l_0F\big)\Big]'\Big|&\leq \big|31\sqrt{\eta}r\cdot \mathfrak{L}'\big|= \big|31\sqrt{\eta}r \cdot \frac{t}{l_t}\big| \nonumber \\
&\leq \big|31\sqrt{\eta}r \cdot \frac{3r}{\eta^{\frac{1}{4}}r}\big|\leq 100\eta^{\frac{1}{4}}r .\nonumber
\end{align}
Take the integral of the above inequality, also note $\mathfrak{L}(0)= F(0)= 0$ and $0\leq t\leq \hat{\rho}(x)- \eta r\leq 3r$, we have 
\begin{align}
\Big|\Big(\frac{1}{2}\alpha_1\mathfrak{L}(t)^2+ l_0\mathfrak{L}(t)\Big)- \Big(\frac{1}{2}\alpha_1F(t)^2+ l_0F(t)\Big)\Big|\leq 300\eta^{\frac{1}{4}}r^2 .\nonumber 
\end{align}

Simplify the above inequality, note $\mathfrak{L}(t)\geq 0$ and $F(t)\geq 0$, using $l_0\geq \eta^{\frac{1}{8}}r$, 
\begin{align}
\big|\mathfrak{L}(t)- F(t)\big|&\leq l_0^{-1}\Big|\Big(\frac{1}{2}\alpha_1\mathfrak{L}(t)^2+ l_0\mathfrak{L}(t)\Big)- \Big(\frac{1}{2}\alpha_1F(t)^2+ l_0F(t)\Big)\Big| \nonumber \\
&\leq 300\eta^{\frac{1}{8}}r .\label{ode ineq 3}
\end{align}

From (\ref{ode ineq 1}), (\ref{ode ineq 2}) and (\ref{ode ineq 3}), we have 
\begin{align}
\big|l_t- f(t)\big|=\Big|\frac{t}{\mathfrak{L}'(t)}- \frac{t}{F(t)'}\Big|\leq \alpha_1\cdot 300\eta^{\frac{1}{8}}r+ 31\sqrt{\eta}r \leq 4831\eta^{\frac{1}{8}}r .\label{ode ineq 4}
\end{align}

Let $t= \hat{\rho}(x)- \eta r$ in (\ref{ode ineq 4}), note $\Big|l_0- d\big(\mathfrak{P}(x), \mathfrak{P}(y)\big)\Big|\leq 2\eta r$ and $l_{\hat{\rho}(x)- \eta r}= d(x, y)$, 
\begin{align}
&\quad \Big|d(x, y)- d\big(\mathfrak{G}_{\rho}(x), \mathfrak{G}_{\rho}(y)\big)\Big|\nonumber \\
&\leq  \big|f(\hat{\rho}(x)- \eta r)- l_{\hat{\rho}(x)- \eta r}\big|+ \Big|d\big(\mathfrak{G}_{\rho}(x), \mathfrak{G}_{\rho}(y)\big)- f(\hat{\rho}(x)- \eta r)\Big| \nonumber \\
&\leq 4831\eta^{\frac{1}{8}}r+ 
\Big|\sqrt{d\big(\mathfrak{P}(x), \mathfrak{P}(y)\big)^2+ \big(\hat{\rho}(y)- \hat{\rho}(x)\big)^2}- \sqrt{l_0^2+ \big(\hat{\rho}(y)- \hat{\rho}(x)\big)^2}\Big| \nonumber \\
&\leq 4831\eta^{\frac{1}{8}}r+ 2\eta r\leq 5000\eta^{\frac{1}{8}}r . \nonumber 
\end{align}
}
\qed

\begin{cor}\label{cor dist of points in geodesic balls}
{For any $\eta\in (0, 1)$, assume (\ref{hat of rho assume}) and 
\begin{align}
&\sup_{B_{16r}(p)} |\nabla \mathbf{b}|\leq 2\ ,  \quad \quad \quad \quad \sup_{B_r(p)}|\mathbf{b}- \rho|\leq \eta r \nonumber \\
&\fint_{B_{16r}(p)} |\nabla \mathbf{b}- \nabla \rho| \leq 2^{-n}\eta^{3n} , \nonumber  
\end{align}
then there exists $\delta_1= n^{21}\eta^{\frac{1}{2n}}$ such that $\sup\limits_{x, y\in T_{\eta, \mathbf{b}}^{r, \rho}\cap Q_{\eta, \mathbf{b}}^{r, \rho}\cap B_{(1- \delta_1)r}(p)}\Big|d(x, y)- d\big(\mathfrak{G}_{\rho}(x), \mathfrak{G}_{\rho}(y)\big)\Big|\leq n^{22}\eta^{\frac{1}{3n}}r$.
}
\end{cor}

\pf
{From Lemma \ref{lem lower bound of good local integra set}, for any $x, y\in T_{\eta, \mathbf{b}}^{r, \rho}\cap Q_{\eta, \mathbf{b}}^{r, \rho}\cap B_{(1- \delta_1)r}(p)$, we have
\begin{align}
&\quad V\Big(B_r(p)\backslash \big(T_{\eta, \mathbf{b}}^{r, \rho}\cap Q_{\eta, \mathbf{b}}^{r, \rho}\cap T_{\eta, \mathbf{b}}^{r, \rho}(x)\cap T_{\eta, \mathbf{b}}^{r, \rho}(y)\big)\Big) \nonumber \\
&= V\Big(\big(B_r(p)\backslash T_{\eta, \mathbf{b}}^{r, \rho}\big)\cup \big(B_r(p)\backslash Q_{\eta, \mathbf{b}}^{r, \rho}\big)\cup \big(B_r(p)\backslash T_{\eta, \mathbf{b}}^{r, \rho}(x)\big)\cup \big(B_r(p)\backslash T_{\eta, \mathbf{b}}^{r, \rho}(y)\big)\Big) \nonumber \\
&\leq V\big(B_r(p)\backslash T_{\eta, \mathbf{b}}^{r, \rho}\big)+ V\big(B_r(p)\backslash Q_{\eta, \mathbf{b}}^{r, \rho}\big)+ V\big(B_r(p)\backslash T_{\eta, \mathbf{b}}^{r, \rho}(x)\big)+ V\big(B_r(p)\backslash T_{\eta, \mathbf{b}}^{r, \rho}(y)\big) \nonumber \\
&\leq 2\big(n^{20n}\eta^{-n}\sqrt{2^{-n}\eta^{3n}}+ \sqrt{\eta}\big)V\big(B_r(p)\big)\leq n^{20n}\eta^{\frac{1}{2}}V\big(B_r(p)\big) .\nonumber 
\end{align}

We claim that $T_{\eta, \mathbf{b}}^{r, \rho}\cap Q_{\eta, \mathbf{b}}^{r, \rho}\cap T_{\eta, \mathbf{b}}^{r, \rho}(x)\cap T_{\eta, \mathbf{b}}^{r, \rho}(y)\cap B_{\delta_1 r}(x)\neq \emptyset$, otherwise
\begin{align}
B_{\delta_1 r}(x)\subset B_r(p)\backslash \big(T_{\eta, \mathbf{b}}^{r, \rho}\cap Q_{\eta, \mathbf{b}}^{r, \rho}\cap T_{\eta, \mathbf{b}}^{r, \rho}(x)\cap T_{\eta, \mathbf{b}}^{r, \rho}(y)\big) ,\nonumber 
\end{align}
which implies 
\begin{align}
\frac{V\big(B_{\delta_1 r}(x)\big)}{V\big(B_r(p)\big) }\leq \frac{V\Big(B_r(p)\backslash \big(T_{\eta, \mathbf{b}}^{r, \rho}\cap Q_{\eta, \mathbf{b}}^{r, \rho}\cap T_{\eta, \mathbf{b}}^{r, \rho}(x)\cap T_{\eta, \mathbf{b}}^{r, \rho}(y)\big)\Big)}{V\big(B_r(p)\big)} \leq n^{20n}\eta^{\frac{1}{2}} . \label{cont 0}
\end{align}
On the other hand, from Bishop-Gromov Comparison Theorem and $B_r(p)\subset B_{2r}(x)$, using the definition of $\delta_1$, we have 
\begin{align}
\frac{V\big(B_{\delta_1 r}(x)\big)}{V\big(B_r(p)\big) }\geq \frac{V\big(B_{\delta_1 r}(x)\big)}{V\big(B_{2r}(x)\big)}\geq \big(\frac{\delta_1}{2}\big)^n> n^{20n}\eta^{\frac{1}{2}} ,\nonumber 
\end{align}
which is the contradiction. Hence we can find
\begin{align}
z\in T_{\eta, \mathbf{b}}^{r, \rho}\cap Q_{\eta, \mathbf{b}}^{r, \rho}\cap T_{\eta, \mathbf{b}}^{r, \rho}(x)\cap T_{\eta, \mathbf{b}}^{r, \rho}(y)\cap B_{\delta_1 r}(x) .\nonumber 
\end{align}

Now apply Lemma \ref{lem property 2 of G-H appr} to $y, z$ and $x,  z$ respectively, then we have 
\begin{align}
&\Big|d(x, y)- d\big(\mathfrak{G}_{\rho}(x), \mathfrak{G}_{\rho}(y)\big)\Big|\leq \big|d(x, y)- d(z, y)\big|+ \Big|d(z, y)- d\big(\mathfrak{G}_{\rho}(z), \mathfrak{G}_{\rho}(y)\big)\Big| \nonumber \\
&\quad \quad + \Big|d\big(\mathfrak{G}_{\rho}(z), \mathfrak{G}_{\rho}(y)\big)- d\big(\mathfrak{G}_{\rho}(x), \mathfrak{G}_{\rho}(y)\big)\Big| \nonumber \\
&\leq d(x, z)+ 5000\eta^{\frac{1}{8}} r + d\big(\mathfrak{G}_{\rho}(x), \mathfrak{G}_{\rho}(z)\big) \nonumber \\
&\leq 10^4\eta^{\frac{1}{8}} r + 2d(x, z)\leq 10^4\eta^{\frac{1}{8}} r+ 2\delta_1 r \nonumber \\
&\leq n^{22}\eta^{\frac{1}{3n}}r . \nonumber 
\end{align} 
}
\qed

\pf [of Theorem \ref{thm existence of harmonic map components imply GH-dist is small}]
{\textbf{Step (1)}. We firstly deal with the case $k= 1$. From Proposition \ref{prop existence of almost linear function}, there exists $\delta\leq  2^{-100n^2}\epsilon_1^{2n^2}$, where $\epsilon_1> 0$ is a constant to be determined later, such that if (\ref{assump on harmonic linear function}) holds for $\delta$, we can find two functions $\rho_1, \tilde{\rho}_1$ satisfying the following
\begin{align}
&\rho_1(x)= \hat{\rho}_1(x)+ t_0= d\big(x, \rho^{-1}(t_0)\big)+ t_0 \ , \quad  \quad \quad \tilde{\rho}_1(x)= t_1- d\big(x, \tilde{\rho}^{-1}(t_1)\big) \label{rho-1 and tilde rho-1} \\
&\frac{r}{320}\leq \hat{\rho}_1(x)\leq  \frac{3r}{320}\ , \quad \quad \quad \quad \forall x\in B_{\frac{r}{320}}(p) \label{hat of rho assume-1} \\
&\sup_{B_{\frac{1}{20}r}(p)}\big|\mathbf{b}_1- \rho_1\big|\leq \epsilon_1\cdot r \quad \quad and \quad \quad \fint_{B_{\frac{1}{20}r}(p)} \big|\nabla (\mathbf{b}_1- \rho_1)\big|\leq \epsilon_1 \label{assumption 1 for almost split} \\
&\rho_1(x)\leq \tilde{\rho}_1(x)+ \frac{\epsilon_1 r}{2} \ , \quad \quad\quad \quad \forall \ x\in B_{\frac{r}{10}}(p) .\label{pinch ineq}
\end{align}

From (\ref{hat of rho assume-1}), (\ref{assump on harmonic linear function}) and (\ref{assumption 1 for almost split}), for $\eta\in (0, 1)$ to be determined later, if we assume 
\begin{align}
\epsilon_1\leq 2^{-n}\eta^{3n} ,\label{epsilon_1 need}
\end{align}
apply Corollary \ref{cor dist of points in geodesic balls}, there exist $\delta_1= n^{21}\eta^{\frac{1}{2n}}$ such that 
\begin{align}
&\quad \sup_{x, y\in T_{\eta, \mathbf{b}_1}^{\frac{r}{320}, \rho_1}\cap Q_{\eta, \mathbf{b}_1}^{\frac{r}{320}, \rho_1}\cap B_{\frac{(1- \delta_1)r}{320}}(p)}\Big|d(x, y)- d\big(\mathfrak{G}_{\rho_1}(x), \mathfrak{G}_{\rho_1}(y)\big)\Big|\leq \frac{n^{22}}{320}\eta^{\frac{1}{3n}}r .\label{dist inequality need 1} 
\end{align}

Now choose a maximal collection of disjoint balls with radius $\frac{\delta_1 r}{320}$ centered at $B_{\frac{(1- \delta_1)r}{320}}(p)$, denoted as $\Big\{B_{\frac{\delta_1 r}{320}}(x_i)\Big\}_{i= 1}^m$, then 
\begin{align}
B_{\frac{r}{320}}(p)\subset  \bigcup_{i= 1}^m B_{\frac{3\delta_1 r}{320}}(x_i) \quad \quad and \quad \quad \bigcup_{i= 1}^m B_{\frac{\delta_1 r}{320}}(x_i)\subset B_{\frac{r}{320}}(p) . \label{cover ineq 1} 
\end{align}

From Lemma \ref{lem lower bound of good local integra set}, (\ref{hat of rho assume-1}), (\ref{epsilon_1 need}) and (\ref{assumption 1 for almost split}), we have 
\begin{align}
\frac{V\big(B_{\frac{r}{320}}(p)\backslash (T_{\eta, \mathbf{b}_1}^{\frac{r}{320}, \rho_1}\cap Q_{\eta, \mathbf{b}_1}^{\frac{r}{320}, \rho_1})\big)}{V\big(B_{\frac{r}{320}}(p)\big)}&\leq \frac{V\big(B_{\frac{r}{320}}(p)\backslash T_{\eta, \mathbf{b}_1}^{\frac{r}{320}, \rho_1}\big)}{V\big(B_{\frac{r}{320}}(p)\big)}+ \frac{V\big(B_{\frac{r}{320}}(p)\backslash Q_{\eta, \mathbf{b}_1}^{\frac{r}{320}, \rho_1}\big)}{V\big(B_{\frac{r}{320}}(p)\big)} \nonumber \\
&\leq (12\eta^{-1})^n\epsilon_1+ n^{20}\eta^{-n}\sqrt{\epsilon_1}\leq n^{20n}\eta^{\frac{n}{2}}  .\label{volume ineq 1} 
\end{align}
On the other hand, from Bishop-Gromov Comparison Theorem and (\ref{cover ineq 1}),
\begin{align}
\frac{V\big(B_{\frac{\delta_1 r}{320}}(x_i)\big)}{V\big(B_{\frac{r}{320}}(p)\big)}\geq \frac{V\big(B_{\frac{\delta_1 r}{320}}(x_i)\big)}{V\big(B_{\frac{2r}{320}}(x_i)\big)}\geq \big(\frac{\delta_1}{2}\big)^n> n^{20n}\eta^{\frac{n}{2}}  .\label{volume ineq 2}
\end{align} 

From (\ref{volume ineq 1}) and (\ref{volume ineq 2}), we get 
$V\big(B_{\frac{\delta_1 r}{320}}(x_i)\big)> V\big(B_{\frac{r}{320}}(p)\backslash (T_{\eta, \mathbf{b}_1}^{\frac{r}{320}, \rho_1}\cap Q_{\eta, \mathbf{b}_1}^{\frac{r}{320}, \rho_1})\big)$. So there exists $y_i$ such that  $y_i\in B_{\frac{\delta_1 r}{320}}(x_i)\cap T_{\eta, \mathbf{b}_1}^{\frac{r}{320}, \rho_1}\cap Q_{\eta, \mathbf{b}_1}^{\frac{r}{320}, \rho_1}$, combining (\ref{cover ineq 1}) yields 
\begin{align}
B_{\frac{r}{320}}(p)\subset  \bigcup_{i= 1}^m B_{\frac{\delta_1 r}{80}}(y_i) . \label{good cover by good points}
\end{align}

Now we define $\mathbf{Y}= \{y_1, \cdots, y_m\}\cap B_{\frac{r}{1280}- n^{22}\eta^{\frac{1}{3n}}r}(p)$ and $\mathbf{X}_1= \bigcup\limits_{y_i\in \mathbf{Y}} \mathfrak{P}_1(y_i)\subset \mathbf{X}_{\mathbf{b}_1}$, which is a metric space with distance function $d$ inherited from the metric $g$ of the Riemannian manifold $M^n$.

\textbf{Step (2)}. For any $x\in B_{\frac{r}{1280}}(p)$, we define 
\begin{align}
\mathfrak{g}_1(x)=\Big(\hat{\rho}_1\big(\mathfrak{n}_1(x)\big)- \hat{\rho}_1\big(\mathfrak{n}_1(p)\big), \mathfrak{P}_1\big(\mathfrak{n}_1(x)\big)\Big): B_{\frac{r}{1280}}(p)\rightarrow B_{\frac{r}{1280}}(0, \hat{p})\subset \mathbb{R}\times \mathbf{X}_1 ,\nonumber 
\end{align}
where $\hat{p}= \mathfrak{P}_1\big(\mathfrak{n}_1(p)\big)$ and 
\begin{align}
\mathfrak{n}_1(x)\in \mathcal{D}_x\vcentcolon = \big\{y_i\in \mathbf{Y}| d(y_i, x)= \min_{y\in \mathbf{Y}}d(y, x)\big\} .
\end{align}
If $\mathcal{D}_x$ contains more than one elements, we choose $\mathfrak{n}_1(x)$ from it freely. 

If $x= y_i$ for some $i$, then $\mathfrak{n}_1(x)= x$. If $x\neq y_i$ for any $i$, then there is $\tilde{x}$ such that
\begin{align}
d(\tilde{x}, x)< \frac{\delta_1 r}{80}+ n^{22}\eta^{\frac{1}{3n}}r \quad \quad and \quad \quad d(\tilde{x}, p)< \frac{r}{1280}- \frac{\delta_1 r}{80}- n^{22}\eta^{\frac{1}{3n}}r .\nonumber 
\end{align}

From  (\ref{good cover by good points}), there exists $y_i$ such that $\tilde{x}\in B_{\frac{\delta_1 r}{80}}(y_i)$, then $d(y_i, x)< \frac{\delta_1 r}{40}$ and 
\begin{align}
d(y_i, p)\leq d(y_i, \tilde{x})+ d(\tilde{x}, p)< \frac{r}{1280}- n^{22}\eta^{\frac{1}{3n}}r ,\nonumber 
\end{align}
which implies $y_i\in \mathbf{Y}$. From the definition of $\mathfrak{n}_1(x)$,
\begin{align}
d\big(\mathfrak{n}_1(x), x\big)\leq d(x, y_i)< \frac{\delta_1 r}{40} .\label{diff n(x) and x}
\end{align} 

From the definition of $\mathfrak{g}_1$ and (\ref{dist inequality need 1}), we have 
\begin{align}
d\big(\mathfrak{g}_1(x), \mathfrak{g}_1(p)\big)&= d\Big(\mathfrak{G}_{\rho_1}\big(\mathfrak{n}_1(x)\big), \mathfrak{G}_{\rho_1}\big(\mathfrak{n}_1(p)\big)\Big)\leq d\big(\mathfrak{n}_1(x), \mathfrak{n}_1(p)\big)+ \frac{n^{22}}{160}\eta^{\frac{1}{3n}}r \nonumber \\
&\leq d\big(\mathfrak{n}_1(x), p\big)+ d\big(\mathfrak{n}_1(p), p\big)+ \frac{n^{22}}{160}\eta^{\frac{1}{3n}}r \nonumber \\
&\leq \frac{r}{1280}- n^{22}\eta^{\frac{1}{3n}}r+ \frac{\delta_1 r}{40}+ \frac{n^{22}}{160}\eta^{\frac{1}{3n}}r \nonumber \\
&< \frac{r}{1280} ,\nonumber  
\end{align}
which implies that $\mathfrak{g}_1\big(B_{\frac{r}{1280}}(p)\big)\subset B_{\frac{r}{1280}}(0, \hat{p})$.

Now for any $x_1, x_2\in B_{\frac{r}{1280}}(p)$, from (\ref{dist inequality need 1}) and (\ref{diff n(x) and x}), we have 
\begin{align}
&\quad \Big|d\big(\mathfrak{g}_1(x_1), \mathfrak{g}_1(x_2)\big)- d(x_1, x_2)\Big| \nonumber \\
&= \Big|d\big(\mathfrak{G}_{\rho_1}\circ\mathfrak{n}_1(x_1), \mathfrak{G}_{\rho_1}\circ\mathfrak{n}_1(x_2)\big)- d(x_1, x_2)\Big| \nonumber \\
&\leq \Big|d\big(\mathfrak{G}_{\rho_1}\circ\mathfrak{n}_1(x_1), \mathfrak{G}_{\rho_1}\circ\mathfrak{n}_1(x_2)\big)- d\big(\mathfrak{n}_1(x_1), \mathfrak{n}_1(x_2)\big)\Big|+ \frac{\delta_1 r}{20} \nonumber \\
&\leq \frac{n^{22}}{160}\eta^{\frac{1}{3n}}r+ \frac{\delta_1 r}{20}\leq n^{22}\eta^{\frac{1}{3n}}r .\label{g-close}
\end{align}

Now to show that $\mathfrak{g}_1: B_{\frac{1}{1280}r}(p)\rightarrow B_{\frac{1}{1280}r}(0, \hat{p})$ is an pointed $(\frac{\epsilon r}{3})$-Gromov-Hausdorff approximation, we only need to show that 
\begin{align}
B_{\frac{r}{1280}}(0, \hat{p})\subset \mathbf{U}_{\frac{\epsilon r}{3}}\Big(\mathfrak{g}_1\big(B_{\frac{r}{1280}}(p)\big)\Big) . \label{property 2 of G-H} 
\end{align}

\textbf{Step (3)}. For any $(t, \hat{x})\in B_{\frac{r}{1280}}(0, \hat{p})\subset \mathbb{R}\times \mathbf{X}_1$, there is $y_i\in \mathbf{Y}$ such that $\hat{x}= \mathfrak{P}_1(y_i)$. 

If $t+ \hat{\rho}_1\big(\mathfrak{n}_1(p)\big)> \hat{\rho}_1(y_i)$, then 
\begin{align}
t+ \hat{\rho}_1\big(\mathfrak{n}_1(p)\big)- \hat{\rho}_1(y_i)\leq \frac{r}{1280}+ 2\cdot \frac{r}{320}< \frac{r}{80} . \label{vertical need 1}
\end{align}
For $t_1= \tilde{\rho}_1(p)+ \frac{r}{20}$, note 
\begin{align}
\tilde{\rho}_1(y_i)\leq \tilde{\rho}_1(p)+ d(y_i, p)< \tilde{\rho}_1(p)+ \frac{r}{1280}< t_1- \Big(t- \hat{\rho}_1\big(\mathfrak{n}_1(p)\big)- \hat{\rho}_1(y_i)\Big) .\nonumber 
\end{align}
There is $\tilde{y}\in \tilde{\rho}_1^{-1}(t_1)$ such that $d\big(y_i, \tilde{\rho}_1^{-1}(t_1)\big)= d(y_i, \tilde{y})$. Define 
\begin{align}
z= \gamma_{y_i, \tilde{y}}\big(t+ \hat{\rho}_1\big(\mathfrak{n}_1(p)\big)- \hat{\rho}_1(y_i)\big) ,\nonumber 
\end{align}
then 
\begin{align}
d(z, y_i)= t+ \hat{\rho}_1\big(\mathfrak{n}_1(p)\big)- \hat{\rho}_1(y_i) . \label{dist between z and y_i}
\end{align}
And from (\ref{rho-1 and tilde rho-1}), we have 
\begin{align}
\tilde{\rho}_1(z)= \tilde{\rho}_1(y_i)+ t+ \hat{\rho}_1\big(\mathfrak{n}_1(p)\big)- \hat{\rho}_1(y_i) .\label{lower barrier linear equation}
\end{align}

Using (\ref{lower barrier linear equation}) and (\ref{pinch ineq}), we get   
\begin{align}
\big|\hat{\rho}_1(z)- \hat{\rho}_1\big(\mathfrak{n}_1(p)\big)- t\big|&= \big|\rho_1(z)- t_0- \hat{\rho}_1\big(\mathfrak{n}_1(p)\big)- t\big|\leq \big|\tilde{\rho}_1(z)- t_0- \hat{\rho}_1\big(\mathfrak{n}_1(p)\big)- t\big|+ \epsilon_1 r \nonumber \\
&= \big|\tilde{\rho}_1(y_i)- \hat{\rho}_1(y_i)- t_0\big|+ \epsilon_1 r \leq  \big|\rho_1(y_i)- \hat{\rho}_1(y_i)- t_0\big|+ 2\epsilon_1 r = 2\epsilon_1 r .\label{verctical differ-1.1}
\end{align}

Note 
\begin{align}
d(z, p)&\leq d(p, y_i)+ d(y_i, z)\leq \frac{r}{1280}+ \big(t+ \hat{\rho}_1\big(\mathfrak{n}_1(p)\big)- \hat{\rho}_1(y_i)\big) \nonumber \\
&< \frac{r}{1280}+ \big(\frac{r}{1280}+ \frac{r}{1280}\big) \leq \frac{r}{320} .\nonumber 
\end{align}
From (\ref{good cover by good points}), there is $y_{i_1}$ such that $d(z, y_{i_1})< \frac{\delta_1}{80}r$, using (\ref{verctical differ-1.1}) and (\ref{dist between z and y_i}),  
\begin{align}
d\big(\mathfrak{P}_1(y_{i_1}), \hat{x}\big)&= \sqrt{d\big(\mathfrak{G}_{\rho_1}(y_{i_1}), \mathfrak{G}_{\rho_1}(y_i)\big)^2- \big|\hat{\rho}_1(y_{i_1})- \hat{\rho}_1(y_i)\big|^2} \nonumber \\
&\leq \sqrt{d(y_{i_1}, y_i)^2- \big|\hat{\rho}_1(z)- \hat{\rho}_1(y_i)\big|^2}+ \frac{n^{22}}{160}\eta^{\frac{1}{3n}}r+ \frac{\delta_1}{80} r \nonumber \\
&\leq \sqrt{d(z, y_i)^2- \big|t+ \hat{\rho}_1\big(\mathfrak{n}_1(p)\big)- \hat{\rho}_1(y_i)\big|^2}+ n^{22}\eta^{\frac{1}{3n}}r= n^{22}\eta^{\frac{1}{3n}}r .\label{horizontal differ-1.2}
\end{align}

Using (\ref{diff n(x) and x}), (\ref{dist inequality need 1}), (\ref{verctical differ-1.1}) and (\ref{horizontal differ-1.2}), we have 
\begin{align}
d(z, p)&\leq d\big(y_{i_1}, \mathfrak{n}_1(p)\big)+ \frac{\delta_1 r}{20} \leq d\Big(\mathfrak{G}_{\rho_1}(y_{i_1}), \mathfrak{G}_{\rho_1}\big(\mathfrak{n}_1(p)\big)\Big)+ \frac{\delta_1 r}{20}+ \frac{n^{22}}{160}\eta^{\frac{1}{3n}}r \nonumber \\
&\leq \sqrt{\big|\hat{\rho}_1(y_{i_1})- \hat{\rho}_1\big(\mathfrak{n}_1(p)\big)\big|^2+ d\big(\mathfrak{P}_1(y_{i_1}), \hat{p}\big)^2}+ \frac{\delta_1 r}{20}+ \frac{n^{22}}{160}\eta^{\frac{1}{3n}}r \nonumber \\
&\leq \sqrt{\big|\hat{\rho}_1(z)- \hat{\rho}_1\big(\mathfrak{n}_1(p)\big)\big|^2+ d\big(\hat{x}, \hat{p}\big)^2}+ \frac{\delta_1 r}{10}+ 2\cdot n^{22}\eta^{\frac{1}{3n}}r \nonumber \\
&\leq \sqrt{t^2+ d\big(\hat{x}, \hat{p}\big)^2}+ n^{23}\eta^{\frac{1}{3n}}r \nonumber \\
&< \frac{r}{1280}+ n^{23}\eta^{\frac{1}{3n}}r .\label{7.45.1}
\end{align}
Hence we can find $y$ such that 
\begin{align}
d(z, y)\leq n^{23}\eta^{\frac{1}{3n}}r+ \delta_1 r \quad \quad and \quad \quad d(y, p)< \frac{r}{1280}- \delta_1 r . \nonumber 
\end{align}

From (\ref{good cover by good points}), there is $y_{i_2}$ such that $y\in B_{\frac{\delta_1 r}{80}}(y_{i_2})$, then $y_{i_2}\in B_{\frac{r}{1280}}(p)$ and 
\begin{align}
d(z, y_{i_2})< n^{23}\eta^{\frac{1}{3n}}r+ \delta_1 r+ \frac{\delta_1 r}{80}\leq n^{24}\eta^{\frac{1}{3n}}r .\nonumber  
\end{align}
 
Now using (\ref{verctical differ-1.1}), similar as (\ref{horizontal differ-1.2}), we get 
\begin{align}
d\big(\mathfrak{g}_1(y_{i_2}), (t, \hat{x})\big)&= d\big(\mathfrak{G}_{\rho_1}(y_{i_2}), (t+ \hat{\rho}_1\big(\mathfrak{n}_1(p)\big), \hat{x})\big) \nonumber \\
&\leq \big|\hat{\rho}_1(y_{i_2})- t- \hat{\rho}_1\big(\mathfrak{n}_1(p)\big)\big|+ d\big(\mathfrak{P}_1(y_{i_2}), \hat{x}\big) \nonumber \\
&\leq \big|\hat{\rho}_1(z)- t- \hat{\rho}_1\big(\mathfrak{n}_1(p)\big)\big|+ 2\cdot n^{24}\eta^{\frac{1}{3n}}r\leq n^{25}\eta^{\frac{1}{3n}}r . \label{dense-1}
\end{align}

If $t+ \hat{\rho}_1\big(\mathfrak{n}_1(p)\big)\leq \hat{\rho}_1(y_i)$, we get $\mathfrak{G}_{\rho_1}\Big(\gamma_{\mathfrak{P}_1(y_i), y_i}\big(t+ \hat{\rho}_1\circ \mathfrak{n}_1(p)\big)\Big)= (t+ \hat{\rho}_1\big(\mathfrak{n}_1(p)\big), \hat{x})$, let $z_0= \gamma_{\mathfrak{P}_1(y_i), y_i}(t+ \hat{\rho}_1\big(\mathfrak{n}_1(p)\big))$, then 
\begin{align}
d(z_0, y_i)&= \hat{\rho}_1(z_0)- \hat{\rho}_1(y_i)= t+ \hat{\rho}_1\big(\mathfrak{n}_1(p)\big) \label{dist along geodesic} \\
d(z_0, p)&\leq d(p, y_i)+ d(y_i, z_0)< \frac{r}{320} .\nonumber 
\end{align}
Now using (\ref{diff n(x) and x}) and (\ref{dist along geodesic}), we have 
\begin{align}
d\big(\mathfrak{P}_1\circ \mathfrak{n}_1(z_0), \hat{x}\big)&= \sqrt{d\big(\mathfrak{G}_{\rho_1}\circ \mathfrak{n}_1(z_0), \mathfrak{G}_{\rho_1}(y_i)\big)^2- \big|\hat{\rho}_1\circ \mathfrak{n}_1(z_0)- \hat{\rho}_1(y_i)\big|^2} \nonumber \\
&\leq \sqrt{d(\mathfrak{n}_1(z_0), y_i)^2- \big|\hat{\rho}_1(z_0)- \hat{\rho}_1(y_i)\big|^2}+ \frac{n^{22}}{160}\eta^{\frac{1}{3n}}r+ \frac{\delta_1 r}{80} \nonumber \\
&\leq \sqrt{d(z_0, y_i)^2- |t+ \hat{\rho}_1\big(\mathfrak{n}_1(p)\big)|^2}+ n^{22}\eta^{\frac{1}{3n}}r= n^{22}\eta^{\frac{1}{3n}}r .\label{horizontal differ-1.2.0}
\end{align}
Then using (\ref{horizontal differ-1.2.0}), similar as (\ref{7.45.1}) we get $d(z_0, p)< \frac{r}{1280}+ n^{23}\eta^{\frac{1}{3n}}r$. The rest argument is similar to get (\ref{dense-1}), we can find $y_{j_2}\in B_{\frac{r}{1280}}(p)$ such that 
\begin{align}
d\big(\mathfrak{g}_1(y_{j_2}), (t, \hat{x})\big)\leq n^{25}\eta^{\frac{1}{3n}}r  .\label{dense-2}
\end{align}

From above, choose $\eta= n^{-100n}\epsilon^{3n}$, $\epsilon_1= n^{-301n^2}\epsilon^{9n^2}< 2^{-n}\eta^{3n}$, and finally $\delta= n^{-700n^4}\epsilon^{18n^4}< 2^{-100n^2}\epsilon_1^{2n^2}$, we get the pointed $(\frac{\epsilon r}{3})$-Gromov-Hausdorff approximation $\mathfrak{g}_1: B_{\frac{r}{1280}}(p)\rightarrow B_{\frac{r}{1280}}(0, \hat{p})$. 

\textbf{Step (4)}. Note $\mathbf{b}_1(p)= 0$, from (\ref{assumption 1 for almost split}), we have $\big|\rho_1(p)\big|< \epsilon_1 r$, then 
\begin{align}
\big|\hat{\rho}_1\circ \mathfrak{n}_1(p)+ t_0\big|&\leq \big|\hat{\rho}_1(p)+ t_0\big|+ \Big|d\big(p, \mathfrak{n}_1(p)\big)\Big|\leq \big|\rho_1(p)\big|+ \frac{\delta_1 r}{40}\leq \epsilon_1 r+ \frac{\delta_1 r}{40} \nonumber \\
&< \frac{1}{2}n^{-29}\epsilon r . \label{value at p is small}
\end{align}

Define $\mathcal{P}_1= \mathfrak{P}_1\circ \mathfrak{n}_1$, for $f_1= (\mathbf{b}_1, \mathcal{P}_1)$, then from (\ref{value at p is small}) and (\ref{assumption 1 for almost split}), 
\begin{align}
d\big(f_1(x), \mathfrak{g}_1(x)\big)&= \big|\mathbf{b}_1(x)- \hat{\rho}_1\circ \mathfrak{n}_1(x)+ \hat{\rho}_1\circ \mathfrak{n}_1(p)\big| \nonumber \\
&\leq \big|\mathbf{b}_1(x)- \rho_1(x)\big|+ \big|\rho_1(x)- \hat{\rho}_1\circ \mathfrak{n}_1(x)+ \hat{\rho}_1\circ \mathfrak{n}_1(p)\big| \nonumber \\
&\leq \epsilon_1 r+ \big|\hat{\rho}_1(x)+ t_0- \hat{\rho}_1\circ \mathfrak{n}_1(x)+ \hat{\rho}_1\circ \mathfrak{n}_1(p)\big| \nonumber \\
&\leq \epsilon_1r+ \big|\hat{\rho}_1(x)- \hat{\rho}_1\circ \mathfrak{n}_1(x)\big|+ \big|t_0+ \hat{\rho}_1\circ \mathfrak{n}_1(p)\big| \nonumber \\
&\leq \epsilon_1r+ d\big(x, \mathfrak{n}_1(x)\big)+ \frac{1}{2}n^{-29}\epsilon r \nonumber \\
&< n^{-29}\epsilon r .\label{f and g is close}
\end{align}

Now from (\ref{g-close}) and (\ref{f and g is close}), 
\begin{align}
\Big|d\big(f_1(x_1), f_1(x_2)\big)- d(x_1, x_2)\Big|&\leq \Big|d\big(\mathfrak{g}_1(x_1), \mathfrak{g}_1(x_2)\big)- d(x_1, x_2)\Big|+ 2\cdot n^{-29}\epsilon r \nonumber \\
&\leq n^{22}\eta^{\frac{1}{3n}}r+ n^{-29}\epsilon r \leq \frac{1}{3}\epsilon r .\label{dist close}
\end{align}

From (\ref{dense-1}), (\ref{dense-2}) and (\ref{f and g is close}), we have
\begin{align}
d\big(f_1(y_{i_2}), (t, \hat{x})\big)\leq n^{-29}\epsilon r+ n^{25}\eta^{\frac{1}{3n}}r\leq \frac{\epsilon}{3}r .\label{f-dense}
\end{align}

From (\ref{dist close}) and (\ref{f-dense}), we obtain that $f_1$ is an $(\frac{\epsilon}{3}r)$-Gromov-Hausdorff approximation.

\textbf{Step (5)}. We will only prove the case $k= 2$, the rest cases have similar argument. From the case $k= 1$, for $\mathbf{b}_i, i= 1, 2$, there is corresponding pointed $\big(\frac{\epsilon r}{3}\big)$-Gromov-Hausdorff approximation $f_{\mathbf{b}_i}: B_{\frac{r}{1280}}(p)\rightarrow B_{\frac{r}{1280}}(0, \hat{p}_i)\subset \mathbb{R}\times \mathbf{X}_{\mathbf{b}_i}$. We assume $f_{\mathbf{b}_i}= (\mathbf{b}_i, \mathcal{P}_{\mathbf{b}_i})$. We will prove the map $f_2(x)= \Big(\mathbf{b}_1(x), \mathbf{b}_2\big(x\big), \mathcal{P}_{\mathbf{b}_2}\circ \mathcal{P}_{\mathbf{b}_1}(x)\Big): B_{\frac{r}{1280}}(p)\rightarrow B_{\frac{r}{1280}}(0, \hat{p}_i)\subset \mathbb{R}^2\times \mathbf{X}_2$ is an pointed $(\frac{2\epsilon r}{3})$-Gromov-Hausdorff approximation.

From the first variation formula of arc length, for $x\in P_{\eta, \mathbf{b}_1, \mathbf{b}_2}^{r, \rho_1}\cap Q_{\eta, \mathbf{b}_1}^{r, \rho_1}\cap Q_{\eta, \mathbf{b}_2}^{r, \rho_2, \rho_1}$, let $\check{x}= \mathfrak{n}_1(x)$, then we have 
\begin{align}
&\quad \Big|\mathbf{b}_2(x)- \mathbf{b}_2\circ \mathcal{P}_{\mathbf{b}_1}(x) \Big|\leq \Big|\rho_2(x)- \rho_2\circ \mathfrak{P}_1(\check{x})\Big|+ 2\epsilon_1 r  \nonumber \\
&\leq \Big|\rho_2(\check{x})- \rho_2\circ \mathfrak{P}_1(\check{x})\Big|+ 2\epsilon_1 r+ \frac{\delta_1 r}{40} \nonumber \\
&= \Big|\hat{\rho}_2(\check{x})- \hat{\rho}_2\big(\sigma_{\check{x}, 1}(0)\big)\Big|+ 2\epsilon_1 r+ \frac{\delta_1 r}{40} \nonumber \\
&= \Big|\int_0^{\hat{\rho}_1(\check{x})- \eta r} \big\langle \nabla \hat{\rho}_1, \nabla \hat{\rho}_2\big\rangle\big(\sigma_{\check{x}, 1}(s)\big) ds\Big|+ 2\epsilon_1 r+ \frac{\delta_1 r}{40} \nonumber \\
&\leq  \int_0^{\hat{\rho}_1(\check{x})- \eta r} \Big(\big|\langle \nabla \mathbf{b}_1, \nabla \mathbf{b}_2\rangle \big|+ \big|\nabla(\hat{\rho}_1- \mathbf{b}_1)\big|+  \big|\nabla(\hat{\rho}_2- \mathbf{b}_2)\big|\Big)\big(\sigma_{\check{x}, 1}(s)\big)ds 
\nonumber \\
&\quad + 2\epsilon_1 r+ \frac{\delta_1 r}{40}   \nonumber \\
&\leq 3\eta r+ 2\epsilon_1 r+ \frac{\delta_1 r}{40}\leq \frac{1}{2n}\epsilon r .\nonumber 
\end{align}

Let $\mathfrak{w}(x)= \Big(\mathbf{b}_1(x), \mathbf{b}_2\big(\mathcal{P}_{\mathbf{b}_1}(x)\big), \mathcal{P}_{\mathbf{b}_2}\big(\mathcal{P}_{\mathbf{b}_1}(x)\big)\Big)$ for any $x, y\in B_{\frac{r}{1280}}(p)$, 
\begin{align}
\Big|d\big(\mathfrak{w}(x), \mathfrak{w}(y)\big)- d(x, y)\Big|&\leq \Big|d\big(\mathfrak{w}(x), \mathfrak{w}(y)\big)- d\big(f_{\mathbf{b}_1}(x), f_{\mathbf{b}_1}(y)\big)\Big|+ \frac{\epsilon r}{3} \nonumber \\
&\leq \frac{\Big|d\big(f_{\mathbf{b}_2}\circ\mathcal{P}_{\mathbf{b}_1}(x), f_{\mathbf{b}_2}\circ\mathcal{P}_{\mathbf{b}_1}(y)\big)^2- d\big(\mathcal{P}_{\mathbf{b}_1}(x), \mathcal{P}_{\mathbf{b}_1}(y)\big)^2\Big|}{d\big(\mathfrak{w}(x), \mathfrak{w}(y)\big)+ d\big(f_{\mathbf{b}_1}(x), f_{\mathbf{b}_1}(y)\big)}
+ \frac{\epsilon r}{3} \nonumber \\
&\leq \Big|d\big(f_{\mathbf{b}_2}\circ\mathcal{P}_{\mathbf{b}_1}(x), f_{\mathbf{b}_2}\circ\mathcal{P}_{\mathbf{b}_1}(y)\big)- d\big(\mathcal{P}_{\mathbf{b}_1}(x), \mathcal{P}_{\mathbf{b}_1}(y)\big)\Big|
+ \frac{\epsilon r}{3} \nonumber \\
&\leq \frac{2\epsilon r}{3} .\nonumber 
\end{align}

Now we have 
\begin{align}
&\quad \Big|d\big(f_2(x), f_2(y)\big)- d(x, y)\Big| \nonumber \\
&\leq \Big|d\big(f_2(x), f_2(y)\big)- d\big(\mathfrak{w}(x), \mathfrak{w}(y)\big)\Big|+ \Big|d\big(\mathfrak{w}(x), \mathfrak{w}(y)\big)- d(x, y)\Big| \nonumber \\
&\leq d\big(f_2(x), \mathfrak{w}(x)\big)+ d\big(f_2(y), \mathfrak{w}(y)\big)+ \frac{2\epsilon r}{3} \nonumber \\
&\leq \Big|\mathbf{b}_2(x)- \mathbf{b}_2\circ \mathcal{P}_{\mathbf{b}_1}(x) \Big|+ \Big|\mathbf{b}_2(y)- \mathbf{b}_2\circ \mathcal{P}_{\mathbf{b}_1}(y) \Big|+ \frac{2\epsilon r}{3} \nonumber \\
&\leq \frac{\epsilon r}{n}+ \frac{2\epsilon r}{3}\leq \epsilon r .\nonumber 
\end{align}

The rest argument is similar as in the case $k= 1$.
}
\qed

%$+ \sum_{i= 1}^2\big|\nabla (\mathbf{b}_i- \rho_i)\big|^2$

\section*{Part III. Covering groups of Riemannian manifolds with $Rc\geq 0$}

\section{Squeeze Lemma and Dimension induction on harmonic functions}

In this section, unless otherwise mentioned, we assume $\varphi: \tilde{M}^n\rightarrow M^n$ is the covering map with covering group $\Gamma$ such that $M^n= \tilde{M}^n/\Gamma$, where $(\tilde{M}^n, \tilde{g})$ and $(M^n, g)$ are two complete Riemannian manifolds and the metric $g$ is the quotient metric of $\tilde{g}$ with respect to group action of $\Gamma$.

We include the definition of quotient metric here for convenience.
\begin{definition}\label{def quotient metric}
{Consider a subgroup $G\subset \mathrm{Isom}(\mathbf{X})$, where $\mathbf{X}$ is a metric space, for every $\bar{x}, \bar{y}\in \mathbf{X}/G$, set the \textbf{quotient metric} $\bar{d}$ on $\mathbf{X}/G$ by:  
\begin{align}
\bar{d}(\bar{x}, \bar{y})= \inf\big\{d(x, y): x\in \bar{x}, y\in \bar{y}\big\}= \inf\big\{d(x, gy): g\in G\big\} .\nonumber 
\end{align}
}
\end{definition}

For any function $f$ defined on a domain $\mathbf{D}\subset M^n$, we define 
\begin{align}
\mathcal{L}(f)(\tilde{x})\vcentcolon = f\big(\varphi(\tilde{x})\big)\ , \quad \quad \quad \quad \forall\ \tilde{x}\in \varphi^{-1}(\mathbf{D}) . \nonumber 
\end{align}

For $\tilde{p}, \tilde{q}\in \tilde{M}^n$ and any $s> 0$, we define 
\begin{align}
\Gamma_{\tilde{p}}(s)= \big\{\gamma\in \Gamma|\ d(\tilde{p}, \gamma\tilde{p})\leq s\big\} .
\end{align}
Similarly we can define $\Gamma_{\tilde{q}}(s)$. When the point is fixed and clear in the context, we use $\Gamma(s)$ instead of $\Gamma_{\tilde{p}}(s)$ for simplicity.

The following Lemma is motivated by the use of the canonical fundamental domain $\mathcal{F}$ in \cite{Anderson}, and is needed for the proof of the squeeze lemma. 
\begin{lemma}\label{lem average integral on covering spaces}
{For any function $f\geq 0$ defined on $B_r(p)\subset M^n$, $\tilde{p}$ is one lift of $p$ and $B_r(\tilde{p})$ is the geodesic ball centered at $\tilde{p}$ with radius $r$ in $\tilde{M}^n$. Then 
\begin{align}
B_r(\tilde{p})\subset \varphi^{-1}\big(B_r(p)\big)\quad \quad and \quad \quad \fint_{B_r(\tilde{p})} \mathcal{L}(f)\leq 4^n\fint_{B_r(p)} f .\nonumber 
\end{align}
}
\end{lemma}

\pf
{Let $\tilde{\Omega}= \varphi^{-1}\big(B_r(p)\big)$. For any $\tilde{y}\in B_r(\tilde{p})$, 
\begin{align}
d\big(p, \varphi(\tilde{y})\big)= \inf_{h\in \Gamma}d(\tilde{p}, h\tilde{y})\leq d(\tilde{p}, \tilde{y})< r ,\nonumber 
\end{align}
which implies $\varphi\big(B_r(\tilde{p})\big)\subset B_r(p)$, then
\begin{align}
B_r(\tilde{p})\subset \varphi^{-1}\big(B_r(p)\big)= \tilde{\Omega} . \nonumber 
\end{align}

We choose a measurable section $s: B_r(p)\rightarrow B_r(\tilde{p})$, such that $\varphi(s(x))= x$ for any $x\in B_r(p)$. Let $T= s(B_r(p))$, then we have $\fint_{B_r(p)}f= \fint_T \mathcal{L}(f)$. Let $S$ be the union of $g(T)$ over all $g\in \Gamma$ such that $g(T)\cap B_r(\tilde{p})\neq \emptyset$. Then from $\mathrm{diam}(T)\leq 2r$, we obtain $S\subseteq B_{3r}(\tilde{p})$. And we also have $\fint_{B_r(p)}f= \fint_{S} \mathcal{L}(f)$ and $B_r(\tilde{p})\subseteq S$. 

From the Bishop-Gromov volume comparison theorem, we know $V(B_{3r}(\tilde{p}))\leq 3^nV(B_r(\tilde{p}))$. Note $B_r(\tilde{p})\subseteq S\subseteq B_{3r}(\tilde{p})$, we have $V(S)\leq 3^nV(B_r(\tilde{p}))$, hence 
\begin{align}
\fint_{B_r(\tilde{p})}\mathcal{L}(f)\leq \frac{1}{V(B_r(\tilde{p}))}\int_S \mathcal{L}(f)= \frac{V(S)}{V(B_r(\tilde{p}))}\fint_S \mathcal{L}(f)\leq 4^n\fint_{B_r(p)}f .\nonumber 
\end{align}
}
\qed

Before we state and prove our squeeze lemma, we would like to include the following well-known result and its proof here, because the squeeze lemma can be looked at the Gromov-Hausdorff perturbation version of the following Lemma.

\begin{lemma}\label{lem generator length for cocompact group action}
{For a locally compact, pointed length space $(\mathcal{Z}, q)$, assume $G$ is a closed subgroup of the isometry group $\mathrm{Isom}(\mathcal{Z})$ and $\mathcal{Z}/G= \mathbf{K}$, where $\mathbf{K}$ is a compact metric space with $\mathrm{diam}(\mathbf{K})= r_0$. Then for any $\epsilon> 0$, we have 
\begin{align}
G= \langle G(2r_0+ \epsilon)\rangle .\nonumber 
\end{align}
}
\end{lemma}

\pf
{For any $g_1\in G$, because $\mathcal{Z}$ is a length space, there is a segment $\gamma_{q, g_1q}$ from $q$ to $g_1q$. Then one can find a middle point $x\in \mathcal{Z}$ of the segment $\gamma_{q, g_1q}$ such that $d(x, q)= d(x, g_1q)= \frac{1}{2}d(q, g_1q)$. Since $\mathcal{Z}/G= \mathbf{K}$ and $G$ is a closed subgroup, one can find $g_2\in G$ such that $d(g_2q, x)\leq r_0$. Then we have 
\begin{align}
&d(q, g_2q)\leq d(q, x)+ d(x, g_2q)\leq 
\frac{1}{2}d(q, g_1q)+ r_0 \nonumber \\
&d(q, g_2^{-1}g_1q)= d(g_2q, g_1q)\leq d(g_2q, x)+ d(x, g_1q)\leq r_0+ \frac{1}{2}d(q, g_1q)  .\nonumber 
\end{align}

Assume $r= d(q, g_1q)$, then we have 
\begin{align}
g_2\in G\big(\frac{r}{2}+ r_0\big)\ , \quad \quad g_2^{-1}g_1\in G\big(\frac{r}{2}+ r_0 \big) ,\nonumber 
\end{align}
which implies $g_1\in \langle G\big(\frac{r}{2}+ r_0 \big)\rangle$. Using that $G$ is a closed subgroup, by induction on $r$ we get $g_1\in \big\langle G(2r_0)\big\rangle$, and the conclusion follows.
}
\qed

\begin{lemma}[Squeeze Lemma]\label{lem first squeeze big Eucliean directions}
{For any $\epsilon> 0, \delta= n^{-900n^4}\epsilon^{40n^4}$, integer $0\leq k\leq n$, if there exist harmonic functions $\big\{\mathbf{b}_i\big\}_{i= 1}^k$ defined on $B_r(p)$, satisfying $\mathbf{b}_i(p)= 0$ and 
\begin{align}
\sup_{B_{r}(p)\atop i= 1, \cdots ,k}|\nabla \mathbf{b}_i|\leq 2\ , \quad \quad  \quad \quad \sup_{t\leq r}\fint_{B_{t}(p)}\sum_{i, j= 1}^{k} \big|\langle \nabla \mathbf{b}_i, \nabla \mathbf{b}_j\rangle- \delta_{ij}\big|\leq \delta ,\label{average integral assumption}
\end{align}
then there is a family of $(\epsilon s)$-Gromov-Hausdorff approximation for $s\in (0, 10r_c]$,  
\begin{align}
f_s= (\mathbf{b}_1, \cdots, \mathbf{b}_k, \mathcal{P}_s): B_s(p)\rightarrow B_s(0, \hat{p}_s)\subset \mathbb{R}^k\times \mathbf{X}_{k, s} ,\nonumber 
\end{align}
where $r_c= \frac{r}{12800}$. And let $\mathrm{diam}\big(B_{r_c}(\hat{p}_{10r_c})\big)= r_0$, we have 
\begin{align}
\Gamma(r_c)\subset \big\langle \Gamma(\epsilon r_c+ 2r_0)\big\rangle  . \nonumber 
\end{align}
}
\end{lemma}

\pf
{In the proof, we can assume that $r_0\leq \frac{1}{2}r_c$, otherwise the conclusion follows directly. Let $\tilde{\mathbf{b}}_i= \mathcal{L}(\mathbf{b}_i), i= 1, \cdots, k$, from (\ref{average integral assumption}) and Lemma \ref{lem average integral on covering spaces}, we have $\mathbf{b}_i(p)= \tilde{\mathbf{b}}_i(\tilde{p})= 0$ and 
\begin{align}
\sup_{B_{r}(\tilde{p})\atop i= 1, \cdots ,k}|\nabla \tilde{\mathbf{b}}_i|\leq 2 \ , \quad \quad \quad \quad \sup_{t\leq r}\fint_{B_{t}(\tilde{p})}\sum_{i, j= 1}^{k} \big|\langle \nabla \tilde{\mathbf{b}}_i, \nabla \tilde{\mathbf{b}}_j\rangle- \delta_{ij}\big|\leq 4^n\delta .\nonumber 
\end{align}

For $0< \epsilon_1< 1$ to be determined later, from Theorem \ref{thm existence of harmonic map components imply GH-dist is small}, 
if we assume 
\begin{align}
\delta\leq 4^{-n}n^{-700n^4}\epsilon_1^{18n^4} ,\label{delta assumption squeeze}
\end{align}
there are two family of $(\epsilon_1 s)$-Gromov-Hausdorff approximation for $s\in (0, 10r_c]$,  
\begin{align}
f_s&= (\mathbf{b}_1, \cdots, \mathbf{b}_k, \mathcal{P}_s): B_s(p)\rightarrow B_s(0, \hat{p}_s)\subset \mathbb{R}^k\times \mathbf{X}_{k, s} \nonumber \\
\tilde{f}_s&= (\tilde{\mathbf{b}}_1, \cdots, \tilde{\mathbf{b}}_k, \tilde{\mathcal{P}_s}): B_s(\tilde{p})\rightarrow B_s(0, \check{p}_s)\subset \mathbb{R}^k\times \tilde{\mathbf{X}}_{k, s} .\nonumber 
\end{align}

From Lemma \ref{lem one side G-H map imply G-H dist}, there is an $(30\epsilon_1r_c)$-Gromov-Hausdorff approximation
\begin{align}
\Phi: B_{10r_c}(0, \check{p}_{10r_c})\rightarrow B_{10r_c}(\tilde{p}) , \nonumber 
\end{align}
where $B_{10r_c}(0, \check{p}_{10r_c})\subset \mathbb{R}^k\times \tilde{\mathbf{X}}_{k, 10r_c}$.

For any $\gamma\in \Gamma(r_c)$, we have 
\begin{align}
d\big(\tilde{\mathcal{P}}_{10r_c}(\tilde{p}), \tilde{\mathcal{P}}_{10r_c}(\gamma\tilde{p})\big)&\leq d\big(\tilde{f}_{10r_c}(\tilde{p}), \tilde{f}_{10r_c}(\gamma\tilde{p})\big)< d(\tilde{p}, \gamma\tilde{p})+ 10\epsilon_1r_c \nonumber \\
&\leq r_c+ 10\epsilon_1r_c< 2r_c .\nonumber 
\end{align}
Also note $\check{p}_{10r_c}= \tilde{\mathcal{P}}_{10r_c}(\tilde{p})$, then we have $\tilde{\mathcal{P}}_{10r_c}(\tilde{p}), \tilde{\mathcal{P}}_{10r_c}(\gamma\tilde{p})\in B_{2r_c}(\check{p}_{10r_c})$.

Apply the argument of Lemma \ref{lem almost mid pt} to $\tilde{\mathcal{P}}_{10r_c}(\tilde{p})$ and $\tilde{\mathcal{P}}_{10r_c}(\gamma\tilde{p})$, we can obtain $\tilde{z}\in B_{4r_c}(\check{p}_{10r_c})$ such that 
\begin{align}
\Big|d\big(\tilde{z}, \tilde{\mathcal{P}}_{10r_c}(\tilde{p})\big)- \frac{1}{2}d\big(\tilde{\mathcal{P}}_{10r_c}(\tilde{p}), \tilde{\mathcal{P}}_{10r_c}(\gamma\tilde{p})\big)\Big|\leq 10\sqrt{30\epsilon_1}r_c \label{almost mid ineq-1} \\
\Big|d\big(\tilde{z}, \tilde{\mathcal{P}}_{10r_c}(\gamma \tilde{p})\big)- \frac{1}{2}d\big(\tilde{\mathcal{P}}_{10r_c}(\tilde{p}), \tilde{\mathcal{P}}_{10r_c}(\gamma\tilde{p})\big)\Big|\leq 10\sqrt{30\epsilon_1}r_c . \label{almost mid ineq-2}
\end{align}
Then note $\tilde{\mathbf{b}}_i(\tilde{p})= \tilde{\mathbf{b}}_i(\gamma\tilde{p})= 0$, we have 
\begin{align}
d(\tilde{z}, \check{p}_{10r_c})&= d\big(\tilde{z}, \tilde{\mathcal{P}}_{10r_c}(\tilde{p})\big)\leq \frac{1}{2}d\big(\tilde{\mathcal{P}}_{10r_c}(\tilde{p}), \tilde{\mathcal{P}}_{10r_c}(\gamma\tilde{p})\big)+ 10\sqrt{30\epsilon_1}r_c \nonumber \\
&= \frac{1}{2}d\big(\tilde{f}_{10r_c}(\tilde{p}), \tilde{f}_{10r_c}(\gamma\tilde{p})\big)+ 10\sqrt{30\epsilon_1}r_c \nonumber \\
&< \frac{1}{2}d(\tilde{p}, \gamma\tilde{p})+ 70\sqrt{\epsilon_1}r_c \label{need 6.13.1} \\
&\leq \frac{1}{2} r_c+ 70\sqrt{\epsilon_1}r_c .\label{radius upper bound} 
\end{align}

Similarly, we have 
\begin{align}
d\big(\tilde{z}, \tilde{\mathcal{P}}_{10r_c}(\gamma\tilde{p})\big)< \frac{1}{2}d(\tilde{p}, \gamma\tilde{p})+ 70\sqrt{\epsilon_1}r_c  . \label{need 6.13.2}
\end{align}

From (\ref{radius upper bound}), we get that $\big(0, \tilde{z}\big)\in B_{\frac{1}{2} r_c+ 70\sqrt{\epsilon_1}r_c}(0, \check{p}_{10r_c})$. Because $\tilde{f}_{10r_c}$ is an $(10\epsilon_1 r_c)$-Gromov-Hausdorff approximation, there is $z_1\in B_{10r_c}(\tilde{p})$ such that 
\begin{align}
d\Big(\tilde{f}_{10r_c}(z_1), \big(0, \tilde{z}\big)\Big)< 10\epsilon_1 r_c . \label{lg 1.5}
\end{align}

Now we have  
\begin{align}
d(z_1, \tilde{p})&\leq d\big(\tilde{f}_{10r_c}(z_1), \tilde{f}_{10r_c}(\tilde{p})\big)+ 10\epsilon_1 r_c \nonumber \\
&\leq d\big(\tilde{f}_{10r_c}(z_1), (0, \tilde{z})\big)+ d\big((0, \tilde{z}), \tilde{f}_{10r_c}(\tilde{p})\big)+ 10\epsilon_1 r_c \nonumber \\
&< \frac{1}{2} r_c+ 70\sqrt{\epsilon_1}r_c+ 20\epsilon_1 r_c \nonumber  \\
&\leq \frac{1}{2} r_c+ 90\sqrt{\epsilon_1}r_c .\label{z_1 is in the small ball}
\end{align}

Assume $z_0= \varphi(z_1)$, we have $d(z_0, p)\leq d(z_1, \tilde{p})< \frac{1}{2}r_c+ 90\sqrt{\epsilon_1}r_c$. Then 
\begin{align}
d\big(\mathcal{P}_{10r_c}(z_0), \hat{p}_{10r_c}\big)&= d\big(\mathcal{P}_{10r_c}(z_0), \mathcal{P}_{10r_c}(p)\big)\leq d\big(f_{10r_c}(z_0), f_{10r_c}(p)\big) \nonumber \\
&\leq d(z_0, p)+ \epsilon_1\cdot (10r_c)< \frac{1}{2}r_c+ 100\sqrt{\epsilon_1}r_c< r_c ,\nonumber 
\end{align}
which implies $\mathcal{P}_{10r_c}(z_0)\in B_{r_c}(\hat{p}_{10r_c})$. From the assumption $\mathrm{diam}\big(B_{r_c}(\hat{p}_{10r_c})\big)= r_0$, 
\begin{align}
d\big(\mathcal{P}_{10r_c}(z_0), \hat{p}_{10r_c}\big)\leq r_0 .\label{use the diam}
\end{align}

Let $\tilde{\mathbf{b}}= (\tilde{\mathbf{b}}_1, \cdots, \tilde{\mathbf{b}}_k)$ and $\mathbf{b}= (\mathbf{b}_1, \cdots, \mathbf{b}_k)$, from (\ref{use the diam}) and (\ref{lg 1.5}), we obtain
\begin{align}
d(z_0, p)&\leq 10\epsilon_1 r_c+ d\big(f_{10r_c}(z_0), f_{10r_c}(p)\big)
\nonumber \\
&\leq 10\epsilon_1 r_c+ d\big(\mathbf{b}(z_0), 0\big)+ d\big(\mathcal{P}_{10r_c}(z_0), \mathcal{P}_{10r_c}(p)\big) \nonumber \\
&\leq 10\epsilon_1 r_c+ d\big(\tilde{\mathbf{b}}(z_1), 0\big)+ d\big(\mathcal{P}_{10r_c}(z_0), \hat{p}_{10r_c}\big) \nonumber\\
&\leq 20\epsilon_1 r_c+ r_0 .\label{need 6.16.1} 
\end{align}

There exists $\gamma_2\in \Gamma$ such that $d(z_1, \gamma_2\tilde{p})= d(z_0, p)< r_c$, then from (\ref{z_1 is in the small ball}), 
\begin{align}
d(\tilde{p}, \gamma_2\tilde{p})\leq d(\tilde{p}, z_1)+ d(z_1, \gamma_2\tilde{p}) < 2r_c .\nonumber 
\end{align}
We have $\gamma_2\tilde{p}\in B_{2r_c}(\tilde{p})$, from (\ref{lg 1.5}) and (\ref{need 6.16.1}), we obtain
\begin{align}
d\Big(\big(0, \tilde{z}\big), \tilde{f}_{10r_c}(\gamma_2\tilde{p})\Big)&< 10\epsilon_1 r_c+ d\big(\tilde{f}_{10r_c}(z_1), \tilde{f}_{10r_c}(\gamma_2\tilde{p})\big)\leq 20\epsilon_1 r_c+  d(z_1, \gamma_2\tilde{p}) \nonumber \\
&= 20\epsilon_1 r_c+ d(z_0, p)\leq 40\epsilon_1 r_c+ r_0 .\label{need 6.17.1}
\end{align}

Now from (\ref{need 6.13.1}) and (\ref{need 6.17.1}), we obtain
\begin{align}
d(\tilde{p}, \gamma_2\tilde{p})&\leq d\big(\tilde{f}_{10r_c}(\tilde{p}), \tilde{f}_{10r_c}(\gamma_2\tilde{p})\big)+ 10\epsilon_1 r_c \nonumber \\
&\leq d\Big(\tilde{f}_{10r_c}(\tilde{p}), \big(0, \tilde{z}\big)\Big)+ d\Big(\big(0, \tilde{z}\big), \tilde{f}_{10r_c}(\gamma_2\tilde{p})\Big)+ 10\epsilon_1 r_c \nonumber \\
&< d(\tilde{z}, \check{p}_{10r_c})+ (40\epsilon_1 r_c+ r_0)+ 10\epsilon_1 r_c  \nonumber \\
&\leq \frac{1}{2}d(\tilde{p}, \gamma\tilde{p})+ (120\sqrt{\epsilon_1} r_c+ r_0) . \label{required 6.2.1}
\end{align}

Similarly from (\ref{need 6.13.2}) and (\ref{need 6.17.1}), we can get 
\begin{align}
d(\gamma_2\tilde{p}, \gamma\tilde{p})&= d\big(\tilde{f}_{10r_c}(\gamma_2\tilde{p}), \tilde{f}_{10r_c}(\gamma\tilde{p})\big)+ 10\epsilon_1 r_c \nonumber \\
&\leq d\Big(\tilde{f}_{10r_c}(\gamma\tilde{p}), \big(0, \tilde{z}\big)\Big)+ d\Big(\big(0, \tilde{z}\big), \tilde{f}_{10r_c}(\gamma_2\tilde{p})\Big)+ 10\epsilon_1 r_c \nonumber \\
&< d\big(\tilde{z}, \tilde{\mathcal{P}}_{10r_c}(\gamma\tilde{p})\big)+ (40\epsilon_1 r_c+ r_0)+ 10\epsilon_1 r_c  \nonumber \\
&\leq \frac{1}{2}d(\tilde{p}, \gamma\tilde{p})+ (120\sqrt{\epsilon_1} r_c+ r_0) .\label{required 6.2.2} 
\end{align}

Let $\epsilon_1= \frac{\epsilon^2}{90000}$, from (\ref{required 6.2.1}) and (\ref{required 6.2.2}), we get $\gamma\in \Big\langle \Gamma\big(\frac{1}{2}d(\tilde{p}, \gamma\tilde{p})+ r_0+ \frac{120}{300}\epsilon r_c\big)\Big\rangle$.

From (\ref{delta assumption squeeze}), we can choose $\delta= n^{-900n^4}\epsilon^{40n^4}$. Now by induction, we have $\gamma\in \big\langle \Gamma(\epsilon r_c+ 2r_0)\big\rangle$, which implies the conclusion.
}
\qed

\begin{lemma}\label{lem at most n almost o.n. vec}
{If $\epsilon< \frac{1}{10n^2}$, there do not exist $\{\nu_i\}_{i= 1}^{n+ 1}\subset \mathbb{R}^n$ such that
\begin{align}
\sup_{1\leq i< j\leq n+ 1 }\big|\langle \nu_i, \nu_j\rangle\big|\leq \epsilon   \quad \quad and \quad \quad \sup_{1\leq i\leq n+ 1}\big||\nu_i|- 1\big|\leq \epsilon  .\nonumber 
\end{align}
}
\end{lemma}

\pf
{To prove the conclusion, we only need to show that $\nu_1, \cdots, \nu_{n+ 1}$ are linearly independent. By contradiction, otherwise, without loss of generality, one can assume $\nu_{n+ 1}= \sum_{i= 1}^n a_i\nu_i$ and $|a_1|= \max_{i= 1,\cdots, n}|a_i|$, where $a_i\in \mathbb{R}$. Then from the assumption,
\begin{align}
\epsilon\geq \big|\langle \nu_{n+ 1}, \nu_1\rangle\big|\geq |a_1|\cdot |\nu_1|^2- \sum_{i= 2}^n |a_i|\cdot \big|\langle \nu_{i}, \nu_1 \rangle\big| \geq |a_1|\cdot \big[(1- \epsilon)- (n- 1)\epsilon\big] ,\nonumber 
\end{align}
which implies 
\begin{align}
|a_1|\leq \frac{\epsilon}{1- n\epsilon} .\label{a_1 upper bound}
\end{align}

On the other hand, from (\ref{a_1 upper bound}) and $\epsilon< \frac{1}{10n^2}$, we can get 
\begin{align}
|\nu_{n+ 1}|\leq |a_1|\cdot \sum_{i= 1}^n |\nu_i|\leq \frac{\epsilon}{1- n\epsilon}\cdot n(1+ \epsilon)< 1- \epsilon ,\nonumber 
\end{align}
which is the contradiction to the assumption.
}
\qed

\begin{prop}\label{prop existence of n-dim harmonic map imply close to Euclidean ball}
{For complete Riemannian manifold $(M^n, g)$ with $Rc\geq 0$, if there are harmonic functions $\big\{\mathbf{b}_i\big\}_{i= 1}^n$, defined on $B_r(p)$, satisfying $\mathbf{b}_i(p)= 0$ and 
\begin{align}
\sup_{B_{r}(p)\atop i= 1, \cdots ,n}|\nabla \mathbf{b}_i|\leq 2 \ , \quad \quad \quad \quad \sup_{s\leq r}\fint_{B_{s}(p)}\sum_{i, j= 1}^{n} \big|\langle \nabla \mathbf{b}_i, \nabla \mathbf{b}_j\rangle- \delta_{ij}\big|\leq n^{-80000n^7} ,\label{n-dim harmonic map assumption}
\end{align}
then $\Gamma_{\tilde{p}}(\frac{r}{12800})= \{e\}$.
}
\end{prop}

\pf
{Let $\delta= n^{-80000n^7}$ in the rest proof. From (\ref{n-dim harmonic map assumption}) and Lemma \ref{lem average integral on covering spaces}, let $\tilde{\mathbf{b}}_i= \mathcal{L}(\mathbf{b}_i)$, then $\tilde{\mathbf{b}}_i(\tilde{p})= \mathbf{b}_i(p)= 0$ and  
\begin{align}
\sup_{B_{r}(\tilde{p})\atop i= 1, \cdots ,n}|\nabla \tilde{\mathbf{b}}_i|\leq 2\ ,  \quad \quad \quad \quad \sup_{s\leq r}\fint_{B_{s}(\tilde{p})}\sum_{i, j= 1}^{n} \big|\langle \nabla \tilde{\mathbf{b}}_i, \nabla \tilde{\mathbf{b}}_j\rangle- \delta_{ij}\big|\leq 4^n\delta .\label{n-dim harmonic map upper}
\end{align}

From (\ref{n-dim harmonic map assumption}), (\ref{n-dim harmonic map upper}) and Theorem \ref{thm existence of harmonic map components imply GH-dist is small}, if $\delta_1> 0$ is to be determined later such that 
\begin{align}
\delta\leq 4^{-n}\cdot n^{-700n^4}\delta_1^{18n^4} ,\label{delta need 1.1}
\end{align}
then there are two family of $(\delta_1\cdot  s)$-Gromov-Hausdorff approximation for any $0< s\leq 10r_c$,  
\begin{align}
f_s&= (\mathbf{b}_1, \cdots, \mathbf{b}_n, \mathcal{P}_s): B_s(p)\rightarrow B_s(0, \hat{p}_s)\subset \mathbb{R}^n\times \mathbf{X}_{n, s} \nonumber \\
\tilde{f}_s&= (\tilde{\mathbf{b}}_1, \cdots, \tilde{\mathbf{b}}_n, \tilde{\mathcal{P}}_s): B_s(\tilde{p})\rightarrow B_s(0, \check{p}_s)\subset \mathbb{R}^n\times \tilde{\mathbf{X}}_{n, s} ,\nonumber 
\end{align}
where $r_c= \frac{1}{12800}r$.

For any $0< t\leq r_c$, let $\mathrm{diam}\big(B_t(\hat{p}_{10t})\big)= t_1$. If $t_1\geq \frac{1}{4}t$, note 
\begin{align}
f_{10t}= (\mathbf{b}_1, \cdots, \mathbf{b}_n, \mathcal{P}_{10t}): B_{10t}(p)\rightarrow B_{10t}(0, \hat{p}_{10t})\subset \mathbb{R}^n\times \mathbf{X}_{n, 10t} , \nonumber 
\end{align}
is an $(10\delta_1 t)$-Gromov-Hausdorff approximation. 

For $\epsilon> 0$ to be determined later, assume 
\begin{align}
\delta_1= n^{-3400n^3}\epsilon^{110n} .\label{delta_1 and epsilon equ}
\end{align}
Apply Theorem \ref{thm AG-triple imply one more splitting}, there exist harmonic functions $\big\{\mathbf{b}_i\big\}_{i= 1}^{n+ 1}$ defined on $B_{l_1}(q)\subset B_{10t}(p)$ such that 
\begin{align}
\sup_{B_{l_1}(q)\atop i= 1, \cdots ,n+ 1}|\nabla \mathbf{b}_i|\leq 2 \quad \quad and \quad \quad \sup_{l\leq l_1}\fint_{B_{l}(q)}\sum_{i, j= 1}^{n+ 1} \big|\langle \nabla \mathbf{b}_i, \nabla \mathbf{b}_j\rangle- \delta_{ij}\big|\leq \epsilon ,\nonumber 
\end{align}
where $l_1= n^{-320n^3}\epsilon^{10n}t$.

Let $l\rightarrow 0$ in the above, we get $\{\nu_i\}_{i= 1}^{n+ 1}$, where $\nu_i\in T_qM^n$ such that 
\begin{align}
\sup_{1\leq i< j\leq n+ 1}\big|\langle \nu_i, \nu_j\rangle\big|\leq \epsilon \quad \quad and \quad \quad \sup_{1\leq i\leq n+ 1}\big||\nu_i|- 1\big|\leq \epsilon .\nonumber 
\end{align}
Choose suitable $\epsilon= \frac{1}{20n^2}$, then (\ref{delta need 1.1}) holds because of (\ref{delta_1 and epsilon equ}) and the definition of $\delta$.

Note $\epsilon< \frac{1}{10n^2}$, from Lemma \ref{lem at most n almost o.n. vec}, it is impossible. Hence we have 
\begin{align}
\mathrm{diam}\big(B_t(\hat{p}_{10t})\big)< \frac{1}{4}t \ , \quad \quad \quad \quad \forall t\in (0, r_c] .\label{radius bound of compact part} 
\end{align}

Note the definition of $\delta$ implies 
\begin{align}
\delta\leq n^{-900n^4}4^{-40n^4} .\label{delta assump upper bound-1}
\end{align}
For any $0< s\leq r_c$, from (\ref{delta assump upper bound-1}), (\ref{n-dim harmonic map assumption}), (\ref{radius bound of compact part}) and Lemma \ref{lem first squeeze big Eucliean directions}, we have 
\begin{align}
\Gamma(s)\subset \big\langle \Gamma(4^{-1} s+ 2\cdot \frac{1}{4}s)\big\rangle\subset \big\langle \Gamma(\frac{3}{4}s)\big\rangle \ ,\quad \quad \quad \quad \forall \ 0< s\leq r_c .\label{first squeeze}
\end{align}

From (\ref{first squeeze}), by induction on $s$, for any positive integer $m$, we have
\begin{align}
\big\langle\Gamma_{\tilde{p}}(r_c)\big\rangle= \big\langle\Gamma_{\tilde{p}}(\frac{3}{4}r_c)\big\rangle= \cdots =\Big\langle\Gamma_{\tilde{p}}\Big(\big(\frac{3}{4}\big)^mr_c\Big)\Big\rangle . \nonumber 
\end{align}

Note the group action is discrete, hence for big enough $m$, we have $\Big\langle\Gamma_{\tilde{p}}\Big(\big(\frac{3}{4}\big)^mr_c\Big)\Big\rangle= \{e\}$, which implies $\big\langle\Gamma_{\tilde{p}}(r_c)\big\rangle= \{e\}$.
}
\qed

We define $\mathfrak{ng}\big\langle \Gamma(r) \big\rangle$ as the minimal number of generators in $\Gamma(r)$ needed to generated $\big\langle \Gamma(r)\big\rangle$.

\begin{prop}\label{prop local fund group counting at close different points}
{For $\tilde{p}, \tilde{q}\in (\tilde{M}^n,\tilde{g})$ with $Rc(\tilde{g})\geq 0$, assume $d(\tilde{p}, \tilde{q})\leq \delta< \frac{s}{2}$, then
\begin{align}
&\mathfrak{ng}\big\langle\Gamma_{\tilde{p}}(\lambda_2)\big\rangle \leq \Big(\frac{2\lambda_2+ \lambda_1}{\lambda_1}\Big)^n\mathfrak{ng}\big\langle\Gamma_{\tilde{p}}(\lambda_1)\big\rangle\ , \quad \quad \quad \quad if\  0< \lambda_1\leq \lambda_2\nonumber \\
&\mathfrak{ng}\big\langle\Gamma_{\tilde{p}}(s)\big\rangle\leq \Big(\frac{3s+ 2\delta}{s- 2 \delta}\Big)^n \cdot \mathfrak{ng}\big\langle\Gamma_{\tilde{q}}(s- 2\delta)\big\rangle .\nonumber 
\end{align}
}
\end{prop}

\pf
{Note $\Gamma_{\tilde{p}}(\lambda_2)$ is a finite set, then there are only finite number of subsets $\mathcal{B}_1, \cdots, \mathcal{B}_m\subseteq \Gamma_{\tilde{p}}(\lambda_2)$ satisfying 
\begin{align}
\{e\}\subseteq \mathcal{B}_i \quad \quad and \quad \quad \{h_1^{-1}h_2: h_1\neq h_2\in \mathcal{B}_i\}\in \Gamma_{\tilde{p}}(\lambda_2)- \Gamma_{\tilde{p}}(\lambda_1) \quad \forall \ 1\leq i\leq m .\nonumber  
\end{align}

Let $\mathcal{A}_{\lambda_1, \lambda_2}$ be one of the above subsets, such that there does not exist $\mathcal{B}_i$ with $\mathcal{A}_{\lambda_1, \lambda_2}\subsetneqq \mathcal{B}_i$. Then $\mathcal{A}_{\lambda_1, \lambda_2}\subset \Gamma_{\tilde{p}}(\lambda_2)$ satisfy: 
\begin{align}
d(\tilde{p}, h_1^{-1}h_2\tilde{p})\geq \lambda_1 \quad \quad \quad \quad if\ h_1\neq h_2, h_1, h_2\in \mathcal{A}_{\lambda_1, \lambda_2} . \label{need 6.5.1} 
\end{align}
And for any $\gamma\in \Gamma_{\tilde{p}}(\lambda_2)$, there is $h\in \mathcal{A}_{\lambda_1, \lambda_2}$ such that $h^{-1}\gamma\in \Gamma_{\tilde{p}}(\lambda_1)$. 

For any $h\in \mathcal{A}_{\lambda_1, \lambda_2}$, we get $h\tilde{p}\in B_{\lambda_2}(\tilde{p})$. From (\ref{need 6.5.1}),
\begin{align}
B_{\frac{\lambda_1}{2}}(g_i\tilde{p})\cap B_{\frac{\lambda_1}{2}}(g_j\tilde{p})= \emptyset \quad \quad \quad \quad if\ g_1\neq g_2, g_1, g_2\in \mathcal{A}_{\lambda_1, \lambda_2} .\nonumber 
\end{align}

Now let $\#(\mathcal{A}_{\lambda_1, \lambda_2})$ denote the number of the elements in $\mathcal{A}_{\lambda_1, \lambda_2}$, we have 
\begin{align}
V\big(B_{\lambda_2+ \frac{\lambda_1}{2}}(\tilde{p})\big)&\geq V\big(\bigcup_{g_i\in \mathcal{A}_{\lambda_1, \lambda_2}} B_{\frac{\lambda_1}{2}}(g_i\tilde{p})\big)= \sum_{g_i\in \mathcal{A}_{\lambda_1, \lambda_2}}V\big(B_{\frac{\lambda_1}{2}}(g_i\tilde{p})\big) \nonumber \\
&= \sum_{g_i\in \mathcal{A}_{\lambda_1, \lambda_2}}V\big(B_{\lambda_2+ \frac{\lambda_1}{2}}(\tilde{p})\big)\cdot \Big[\frac{V\big(B_{\frac{\lambda_1}{2}}(\tilde{p})\big)}{V\big(B_{\lambda_2+ \frac{\lambda_1}{2}}(\tilde{p})\big)}\Big] \nonumber \\
&\geq \#(\mathcal{A}_{\lambda_1, \lambda_2})\cdot \Big(\frac{\frac{\lambda_1}{2}}{\lambda_2+ \frac{\lambda_1}{2}}\Big)^n\cdot V\big(B_{\lambda_2+ \frac{\lambda_1}{2}}(\tilde{p})\big) .\nonumber 
\end{align}
In the last inequality the Bishop-Gromov Comparison Theorem is used. Then 
\begin{align}
\#(\mathcal{A}_{\lambda_1, \lambda_2})\leq \big(\frac{2\lambda_2+ \lambda_1}{\lambda_1}\big)^n .\label{group A counting}
\end{align}

It is easy to see $\Gamma_{\tilde{p}}(\lambda_2)\subset  \mathcal{A}_{\lambda_1, \lambda_2}\cdot \Gamma_{\tilde{p}}(\lambda_1)$, from (\ref{group A counting}), we have  
\begin{align}
\mathfrak{ng}\big\langle\Gamma_{\tilde{p}}(\lambda_2)\big\rangle \leq \#\big(\mathcal{A}_{\lambda_1, \lambda_2}\big)\cdot \mathfrak{ng}\big\langle\Gamma_{\tilde{p}}(\lambda_1)\big\rangle\leq 
\Big(\frac{2\lambda_2+ \lambda_1}{\lambda_1}\Big)^n\mathfrak{ng}\big\langle\Gamma_{\tilde{p}}(\lambda_1)\big\rangle . \nonumber 
\end{align}

For $\gamma\in \Gamma_{\tilde{p}}(s)$, we have 
\begin{align}
d(\gamma\tilde{q}, \tilde{q})\leq d(\gamma\tilde{q}, \gamma\tilde{p})+ d(\gamma\tilde{p}, \tilde{p})+ d(\tilde{p}, \tilde{q})\leq 2d(\tilde{p}, \tilde{q})+ d(\gamma\tilde{p}, \tilde{p})\leq s+ 2\delta ,\nonumber 
\end{align}
which implies $\Gamma_{\tilde{p}}(s)\subset \Gamma_{\tilde{q}}(s+ 2\delta)$. Then
\begin{align}
\mathfrak{ng}\big\langle\Gamma_{\tilde{p}}(s)\big\rangle &\leq \mathfrak{ng}\big\langle\Gamma_{\tilde{q}}(s+ 2\delta)\big\rangle\leq \Big[\frac{2(s+ 2\delta)+ (s- 2\delta)}{s- 2\delta}\Big]^n\cdot \mathfrak{ng}\big\langle\Gamma_{\tilde{q}}(s- 2\delta)\big\rangle \nonumber \\
&= \Big(\frac{3s+ 2\delta}{s- 2 \delta}\Big)^n \cdot \mathfrak{ng}\big\langle\Gamma_{\tilde{q}}(s- 2\delta)\big\rangle  .\nonumber 
\end{align}
}
\qed

\begin{prop}\label{prop splitting one dim by group generator}
{For any $\epsilon\in (0, 1)$ and $0\leq k< n$, there exists $\delta= (n^{-30}\epsilon)^{5000n^7}$, such that if there exist harmonic functions $\big\{\mathbf{b}_i\big\}_{i= 1}^k$, satisfying
\begin{align}
\sup_{B_r(q)\atop i= 1, \cdots ,k}|\nabla \mathbf{b}_i|\leq 2 \quad \quad and \quad \quad \sup_{s\leq r}\fint_{B_s(q)}\sum_{i, j= 1}^k \big|\langle \nabla \mathbf{b}_i, \nabla \mathbf{b}_j\rangle- \delta_{ij}\big|\leq \delta .\label{assumption 4.2.6.1} 
\end{align}
Then there are harmonic functions $\big\{\mathbf{b}_i\big\}_{i= 1}^{k+ 1}$ defined on $B_{r_1}(q_1)\subset B_r(q)$, such that 
\begin{align}
&\mathfrak{ng}\big\langle \Gamma_{\tilde{q}}(r)\big\rangle\leq n^{340n^4}\epsilon^{-10n^2}\cdot \mathfrak{ng}\big\langle \Gamma_{\tilde{q}_1}(r_1)\big\rangle \label{generator squeeze est} \\
&\sup_{B_{r_1}(q_1)\atop i= 1, \cdots ,k+ 1}|\nabla \mathbf{b}_i|\leq 2 \quad \quad and \quad \quad \sup_{t\leq r_1}\fint_{B_{t}(q_1)}\sum_{i, j= 1}^{k+ 1} \big|\langle \nabla \mathbf{b}_i, \nabla \mathbf{b}_j\rangle- \delta_{ij}\big|\leq \epsilon .\nonumber 
\end{align}
}
\end{prop}

\pf
{Set $r_c= \frac{1}{12800}r$, from Proposition \ref{prop local fund group counting at close different points}, 
\begin{align}
\mathfrak{ng}\big\langle\Gamma_{\tilde{q}}(r)\big\rangle\leq \big(\frac{2r+ r_c}{r_c}\big)^n\cdot \mathfrak{ng}\big\langle\Gamma_{\tilde{q}}(r_c)\big\rangle .\label{require 6.6.1}
\end{align}

Define $\hat{r}_1= \min \Big\{s\geq 0|\ \Gamma_{\tilde{q}}(r_c)\subset \big\langle \Gamma_{\tilde{q}}(s)\big\rangle \Big\}$, then $\hat{r}_1\leq r_c$ and 
\begin{align}
\Gamma_{\tilde{q}}(r_c) \subset \big\langle\Gamma_{\tilde{q}}(\hat{r}_1)\big\rangle .\label{group subset 4.2.6.0}
\end{align}

For $\delta_1> 0$ to be determined later, we choose $\delta\leq n^{-900n^4}\delta_1^{40n^4}$. From (\ref{assumption 4.2.6.1}) and apply Lemma \ref{lem first squeeze big Eucliean directions} on $B_{10\hat{r}_1}(q)$, we have a family of $\big(\delta_1\cdot s\big)$-Gromov-Hausdorff approximation for any $0< s\leq 10\hat{r}_1$, 
\begin{align}
f_s: B_s(q)\rightarrow B_s(0, \hat{q}_s)\subset \mathbb{R}^k\times \mathbf{X}_{k, s} .\label{any small G-H-appr}
\end{align}
And assume $\mathrm{diam}\big(B_{\hat{r}_1}(\hat{q}_{10\hat{r}_1})\big)= r_0$, we have 
\begin{align}
\Gamma_{\tilde{q}}(\hat{r}_1)\subset \big\langle\Gamma_{\tilde{q}}(\delta_1\hat{r}_1+ 2r_0)\big\rangle .\label{group subset 4.2.6.1}
\end{align}

From (\ref{group subset 4.2.6.0}) and (\ref{group subset 4.2.6.1}), we have $\Gamma_{\tilde{q}}(r_c)\subset \big\langle\Gamma_{\tilde{q}}(\delta_1\hat{r}_1+ 2r_0)\big\rangle$. Then by the definition of $\hat{r}_1$, we get $\hat{r}_1\leq \delta_1\hat{r}_1+ 2r_0$. If we assume $\delta_1\leq \frac{1}{32}$, it yields $r_0\geq \frac{1}{2}\hat{r}_1- \frac{1}{2}\delta_1 \hat{r}_1\geq \frac{1}{4}\hat{r}_1$. Note $\mathrm{diam}\big(B_{\hat{r}_1}(\hat{q}_{10\hat{r}_1})\big)= r_0$, then combining (\ref{any small G-H-appr}), choose $\delta_1= n^{-3400n^3}\epsilon^{110n}$. Apply Theorem \ref{thm AG-triple imply one more splitting}, we can find harmonic functions $\big\{\mathbf{b}_i\big\}_{i= 1}^{k+ 1}$ defined on $B_{r_1}(q_1)\subset B_{10\hat{r}_1}(q)$, such that 
\begin{align}
\sup_{B_{r_1}(q_1)\atop i= 1, \cdots, k+ 1} |\nabla \mathbf{b}_i|\leq 2 \quad \quad and \quad \quad \sup_{t\leq r_1}\fint_{B_{t}(q_1)}\sum_{i, j= 1}^{k+ 1} \big|\langle \nabla \mathbf{b}_i, \nabla \mathbf{b}_j\rangle- \delta_{ij}\big|\leq \epsilon ,\nonumber 
\end{align}
where $r_1= n^{-320n^3}\epsilon^{10n}\hat{r}_1> 0$.

From Proposition \ref{prop local fund group counting at close different points} and $q_1\in B_{10\hat{r}_1}(q)$, we have 
\begin{align}
\mathfrak{ng}\big\langle \Gamma_{\tilde{q}}(21\hat{r}_1)\big\rangle&\leq (83)^n\mathfrak{ng}\big\langle \Gamma_{\tilde{q}_1}(\hat{r}_1)\big\rangle\leq (83)^n\cdot \Big(\frac{2\hat{r}_1+ r_1}{r_1}\Big)^n\cdot \mathfrak{ng}\big\langle \Gamma_{\tilde{q}_1}(r_1)\big\rangle \nonumber \\
&\leq n^{330n^4}\epsilon^{-10n^2}\cdot \mathfrak{ng}\big\langle \Gamma_{\tilde{q}_1}(r_1)\big\rangle .\label{squeeze step 2}
\end{align}

From (\ref{require 6.6.1}), (\ref{group subset 4.2.6.0}) and (\ref{squeeze step 2}), we have 
\begin{align}
\mathfrak{ng}\big\langle\Gamma_{\tilde{q}}(r)\big\rangle\leq \big(\frac{2r+ r_c}{r_c}\big)^n \cdot \mathfrak{ng}\big\langle\Gamma_{\tilde{q}}(\hat{r}_1)\big\rangle\leq n^{340n^4}\epsilon^{-10n^2}  \mathfrak{ng}\big\langle \Gamma_{\tilde{q}_1}(r_1)\big\rangle .\nonumber 
\end{align}
}
\qed

\begin{theorem}\label{thm uniform est on generators of local funda group}
{Assume $\varphi: \tilde{M}^n\rightarrow M^n$ is the covering map with covering group $\Gamma$ such that $M^n= \tilde{M}^n/\Gamma$, where $(\tilde{M}^n, \tilde{g})$ and $(M^n, g)$ are two complete Riemannian manifolds and the metric $g$ is the quotient metric of $\tilde{g}$ with respect to group action of $\Gamma$, furthermore $Rc\geq 0$, then $\mathfrak{ng}\big\langle\Gamma_{\tilde{p}}(1)\big\rangle \leq n^{n^{20n}} $ for any $\tilde{p}\in \tilde{M}^n$.
}
\end{theorem}

\pf
{Apply Proposition \ref{prop splitting one dim by group generator}, firstly when $k= 0$, we get 
\begin{align}
&\mathfrak{ng}\big\langle \Gamma_{\tilde{p}}(1)\big\rangle\leq n^{340n^4}\epsilon_1^{-10n^2} \cdot \mathfrak{ng}\big\langle \Gamma_{\tilde{q}_1}(r_1)\big\rangle \nonumber  \\
&\sup_{B_{r_1}(q_1)\atop i= 1}|\nabla \mathbf{b}_1|\leq 2 \quad \quad and \quad \quad \sup_{t\leq r_1}\fint_{B_{t}(q_1)} \big||\nabla \mathbf{b}_1|^2- 1\big|\leq \epsilon_1 ,\nonumber 
\end{align}
where $\epsilon_1= (n^{-30}\epsilon_2)^{\tau}, \tau\vcentcolon = 5000n^7$ and $\epsilon_2$ is to be determined later.

By induction, apply Proposition \ref{prop splitting one dim by group generator}, for $0\leq k\leq (n- 1)$, we get 
\begin{align}
&\mathfrak{ng}\big\langle \Gamma_{\tilde{p}}(1)\big\rangle\leq n^{340n^4}\epsilon_{k+ 1}^{-10n^2}\cdot \mathfrak{ng}\big\langle \Gamma_{\tilde{q}_{k+ 1}}(r_{k+ 1})\big\rangle \label{original group in induction} \\
&\sup_{B_{r_{k+ 1}}(q_{k+ 1})\atop i= 1, \cdots ,k+ 1}|\nabla \mathbf{b}_i|\leq 2 \ ,  \quad  \quad \quad \sup_{t\leq r_{k+ 1}}\fint_{B_{t}(q_{k+ 1})}\sum_{i, j= 1}^{k+ 1} \big|\langle \nabla \mathbf{b}_i, \nabla \mathbf{b}_j\rangle- \delta_{ij}\big|\leq \epsilon_{k+ 1} .\label{k-dim harmonic map } 
\end{align}
Let $\epsilon_{j}= (n^{-30}\epsilon_{j+ 1})^{\tau}, j= 1, \cdots, n- 1$.

We let $\epsilon_n= n^{-200\tau}$, then from (\ref{original group in induction}) and (\ref{k-dim harmonic map }) for $k= (n- 1)$, apply Proposition \ref{prop existence of n-dim harmonic map imply close to Euclidean ball}, we have $\Gamma_{\tilde{q}_n}(\frac{r_n}{12800})= \{e\}$.

From the induction formula for $\epsilon_j$, we have 
\begin{align}
\epsilon_j= n^{-30\tau\sum_{i= 0}^{n- j- 1}\tau^{i}}(\epsilon_n)^{\tau^{n- j}}\ , \quad \quad \quad 1\leq j\leq (n-1) .\label{induction on epsilon}
\end{align}

Then from (\ref{original group in induction}), (\ref{induction on epsilon}) and Proposition \ref{prop local fund group counting at close different points} yields 
\begin{align}
\mathfrak{ng}\big\langle \Gamma_{\tilde{p}}(1)\big\rangle&\leq n^{340n^4}\epsilon_1^{-10n^2}\mathfrak{ng}\big\langle \Gamma_{\tilde{q}_1}(r_1)\big\rangle \nonumber \\
&\leq (n^{340n^4})^n\cdot (\epsilon_1\cdots \epsilon_n)^{-10n^2}\mathfrak{ng}\big\langle \Gamma_{\tilde{q}_n}(r_n)\big\rangle
\nonumber \\
&\leq n^{340n^5}\cdot n^{1200n^2\tau^{n- 1}}\cdot \epsilon_n^{-20n^2\tau^{n- 1}}\cdot \mathfrak{ng}\big\langle \Gamma_{\tilde{q}_n}(r_n)\big\rangle \nonumber \\
&\leq n^{5000n^2\tau^n}\cdot \Big(\frac{2r_n+ \frac{1}{12800}r_n}{\frac{1}{12800} r_n}\Big)\mathfrak{ng}\Big\langle \Gamma_{\tilde{q}_n}\big(\frac{1}{12800} r_n\big)\Big\rangle \nonumber \\
&\leq n^{n^{20n}} .\nonumber 
\end{align}
}
\qed

\begin{theorem}\label{thm f.g. group has uniform generator bound}
{Suppose $(M^n, g)$ is a complete Riemannian manifold with $Rc\geq 0$, then for any finitely generated subgroup $\Gamma$ of $\pi_1(M^n)$, we have $\mathfrak{ng}(\Gamma)\leq n^{n^{20n}}$.
}
\end{theorem}

\pf
{From \cite[Theorem $82.1$]{Munkres}, there exists a covering map $\varphi: N^n\rightarrow M^n$ such that $\varphi_*\big(\pi_1(N^n, \hat{p})\big)= \Gamma$, where $\varphi(\hat{p})= p$. Now from \cite[Theorem $54.6(a)$]{Munkres}, we get $\pi_1(N^n, \hat{p})\simeq \Gamma$. Hence there exists a complete Riemannian manifold $(N^n, \tilde{g})$ with $Rc(N^n)\geq 0$ and $\pi_1(N^n)= \Gamma$.

Assume $\Gamma= \langle \gamma_1, \cdots, \gamma_k\rangle$ and $\sup_{i= 1, \cdots, k}d(\gamma_i\hat{p}, \hat{p})\leq C_1$. Then by scaling the metric $\tilde{g}$ to $\hat{g}= C_1^{-1}\tilde{g}$, we get 
that $\Gamma= \big\langle\Gamma_{\hat{p}}(1)\big\rangle$, where $\hat{p}\in (N^n, \hat{g})$. From Theorem \ref{thm uniform est on generators of local funda group}, we have $\mathfrak{ng}(\Gamma)\leq \mathfrak{ng}\big\langle\Gamma_{\hat{p}}(1)\big\rangle\leq n^{n^{20n}}$.
}
\qed

\section*{Acknowledgments}
The author thank Aaron Naber for helpful suggestion, which greatly improves the quantitative estimates in the earlier version of this paper. We are grateful to Jiaping Wang for continuous encouragement, William P. Minicozzi II and Christina Sormani for comments. We are indebted to Vitali Kapovitch for several email reply about \cite{KW}, which help us to understand the results there better. %Last but not least, we particularly thank the former anonymous referee, who carefully read the paper and gave detailed helpful suggestions. 

\appendix
  \renewcommand{\appendixname}{Appendix~\Alph{section}}

\end{document}